
\voffset-.5in

\input amstex
\input amsppt.sty

\catcode`\@=11

 \def\AMSTeXfeatures{\Plainheads 
   \let\current@vert=\AMS@vert}

 \def\Plainheads{\sh@ftdiam=0.05em
   \getlabeldims
   \let\vshaftfill=\plnvsolidfill
   \let\hshaftfill=\plnhsolidfill
   \let\th@rhead=\plnrhead
   \let\th@lhead=\plnlhead
   \let\th@dnhead=\plndnhead
   \let\th@uphead=\plnuphead}
 
 \def\glet{\global\let}

 \def\LaTeXfeatures{\catcode`\@=11
   \ifx\@clnwd\undefined \nol@g
      \input ltxcode.tex \dol@g \fi
   \ltxheads \let\current@vert=\new@vert
   \providelto \catcode`\@=\active}

 \def\nol@g{\def\wlog{\edef\garbage}}
 \def\dol@g{\let\wlog=\wl@g} \let\wl@g=\wlog
 \nol@g 

 \newbox\ltobox
 \def\providelto{{\setbox\z@=
   \hbox{$\to$}\minharrlen=\wd\z@
   \global\setbox\ltobox=\hbox{$\activeat>>>$}}
   \def\lto{\mathrel{\copy\ltobox}}}

 \def\ltxheads{\sh@ftdiam=\@wholewidth
   \getlabeldims
   \let\vshaftfill= \ltxvsolidfill
   \let\hshaftfill=\ltxhsolidfill
   \let\th@rhead=\ltxrhead
   \let\th@lhead=\ltxlhead
   \let\th@dnhead=\ltxdnhead
   \let\th@uphead=\ltxuphead}
 {\catcode`\@=\active
   \gdef@#1{\csname #1\string@at\endcsname}
   \glet\activeat=@}
 \def\def@#1{\expandafter\def\csname #1@at\endcsname}

 \def@>#1>#2>{\@rrow R{#1}{#2}}
 \def@<#1<#2<{\@rrow L{#1}{#2}}
 \def@ V#1V#2V{\@rrow V{#1}{#2}}
 \def@ A#1A#2A{\@rrow A{#1}{#2}}
 \def@/#1/#2/#3/{\@rrow{#1}{#2}{#3}}
 \def@.{\ifodd\row\ifmmode\noharrow
     \else\leavevmode.\spacefactor3000 \fi
   \else\novarrow\fi}
 \def@={\ifodd\row\harrow\hequalfill{}{}%
   \else\varrow\vequalfill{}{}\fi}
 \def@:#1{\ifx=#1\harrow\deffill{}{}%
   \else\leavevmode\null:#1\fi}
 \def@|{\current@vert}
  \def\AMS@vert{\varrow\vequalfill{}{}}
  \def\new@vert#1|#2|{\ifodd\row
   \let\nextarrow\vertexvarrow
   \else\let\nextarrow\varrow\fi
   \nextarrow\vshaftfill{#1}{#2}}
 \def@-{\ifmmode\let\next\hl@ne
   \else\let\next\AMSatdash \fi \next}
  \def\hl@ne#1-#2-{\harrow\hshaftfill{#1}{#2}}
  \def\AMSatdash{\let\next\relax\leavevmode
    \def\next@{\ifx\next-%
      \def\next-{\futurelet\next\nextii@}%
     \else\def\next{\hbox{-}}\fi\next}%
    \def\nextii@{\ifx\next-\def\next-{\hbox{---}}%
      \else\def\next{\hbox{--}}\fi\next}%
    \futurelet\next\next@}
 \def@(#1){\tweenarrows{#1}}
 \def@[#1]{\setsp@n#1\relax\activeat}
 \def\fiberbox{\hbox{$\vcenter{\hr@le\hbox{\vr@le
   \kern1ex\vbox{\kern1.2ex}\vr@le}\hr@le}$}}
  \def\hr@le{\hrule height \sh@ftdiam}
  \def\vr@le{\vrule width \sh@ftdiam}
 \def@+#1+#2+#3+{\ifodd\row \harrow{#1}{#2}{#3}%
   \else \varrow{#1}{#2}{#3}\fi}


 \def\Dnarrfill{\vequalfill\Dnhe@d}
 \def\Uparrfill{\Uphe@d\vequalfill}
 
 \def\ontofill{\rtarrfill\kern-0.3em 
   \th@rhead\kern 0.3em} 

 \def\rtarrfill{\hshaftfill\th@rhead}
 \def\ltarrfill{\th@lhead\hshaftfill}
 \def\dnarrfill{\vshaftfill\th@dnhead}
 \def\uparrfill{\th@uphead\vshaftfill}
 \def\hequalfill{\plnhfill=}
 \def\deffill{:\plnhfill=}
 \def\plnvextfill#1{\setbox\z@
   \hbox{\the\textfont3 #1}%
   \dimen@=\dp\z@\advance\dimen@\ht\z@
   \copy\z@ \kern-\dimen@ 
   \cleaders\copy\z@ \vfill
   \kern-\dimen@ 
   \box\z@}
 \def\plnhfill#1{$\m@th\mkern-1.5mu\mathord#1\mkern-6mu
    \cleaders\hbox{$\mkern-2mu\mathord#1\mkern-2mu$}\hfill
    \mkern-6mu\mathord#1\mkern-1.5mu$}
 \def\vequalfill{\plnvextfill{\char'167}}
 \def\plnvsolidfill{\plnvextfill{\char'077}}
 \def\plnhsolidfill{\plnhfill-}
 \def\ltxhsolidfill{\leaders\hrule height\topofshaft depth\botofshaft
   \hfill}
 \def\ltxvsolidfill{\leaders\vrule width\sh@ftdiam\vfill}
 \def\hdashfill{\hd@sh\wd@sh
   \xleaders \hbox{\wd@sh\hd@sh\wd@sh}\hfill
   \wd@sh\hd@sh}
 \def\vdashfill{\vd@sh\wd@sh
   \xleaders \vbox{\wd@sh\vd@sh\wd@sh}\vfill
   \wd@sh\vd@sh}
 \def\dashed{\ifinmeasureCD\else
    \ifodd\row\option{\let\hshaftfill=\hdashfill}%
   \else\option{\let\vshaftfill=\vdashfill}\fi\fi}


 \newdimen\CDstrutht  \newdimen\CDstrutdp
   \CDstrutht=0.875\baselineskip
   \CDstrutdp=0.375\baselineskip
 \newdimen\CDstrutlen \CDstrutlen=\CDstrutht
   \advance\CDstrutlen by \CDstrutdp

 \def\CDstrut{\vrule
   height \ifnum\row=1 \z@\else\CDstrutht \fi
   depth \ifnum\row=\numrows \z@ \else\CDstrutdp \fi
   width\z@}

 \newdimen\CDarrsurr \CDarrsurr=0.375em
 \newdimen\CDdashlen
    \CDdashlen= 0.1875\baselineskip
 \newdimen\CDvarrlen \CDvarrlen=1.5\baselineskip
 \newdimen\minharrlen 
  \setbox\z@\hbox{$\longrightarrow$} \minharrlen=\wd\z@
 \newdimen\minCDharrlen \minCDharrlen=2.5em 
\newdimen \minc@lwd
\def\findminc@lwd{\minc@lwd=2\CDarrsurr
  \advance\minc@lwd\minCDharrlen}

 \newdimen\sh@ftdiam


 \newdimen\labelsurr \labelsurr=1.25 em

\newcount\sp@ncnt \sp@ncnt=\@ne
\newcount\sp@ncnt@ \sp@ncnt@=\@ne
\newdimen\@rrwd \newdimen\@rrdp


 \def\adjustbot#1{\option{\advance\@rrdp#1\relax}}
 
\def\pushvertex#1{\global\p@shlen#1\relax
   \global\let\maybepush=\dopush}


 \newdimen\p@shlen \p@shlen=\z@

 
 \let\maybepush=\relax
 \def\dopush{\ifinmeasureCD 
   \advance\locdimen by -\p@shlen 
   \else\advance \@rrwd by -\p@shlen \fi 
   \global\let\maybepush=\relax \global\p@shlen=\z@\relax}


 \def\span@ne{\global\sp@ncnt=\@ne\relax}
 \def\setsp@n#1#2{\global\sp@ncnt=#1\relax
   \ifx\relax#2\relax\else\global\sp@ncnt@=#2\relax\fi}

 \def\plnrhead{\llap{$\rightarrow\mkern-1.5mu$}}
 \def\plnlhead{\rlap{$\mkern-1.5mu\leftarrow$}}

 \def\clap#1{\hbox to \z@{\hss #1\hss}}

 \def\plndnhead{\hbox{\the\textfont3 \char'171}}
 \def\plnuphead{\hbox{\the\textfont3 \char'170}}
 \def\Dnhe@d{\hbox{\the\textfont3 \char'177}}
 \def\Uphe@d{\hbox{\the\textfont3 \char'176}}

 \def\ltxrhead{\raise\@xisheight
   \llap{\smash{\@linefnt\@getrarrow(1,0)}}}
 \def\ltxlhead{\raise\@xisheight
   \rlap{\@linefnt\@getlarrow(-1,0)}}
 \def\ltxuphead{\setbox\z@=\rlap{%
   \kern\@halfwidth\@linefnt\char'66}%
   \copy\z@\kern-\ht\z@}
 \def\ltxdnhead{\setbox\z@=\rlap{%
   \kern\@halfwidth\@linefnt\char'77}%
   \ht\z@=\z@\box\z@}

 \def\wd@sh{\kern0.5\CDdashlen}
 \def\hd@sh{\vrule height\topofshaft depth\botofshaft
    width\CDdashlen}
 \def\vd@sh{\hrule height\CDdashlen
   depth\z@ width\sh@ftdiam}

\def\xylist{14{3434}13{2414}12{1723}%
  23{1413}34{1153}11{0867}43{0707}%
  32{0580}21{0414}31{0291}41{0}}
\newcount\tgtcnt@
\def\find@xyargs{\dimen@=\@rrdp
  \advance\dimen@ by \CDstrutlen
  \tgtcnt@=\dimen@ \dimen@=\@rrwd 
  \divide\dimen@ by \@m 
  \divide \tgtcnt@ by \dimen@ 
  \expandafter\testxy\xylist\relax
  \unitlength=\@xarg\@rrdp
  \divide\unitlength by\@yarg\relax}
\def\testxy#1#2#3{\ifnum\tgtcnt@>#3
    \@xarg=#1\relax \@yarg=#2\relax
    \let\next=\ignorerest
  \else\let\next\testxy\fi\next}
\def\ignorerest#1\relax{\relax}

\let\scalefactor=\@ne
\def\SWarrow{\find@xyargs\vector
  (-\@xarg,-\@yarg)\scalefactor\hskip-\wd\@linechar}
\def\NWarrow{\find@xyargs\vector
  (-\@xarg,\@yarg)\scalefactor\hskip-\wd\@linechar}
\def\NEarrow{\find@xyargs\vector
  (\@xarg,\@yarg)\scalefactor}
\def\SEarrow{\find@xyargs\vector
  (\@xarg,-\@yarg)\scalefactor}
\def\rightupline{\find@xyargs\@linelen=\scalefactor
     \unitlength\@sline}
\def\rightdownline{\find@xyargs\@yarg=-\@yarg\relax
     \@linelen=\scalefactor\unitlength\@sline}

\def\Sim{\ifodd\row\setbox\z@=\hbox{$\sim$}\dimen@=\ht\z@
 \advance\dimen@ by -\@xisheight
  \vbox{\box\z@\kern-\@xisheight\kern\dimen@}%
  \else\hbox{$\wr$}\fi}

%
\def\harrow#1#2#3{\inmeasureCDtrue\findminarrwd
  {#2}{#3}{\sp@ncnt\minharrlen}\inmeasureCDfalse\span@ne
  \mathrel{\hbox{\options\hplace{#1}\ulabel{#2}\dlabel{#3}}}}

\def\noharrow{\harrow\hfill{}{}}
\def\vertexvarrow#1#2#3{\findarrdp \@rrwd=\z@ \setsp@n\@ne\@ne
  \vbox to \z@{\kern-1.2\CDstrutht
  \rlap{\options\vplace{#1}\llabel{#2}\rlabel{#3}}\vss}}

\newif\ifinmeasureCD
\def\measurelabel#1{\setbox\z@
  \hbox{$\scriptstyle#1\kern\labelsurr$}%
  \ifdim\wd\z@>\@rrwd \@rrwd=\wd\z@\fi}
\def\findminarrwd#1#2#3{\@rrwd=#3\relax
   \measurelabel{#1}\measurelabel{#2}}
\def\findCDarrwd#1#2{\@rrwd=\minCDharrlen
   \measurelabel{#1}\measurelabel{#2}%
  }

\newcount\row \row=\@ne \newcount\col \col=\@ne 
 \newcount\numrows 
\numrows=\@ne
 \newcount\numcols
\newcount\arrspan \newdimen\vrtxhalfwd  \newbox\tempbox

\def\DANABUG{\advance\col by \@ne
 \@rrwd=\minCDharrlen
  \advance\@rrwd by \vrtxhalfwd
  \advance\@rrwd by \CDarrsurr
  \ifnum\col>\numcols \numcols=\col
     \newlocdimen{col\the\col}\locdimen=\@rrwd 
  \else \ifdim\@rrwd>\c@l \c@l=\@rrwd\fi\fi}

\def\drop#1\\{
  \findvrtxhalfsum\DANABUG\advance\row by 2 \measureinit}

\def\measureinit{\col=\@ne \vrtxhalfwd=-\CDarrsurr\arrspan=\@ne\@rrwd=\z@
   \setbox\tempbox=\hbox\bgroup$}
\def\measure{
  \let\harrow\measureCDarrow
  \let\CDCR=\measureCR 
   \findminc@lwd 
  \inmeasureCDtrue
  \row=\@ne \numcols=\z@ \measureinit}

\def\endmeasure{\findvrtxhalfsum\DANABUG
  \numrows=\row 
  \inmeasureCDfalse}




\def\newlocdimen#1{\advance\dimenc@unt by \@ne
  \ifnum\dimenc@unt<\insc@unt
     \else\errmessage{No room for the CD}\fi
  \dimendef\locdimen=\dimenc@unt
  \expandafter\dimendef\csname#1\endcsname=\dimenc@unt}

 \def\r@wc@l{\csname row\the\row col\the\col\endcsname}
 \def\c@l{\csname col\the\col\endcsname}

 \def\findvrtxhalfsum{$\egroup
  \newlocdimen{row\the\row col\the\col}
  \locdimen=\vrtxhalfwd 
  \vrtxhalfwd=0.5\wd\tempbox 
  \advance\vrtxhalfwd by \CDarrsurr
  \advance\locdimen by \vrtxhalfwd 
  \advance\@rrwd by \locdimen 
  \maybepush
  \divide\@rrwd by \arrspan\relax
  \ifdim\@rrwd<\minc@lwd
    \ifnum\col>\@ne \@rrwd=\minc@lwd\fi \fi
  \loop 
    \ifnum\col>\numcols \numcols=\col
       \newlocdimen{col\the\col}
       \locdimen=\@rrwd 
    \else \ifdim\@rrwd>\c@l \c@l=\@rrwd\fi \fi
   \ifnum\arrspan>\@ne
      \advance\arrspan by -1 \advance\col by \@ne
  \repeat }

 \def\measureCDarrow#1#2#3{\findvrtxhalfsum
   \arrspan=\sp@ncnt\relax\global\sp@ncnt=1\relax
   \advance\col by \@ne
   \findCDarrwd{#2}{#3}%
   \setbox\tempbox=\hbox\bgroup$}

 \newcount\dr@tn \dr@tn=\z@
 \def\locate#1:#2{\ifinmeasureCD\else
   \count@=-#1
   \multiply\count@ by 2
   \advance\count@ by #2
   \dimen@=\count@\@rrwd
   \ifnum\dr@tn=\@ne\relax \else\dimen@=-\dimen@ \fi
   \dimen@i=\@rrdp
   \ifnum\dr@tn>\z@\advance\dimen@i by \CDstrutlen \fi
   \dimen@i=\count@\dimen@i
   \count@=#2 \multiply\count@ by 2
   \divide\dimen@ by \count@
   \divide\dimen@i by \count@
   \lift\dimen@i\nudge\dimen@\fi}

\def\betweenCDrows{\advance\row by \@ne \col=\@ne
\options}


\def\hbegin{\hbox\bgroup\kern\c@l \kern-\r@wc@l$}
\def\hend{$\glet\maybepush\relax \CDstrut\egroup}
\def\vbegin{\setbox\tempbox=\hbox\bgroup$}
\def\vend{$\egroup\ht\tempbox=\z@\dp\tempbox\CDvarrlen
  \box\tempbox}
\def\setCD{\let\harrow=\setCDarrow
  \let\CDCR=\setCR 
  \row=\@ne \col=\@ne \hbegin}
\let\endsetCD=\hend 

\def\findarrwd{\@rrwd=\z@ \count@=\col \advance\count@ by\sp@ncnt
  \loop\ifnum\count@>\col \advance\count@ by -1
      \advance\@rrwd by\csname col\the\count@\endcsname\repeat}
\def\setCDarrow#1#2#3{\kern\CDarrsurr\advance\col by \@ne
  \findarrwd \advance\@rrwd by -\r@wc@l  
  \@rrdp=\z@ 
  \maybepush
  \advance\col by -\@ne \advance\col by \sp@ncnt \span@ne
  \hbox to \@rrwd{\options
   \@rrwd=\scalefactor\@rrwd\hss
   \hplace{#1}\ulabel{#2}\dlabel{#3}\hss}%
   \kern\CDarrsurr}

\newdimen\labspacei 
\newdimen\labspaceii 

\newdimen\@xisheight
  \@xisheight=\the\fontdimen22\textfont2
\newdimen\labelskip
  \labelskip=\the\fontdimen10\textfont3 
\newdimen\topofshaft
\newdimen\botofshaft
\newdimen\botofulabel
\newdimen\topofdlabel
\def\getlabeldims{
  \topofshaft=0.5\sh@ftdiam
  \botofshaft=\topofshaft
  \advance\topofshaft by \@xisheight  
  \advance\botofshaft by -\@xisheight  
  \botofulabel=\topofshaft
  \advance\botofulabel by \labelskip
  \topofdlabel=\botofshaft
  \advance\topofdlabel by \labelskip}

\def\ulabel{\ifnum\row=\@ne\let\next\ulabeli
   \else\let\next\ulabellap\fi\next}
\def\ulabeli#1{\vbox{
  \clap{\kern-\@rrwd$\scriptstyle#1$}%
  \kern\botofulabel}\maybeoffset}
\def\ulabellap#1{\vbox to \z@{\vss
  \clap{\kern-\@rrwd$\scriptstyle#1$}%
  \kern\botofulabel}\maybeoffset}
\def\dlabel{\ifnum\row=\numrows\let\next\dlabeli
   \else\let\next\dlabellap\fi\next}
\def\dlabeli#1{\vtop{\kern\topofdlabel
  \clap{\kern-\@rrwd$\scriptstyle#1$}%
  }\maybeoffset}
\def\dlabellap#1{\vbox to \z@{\kern\topofdlabel
  \clap{\kern-\@rrwd$\scriptstyle#1$}%
  \vss}\maybeoffset}
\def\rlabel#1{\vbox to \z@{\vss
  \rlap{\kern\labelskip$\scriptstyle#1$}%
  \vss\kern-\@rrdp}\maybeoffset}
\def\llabel#1{\vbox to \z@{\vss
  \llap{$\scriptstyle#1$\kern\labelskip}%
  \vss\kern-\@rrdp}\maybeoffset}
\def\swlabel#1{\vtop{\kern0.5\@rrdp
  \llap{$\scriptstyle#1$\kern\labelskip\kern-0.5\@rrwd}
  }\maybeoffset}
\def\nwlabel#1{\vbox{
  \llap{$\scriptstyle#1$\kern\labelskip\kern-0.5\@rrwd}%
  \kern-0.5\@rrdp}\maybeoffset}
\def\selabel#1{\vtop{\kern0.5\@rrdp
  \rlap{\kern0.5\@rrwd\kern\labelskip$\scriptstyle#1$}%
  }\maybeoffset}
\def\nelabel#1{\vbox{
  \rlap{\kern0.5\@rrwd\kern\labelskip$\scriptstyle#1$}%
  \kern-0.5\@rrdp}\maybeoffset}
\def\cplace#1{\vbox to \z@{\vss
  \clap{$#1$\kern-\@rrwd}%
  \kern-\@rrdp\vss}\maybeoffset}
\def\hplace#1{\hbox to \@rrwd{#1}\maybeoffset}
\def\vplace#1{\clap{\vbox to \z@{#1\kern-\@rrdp}}\maybeoffset}

\newdimen\nudgeamount \nudgeamount=\z@
\newdimen\liftamount \liftamount=\z@
\let\maybeoffset\relax
\newbox\offsetbox \newdimen\lastheight
\def\dooffset{
  \setbox\offsetbox=\lastbox \lastheight=\ht\offsetbox 
  \setbox\offsetbox=\vbox{\kern-\liftamount\box\offsetbox}%
  \ht\offsetbox=\lastheight
  \kern\nudgeamount\box\offsetbox\kern-\nudgeamount
  \global\nudgeamount=\z@ \global\liftamount=\z@
  \glet\maybeoffset=\relax}
\def\nudge#1{\ifinmeasureCD\else
  \global\advance\nudgeamount#1\relax
  \global\let\maybeoffset\dooffset\fi}
\def\lift#1{\ifinmeasureCD\else
  \global\advance\liftamount#1\relax
  \global\let\maybeoffset\dooffset\fi}

\def\findarrdp{\@rrdp=\CDvarrlen
  \ifnum\sp@ncnt@>1
    \advance\@rrdp by \CDstrutlen
    \multiply\@rrdp by \sp@ncnt@
    \advance\@rrdp by -\CDstrutlen \fi
 }

\def\varrow#1#2#3{\ifnum\sp@ncnt>\@ne 
     \sp@ncnt@=\sp@ncnt\relax\fi
  \findarrdp \@rrwd=\z@ 
  \kern\c@l
   \hbox to \z@{\options
   \@rrdp=\scalefactor\@rrdp
    \hss\vplace{#1}\llabel{#2}\rlabel{#3}\hss}%
  \global\advance\col by \@ne \setsp@n\@ne\@ne
  }

\def\novarrow{\varrow\vfill{}{}}

\def\tweenarrows#1{\findarrwd \findarrdp \setsp@n\@ne\@ne
  \rlap{\options\cplace{#1}}}

\def\usarrow #1#2#3{\dr@tn=\@ne
  \findarrwd \findarrdp \setsp@n\@ne\@ne 
  \rlap{\options\cplace{#1}\nwlabel{#2}\selabel{#3}}%
  \dr@tn=\z@}
\def\dsarrow #1#2#3{\dr@tn=\tw@
  \findarrwd \findarrdp \setsp@n\@ne\@ne 
  \rlap{\options\cplace{#1}\swlabel{#2}\nelabel{#3}}%
  \dr@tn=\z@}
 \def\@rrow#1{\csname #1@rrow\endcsname}
 \def\R@rrow{\harrow \rtarrfill}
 \def\L@rrow{\harrow \ltarrfill}
 \def\V@rrow{\varrow \dnarrfill}
 \def\A@rrow{\varrow \uparrfill}
 \def\SE@rrow{\dsarrow \SEarrow}
 \def\NW@rrow{\dsarrow \NWarrow}
 \def\SW@rrow{\usarrow \SWarrow}
 \def\NE@rrow{\usarrow \NEarrow}
 \def\DS@rrow{\dsarrow \dnslope}
 \def\US@rrow{\usarrow \upslope}
 \def\upslope{\find@xyargs
       \@linelen=\unitlength\@sline}
 \def\dnslope{\find@xyargs\@yarg=-\@yarg\relax
       \@linelen=\unitlength\@sline}

\newtoks\optionlist 
\optionlist={}
\let\options\relax
\def\dooptions{\the\optionlist\global\optionlist={}%
  \glet\options=\relax}
\def\option#1{\ifinmeasureCD\else
  \glet\options=\dooptions
  \global\optionlist=\expandafter{\the\optionlist\relax#1}\fi}
\def\wider#1{\ifinmeasureCD\else
   \option{\advance\@rrwd by #1}\fi}
\def\deeper#1{\ifinmeasureCD\else
   \option{\advance\@rrdp by #1}\fi}


{\def\\{\global\let\sptoken= }\\ }

\def\CR{\futurelet\nexttok\testCR}
\def\testCR{\ifx\nexttok\sptoken
   \let\next\eatspaceCR\else\let\next\CDCR\fi\next}
\def\eatspaceCR#1 {\CR}
\def\measureCR{\ifx\nexttok\endmeasure\let\nextCR\relax
    \else\let\nextCR\drop\fi\nextCR}
\def\setCR{\ifodd\row
  \ifx\nexttok\endsetCD\else\hend\betweenCDrows\vbegin\fi
  \else\vend\betweenCDrows\hbegin\fi}

\countdef\dimenc@unt=11
\def\CD#1\endCD{
   \begingroup\let\\=\CR
  \m@th\offinterlineskip
   \measure#1\endmeasure\null\,\vcenter{\setCD#1\endsetCD}\,
   \endgroup
    }

\ifx\@clnwd\undefined \nol@g\else\catcode`\ =14\relax\fi
 \font\@linefnt=line10 
 \newcount\@tempcnta
 \newcount\@tempcntb
 \newdimen\@tempdima
 \newdimen\@tempdimb
 \newdimen\@wholewidth
 \newdimen\@halfwidth
   \@wholewidth\fontdimen8\@linefnt \@halfwidth .5\@wholewidth
 \newdimen\unitlength
 \newcount\@xarg
 \newcount\@yarg
 \newcount\@yyarg
 \newbox\@linechar
 \newdimen\@linelen
 \newdimen\@clnwd
 \newdimen\@clnht
 \newif\if@negarg
 
 \def\@whilenoop#1{}

 \def\@whiledim#1\do #2{\ifdim #1\relax#2\@iwhiledim{#1\relax#2}\fi}

 \def\@iwhiledim#1{\ifdim #1\let\@nextwhile=\@iwhiledim 
         \else\let\@nextwhile=\@whilenoop\fi\@nextwhile{#1}}

 \def\@sline{\ifnum\@xarg< 0 \@negargtrue \@xarg -\@xarg \@yyarg -\@yarg
   \else \@negargfalse \@yyarg \@yarg \fi
 \ifnum \@yyarg >0 \@tempcnta\@yyarg \else \@tempcnta -\@yyarg \fi
 \ifnum\@tempcnta>6 \@badlinearg\@tempcnta0 \fi
 \ifnum\@xarg>6 \@badlinearg\@xarg 1 \fi
 \setbox\@linechar\hbox{\@linefnt\@getlinechar(\@xarg,\@yyarg)}%
 \ifnum \@yarg >0 \let\@upordown\raise \@clnht\z@
    \else\let\@upordown\lower \@clnht \ht\@linechar\fi
 \@clnwd=\wd\@linechar
 \if@negarg \hskip -\wd\@linechar \def\@tempa{\hskip -2\wd\@linechar}\else
      \let\@tempa\relax \fi
 \@whiledim \@clnwd <\@linelen \do
   {\@upordown\@clnht\copy\@linechar
    \@tempa
    \advance\@clnht \ht\@linechar
    \advance\@clnwd \wd\@linechar}%
 \advance\@clnht -\ht\@linechar
 \advance\@clnwd -\wd\@linechar
 \@tempdima\@linelen\advance\@tempdima -\@clnwd
 \@tempdimb\@tempdima\advance\@tempdimb -\wd\@linechar
 \if@negarg \hskip -\@tempdimb \else \hskip \@tempdimb \fi
 \multiply\@tempdima \@m
 \@tempcnta \@tempdima \@tempdima \wd\@linechar \divide\@tempcnta \@tempdima
 \@tempdima \ht\@linechar \multiply\@tempdima \@tempcnta
 \divide\@tempdima \@m
 \advance\@clnht \@tempdima
 \ifdim \@linelen <\wd\@linechar
    \hskip \wd\@linechar
   \else\@upordown\@clnht\copy\@linechar\fi}
 
 \def\@getlinechar(#1,#2){\@tempcnta#1\relax\multiply\@tempcnta 8
 \advance\@tempcnta -9 \ifnum #2>0 \advance\@tempcnta #2\relax\else
 \advance\@tempcnta -#2\relax\advance\@tempcnta 64 \fi
 \char\@tempcnta}
 
 \def\vector(#1,#2)#3{\@xarg #1\relax \@yarg #2\relax
 \@tempcnta \ifnum\@xarg<0 -\@xarg\else\@xarg\fi
 \ifnum\@tempcnta<5\relax
 \@linelen=#3\unitlength
 \ifnum\@xarg =0 \@vvector 
   \else \ifnum\@yarg =0 \@hvector \else \@svector\fi
 \fi
 \else\@badlinearg\fi}
 
 \def\@svector{\@sline
 \@tempcnta\@yarg \ifnum\@tempcnta <0 \@tempcnta=-\@tempcnta\fi
 \ifnum\@tempcnta <5
   \hskip -\wd\@linechar
   \@upordown\@clnht \hbox{\@linefnt  \if@negarg 
   \@getlarrow(\@xarg,\@yyarg) \else \@getrarrow(\@xarg,\@yyarg) \fi}%
 \else\@badlinearg\fi}
 
 \def\@getlarrow(#1,#2){\ifnum #2 =\z@ \@tempcnta='33\else
 \@tempcnta=#1\relax\multiply\@tempcnta \sixt@@n \advance\@tempcnta
 -9 \@tempcntb=#2\relax\multiply\@tempcntb \tw@
 \ifnum \@tempcntb >0 \advance\@tempcnta \@tempcntb\relax
 \else\advance\@tempcnta -\@tempcntb\advance\@tempcnta 64
 \fi\fi\char\@tempcnta}
 
 \def\@getrarrow(#1,#2){\@tempcntb=#2\relax
 \ifnum\@tempcntb < 0 \@tempcntb=-\@tempcntb\relax\fi
 \ifcase \@tempcntb\relax \@tempcnta='55 \or 
 \ifnum #1<3 \@tempcnta=#1\relax\multiply\@tempcnta
 24 \advance\@tempcnta -6 \else \ifnum #1=3 \@tempcnta=49
 \else\@tempcnta=58 \fi\fi\or 
 \ifnum #1<3 \@tempcnta=#1\relax\multiply\@tempcnta
 24 \advance\@tempcnta -3 \else \@tempcnta=51\fi\or 
 \@tempcnta=#1\relax\multiply\@tempcnta
 \sixt@@n \advance\@tempcnta -\tw@ \else
 \@tempcnta=#1\relax\multiply\@tempcnta
 \sixt@@n \advance\@tempcnta 7 \fi\ifnum #2<0 \advance\@tempcnta 64 \fi
 \char\@tempcnta}
\catcode`\ =10

\dol@g 
\catcode`\@=\active
\LaTeXfeatures

\input definiti
\input mathchar

\magnification 1200
\baselineskip 18pt


\def\pbf{\par\bigpagebreak\flushpar}
\def\pmf{\par\medpagebreak\flushpar}

\def\widetilde{\mathaccent"0365 }

\def\vp{\varphi}
\def\ol{\overline}
\def\em{{_{\Phi \Theta}}}

\def\a{\alpha}

\def\Ker{\hbox{\rm Ker}}
\def\sub{\hbox{\rm Sub}}
\def\Log{\hbox{\rm Log}}
\def\Hal{\hbox{\rm Hal}}
\def\Val{\hbox{\rm Val}}
\def\Aut{\hbox{\rm Aut}}
\def\Hom{\hbox{\rm Hom}}
\def\Hom{\hbox{\rm Hom}}
\def\Alv{\hbox{\rm Alv}}
\def\Cl{\hbox{\rm Cl}}
\def\Ct{\hbox{\rm Ct}}
\def\Sub{\hbox{\rm Sub}}
\def\Know{\hbox{\rm Know}}
\def\knowf{\hbox{\rm Knowf}}
\def\know{\hbox{\rm Know}}
\def\Om{\Omega}
\def\ol{\overline}
\def\om{\omega}
\def\si{\sigma}
\def\Bool{\hbox{\rm Bool}}
\def\G{\Gamma}
\def\g{\gamma}
\define\Th{\Theta}
\define \dl{\delta}

\topmatter
\title
Algebraic geometry in First Order Logic
\endtitle

\author
B. Plotkin, Hebrew University, Jerusalem
\endauthor
\endtopmatter

\bigskip
\centerline{\smc Abstract:}

In every variety of algebras $\Theta$ we can consider its logic
and its algebraic geometry.  In the previous papers geometry in
equational logic, i.e., equational geometry has been studied. Here
we describe an extension  of this theory towards the First Order
Logic (FOL).  The algebraic sets in this geometry are determined
by arbitrary sets of FOL formulas. The principal motivation of
such generalization lies in the area of applications to knowledge
science.

In this paper the FOL formulae  are considered in the context of
algebraic logic.

With this aim we define special Halmos categories. These
categories in the algebraic geometry related to FOL play the same
role as the category of free algebras $\Theta^0$ play in the
equational algebraic geometry.

The paper consists of three parts. Section 1 is of introductory
character.  The first part (sections 2--4) contains background on
algebraic logic in the given variety of algebras $\Theta$. The
second part is devoted to algebraic geometry related to FOL
(sections 5--7). In the last part (sections 8--9) we consider
applications of the previous material to knowledge science.


\par\newpage\par\flushpar
\centerline{\smc Contents}

\bigskip

\pmf {\smc 1. \  Algebra and logic (preliminary remarks)}

1.1 Multi-sorted algebra

1.2 Logic

1.3 Geometrical aspect

\pmf
 {\smc 2. \  Algebraic logic in the given variety of algebras
$\Theta$}

2.1 The main idea

2.2 Halmos categories and multi-sorted Halmos algebras

2.3  The category and algebra $\Hal_\Theta(G)$


\pmf
{3. \ \smc The category and algebra of formulas $\Hal_\Theta (\Phi)$}

3.1 The definition of the algebra $\Hal_\Theta (\Phi)$. Compressed
formulas


3.2 The value of a formula \pmf {4. \ \smc Structure of the
algebra $\Hal_\Theta (\Phi)$}

4.1 Elementary formulas

4.2 Additional remarks

4.3  Some  proofs

\pmf {5. \ \smc Algebraic geometry in FOL}

5.1 Sets of formulas and algebraic sets. The Galois
correspondence.

5.2 The Galois correspondence and morphisms

5.3 Lattices  and topology

5.4  The categories $K_{\Phi\Theta}(f)$ and $C_{\Phi\Theta}(f)$

5.5 The categories $K_{\Phi\Theta}$ and $C_{\Phi\Theta}$

 \pmf {6. \ \smc The Galois--Krasner theory in algebraic logic and algebraic geometry}

6.1 The group $\Aut (\Hal_\Theta(G))$

6.2 Main results

6.3 Proofs

\pmf {7. \ \smc Geometrical properties of models}

7.1 Isomorphisms of categories

7.2 A remark on the category $K_{\Phi\Theta}(G)$

7.3 Geometrical properties

 \pmf {8. \ \smc
Applications to the knowledge science}

8.1 Introduction

8.2 The category of knowledge

8.3 The category $L_\Theta(\Phi)$

8.4 The functor $\Ct_f$

8.5 Functors and homomorphisms

8.6 Knowledge base

 \pmf {9. \ \smc Equivalence of knowledge bases}

9.1. Definition

9.2 The case of finite models

9.3 Main result. Additional remarks
\par\newpage\par\flushpar

\bigskip

\pbf \head 1. \ Algebra and logic
\endhead

\subhead{1.1 Multi-sorted algebra}
\endsubhead
Keeping in mind applications, throughout the paper the term
algebra means {\it multi-sorted}, i.e., not necessarily
one-sorted, algebra.  We fix a set of sorts $\Gamma$.  In the
considered varieties $\Theta$ this set is finite, but it need not
to be finite in general.  We meet infinite $\G$ in the next
section.

For every algebra $G\in \Theta$ we write
$$G = (G_i, i \in \G).
$$
The set of operations $\Omega$ is called the {\it signature of
algebras} in $\Theta$. Every symbol $\omega \in \Omega $ has a
type $\tau = \tau(\omega) = (i_1, \ldots, i_{n}; j), i , j\in \G$.
An operation of type $\tau$ is a mapping
$$
G_{i_1} \times \ldots \times G_{i_n} \to G_j.
$$
All operations of the signature $\Omega $ satisfy some set of
identities.  These identities define  the variety $\Theta$ of
$\G$-sorted $\Omega$-algebras.  Let us consider homomorphisms and
free algebras in $\Theta$ .  A homomorphism of algebras in
$\Theta$ has the form
$$
\mu=(\mu_i, i \in \G) \colon G = (G_i, i \in \G) \to G'= (G'_i, i
\in \G).
$$
Here $\mu_i \colon G_i \to G'_i$ are mappings of sets, coordinated
with operations in $\Omega$. A congruence $Ker \mu = (\Ker {\mu_
i}, i \in \G)$ is the kernel of a homomorphism $\mu$.

We consider multi-sorted sets $X = (X_i, i \in \G)$ and the
corresponding free in $\Theta$ algebras
$$
W = W(X) = (W_i, i \in \G).
$$

A set $X$ and a free algebra $W$ can be presented as the free
union of all $X_i$ and all $W_i$, respectively.

Every (multi-sorted) mapping $\mu: X \to G$ is extended up to a
homomorphism $\mu: W \to G$.  Denote the set of all such $\mu$ by
$\Hom (W, G)$.  If all $X_i$ are finite, we treat this set as an
affine space. Homomorphisms $\mu\colon W \to G$ are points of this
space.

For the given $G = (G_i, i \in \G)$ and $X=(X_i, i \in \G)$ we can
consider the set
$$
G^X = (G^{X_i}_i, i \in \G).
$$
It is the set of mappings
$$
\mu = (\mu_i, i \in \G)\colon X \to G.
$$
There is a  natural bijection $\Hom(W,G) \to G^X$. More
information about multi-sorted algebras can be found in [Pl1].

Now let us turn to the models.  Fix some set of symbols of
relations $\Phi$.  Every $\vp \in \Phi$ has its type $\tau = \tau
(\vp) = (i_1, \ldots, i_n)$.  A relation, corresponding to $\vp$,
is a subset in the Cartesian product $G_{i_1} \times\ldots\times
G_{i_n}$. Denote by $\Phi\Theta$ the class of models $(G, \Phi,
f)$, where $G \in \Theta$, and $f$ is a interpretation of the set
$\Phi$ in $G$. As for homomorphisms of models, they are
homomorphisms of the corresponding algebras which are coordinated
with relations.


\subhead{1.2 \ Logic}
\endsubhead

In the sequel for the sake of clearness of the exposition we
sometimes
 "forget" that algebras are multi-sorted in general and use
 one-sorted language.
We consider logic in the given variety $\Theta$. For every finite
$X$, there is  a logical signature
$$
L = L_X=\{\vee,\wedge,\lnot,\exists x, \; \; x \in X\},
$$
where $X$ is $\mathop\bigcup\limits_{i\in\G} X_i$ for a finite
$\G$. We consider the set (more precisely, the $L$-algebra) of
formulas $L\Phi W$ over the free algebra $W=W(X)$.  This algebra
is an $L$-algebra of formulas of FOL over the given $\Theta$,
$\Phi$, and $X$.

First we define the atomic formulas. They are equalities of the
form $w\equiv w',$ with $w, w'\in W$ of the same sort and the
formulas $\vp (w_1, \ldots, w_n)$, where $w_i \in W,$ and all
$w_i$ are positioned according to the type $\tau=\tau(\vp)$ of the
relations $\vp$ and to the sorts.  The set of all atomic formulas
we denote by $M = M_X$.  Define $L\Phi W$ to be the absolutely
free $L_X$-algebra over $M_X$.

Let us consider another example of an $L_X$-algebra.

Given $W=W(X)$ and $G \in \Theta$, denote by $\Bool (W,G)$ the
Boolean algebra $\Sub (\Hom(W,G))$ of all subsets in $\Hom(W,G)$.
Define  the action of quantifiers in $\Bool (W,G)$.  Let $A$ be a
subset in $\Hom(W,G)$ and $x \in X_i$ be a variable of the sort
$i$.  Then $\mu\colon W\to G$ belongs to the set $\exists x A$ if
there exists $\nu \colon W \to G$ in $A$ such that $\mu (y) = \nu
(y)$ for every $y \in X$ of the sort $j, j \neq i$, and for every
$y \in X_i$, $y \neq x$.  Thus we get an $L$-algebra $\Bool (W,
G)$.

Now let us define a mapping
$$
\Val^X_f\colon  M_X \to \Bool (W,G),
$$
where $f$ is a model (the subject of knowledge), which realizes
the set $\Phi$ in the given $G$. If $w\equiv w'$ is an equality of
the sort $i$, then we set:
$$
\mu: W\to G \in \Val^X_f (w\equiv w') = \Val^X (w\equiv w')
$$
if $\mu_i(w)=\mu_i (w') $ in $G$.  Here the point $\mu$ is a
solution of the equation $w \equiv w'$.  If the formula is of the
form $\vp (w_1, \ldots, w_n)$, then
$$\mu \in \Val^X_f(\vp(w_1,\ldots, w_n))
$$
if $\vp(\mu(w_1),\ldots, \mu(w_n))$ is valid in the model $(G,
\Phi, f)$.  Here $\mu(w_j) = \mu_{i_j} (w_j)$, $i_j$ is the sort
of $w_j$.  The mapping $\Val^X_f$ is uniquely extended up to the
$L$-homomorphism
$$
\Val^X_f \colon L\Phi W\to \Bool (W,G).
$$
Thus, for every formula $u \in L\Phi W$ we defined its value
$\Val_f (u)$ in the model $(G,\Phi, f)$, which is an element in
$\Bool (W,G)$.

Every formula $u \in L\Phi W$ can be viewed as an equation in the
given model.  Then a  point $\mu\colon W\to G$ is the solution of
the ``equation" $u$ if $\mu \in \Val_f (u)$.


\subhead{1.3 Geometrical Aspect}
\endsubhead

In the $L$-algebra of formulas $L\Phi W$, $W = W(X)$, we consider
its various subsets $T$.
 On the other hand, we consider subsets $A$  in the
affine space $\Hom (W, G)$, i.e., elements of the $L$-algebra
$\Bool (W, G)$. For each  model $(G, \Phi, f)$  and for these $T$
and $A$ we establish the following {\it Galois correspondence}
between sets of formulas in  $L$-algebra of formulas $L\Phi W$ and
sets of points in the space $\Hom (W, G)$:
$$
\eqalign{ T^f &= A = \mathop{\bigcap}\limits_{u\in T} \Val_f
(u),\cr
 A^f&=T=\{u|A\subset \Val_f (u) \}.\cr}
$$

Here $A=T^f$ is a locus of all points satisfying the formulas from
$ T$. We regard $T$ also as a system of "equations", where each
"equation" is represented by a formula $u$ from $T$. Every set $A$
of such kind is said to be an {\it algebraic set}
 (or closed set, or  algebraic variety), determined for the given model.
We define {\it knowledge} as
$$
(X, T, A, (G, \Phi, f)).
$$
Here $T$ is a {\it description of knowledge} and $(G, \Phi, f)$ is
a {\it subject of  knowledge}.  $A=T^f$ is a {\it content of
knowledge}, represented as an algebraic variety, $X$ is a {\it
place of knowledge} (the place, where the knowledge is situated).
A set $A$ can be regarded also as a relation between elements of
$G$ derived from equalities and relations of the basic set $\Phi$.
The relation $A=T^f$ belongs to the multi-sorted set
$$
G^X = \{ G^{X_i}_i, \; \; i \in \Gamma\}.
$$
A set $T$ of the form $T=A^f$ for some $A$ is called an $f$-{\it
closed set}.
 For an arbitrary $T$ we have its closure
 $
 T^{ff}=(T^f)^f$ and for every $A \subset \Hom(W, G)$ we have the
 closure $A^{ff} = (A^f)^f$.

 It is easy to understand that the following rule takes place:

\medskip
 {\it A formula $v$ belongs to the set $T^{ff}$ if and only if the
 formula
 $$(\mathop{\wedge}\limits_{u \in T} u) \to v$$ holds in the model
 $(G, \Phi, f)$.}
\medskip

 If the set $T$ is infinite then the corresponding formula is called {\it infinitary}.

 We want to study knowledge with
 different, changing ``places of knowledge" $X$.  In this case
 one should consider different $W = W(X)$, different ``spaces of
 knowledge" $\Hom(W(X), G)$, and different $L\Phi W(X)$.

 Free in $\Theta$ algebras $W(X)$ with finite $X$ are the objects
 of the category, denoted by $\Theta^0$.   Morphisms of this
 category $s\colon W(X) \to W(Y)$ are arbitrary homomorphisms of
 algebras.  The category $\Theta^0$ is a full subcategory in the
 category $\Theta$.

 We intend to build a new category related to the first order logic
 for the given $\Theta$. This category  will play for the  FOL the role
 similar to that of
 the category of free algebras $\Theta^0$ for the equational logic.
 With this end we turn from pure logic to algebraic logic.
 Such a transition will allow us to associate
 description of knowledge with its content in a more interesting way.
The sets of the type $T=A^f$ also look more natural.


\head
 2. Algebraic logic
 \endhead

\subhead{2.1 The main idea}
\endsubhead
Algebraic logic deals with algebraic structures, related to
 various logical structures which correspond to different  logical calculi.
For example, Boolean algebras are associated with classical
propositional logic, Heyting algebras are associated
 with non-classical propositional logic, Tarski cylindric algebras
 and Halmos polyadic algebras are associated with FOL.

 Every logical calculus assumes that there are  formulas of the calculus,
 axioms of logic and rules of inference.  On this basis a
syntactical  equivalence of formulas compatible with their
semantical equivalence is defined.  The transition from pure logic
to algebraic logic is grounded on treating logical formulas up to
a certain equivalence. We call the corresponding classes the {\it
compressed formulas}. This transition leads to various special
algebraic structures, in particular to the structures mentioned
above.

 Logical calculi  are usually associated with some
 infinite set of variables.  Denote such a set by $X^0$.  In our
 situation it is a multi-sorted set $X^0 = (X^0_i, i \in \Gamma)$.
 Keeping in mind theory of knowledge and its geometrical aspect we
 will use a system of all finite subsets $X = (X_i, \; i \in \G)$
 of $X^0$ instead of this infinite universum.  This gives rise to
 multi-sorted logic and multi-sorted algebraic logic.  Every formula
 has a definite type (sort) $X$.  Denote the new set of sorts by
 $\G^0$.  It is a set of all finite subsets of the initial set
 $X^0$.

\subhead{2.2 Halmos Categories and multi-sorted Halmos algebras}
\endsubhead

 Fix some variety of algebras $\Theta$.   This means that a
 finite set of sorts $\G$, a signature $\Omega = \Omega (\Theta)$
 related to $\G$, and a system of identities $Id(\Theta)$ are
 given.

 Define Halmos categories for the given $\Theta$.

 First, for the given Boolean algebra $B$ we define its
 existential quantifiers [HMT].
 Existential quantifiers are the mappings $\exists\colon B \to B$
 with the conditions:

 1) \ $\exists 0 = 0$,

 2) \ $a < \exists a$,

 3) \ $\exists (a \wedge \exists b) = \exists a \wedge \exists b$,
 $0, a, b \in B$.

 The universal quantifier $\forall \colon B\to B$ is defined
 dually:

 1) \ $\forall 1  = 1$,

 2) \ $a > \forall a$,

 3) \ $\forall (a\vee \forall b) = \forall a \vee \forall b.$



 Let $B$ be a Boolean algebra and $X$ a set.  We say that $B$ is a
 {\it quantifier $X$-algebra} if a quantifier $\exists
 x \colon B \to B$ is defined  for every $x \in X$
 and for every two elements $x, y \in
 X$ the equality $\exists x \exists y = \exists y \exists x$ holds.

 One may consider also {\it quantifier $X$-algebras $B$ with equalities}
over $W(X)$. In such algebras, to
 each pair of elements $w, w' \in W(X)$ of the same sort it
 corresponds an element $w \equiv w' \in B$ satisfying the
 conditions

1)  $w\equiv w$ is the unit in $B,$

2)  $(w_1\equiv w'_1 \wedge\ldots \wedge w_n \equiv w'_n) < (w_1
\ldots w_n \omega \equiv w'_1\ldots w'_n \omega) $ where $\omega$
is an operation in $\Omega$ and everything is compatible with the
type of operation.

Now we will give the general definition of the Halmos category for
the given $\Theta$, which will be followed by examples.

{\it Halmos category} $H$ for an arbitrary finite $X=(X_i, i \in
\Gamma)$ fixes some quantifier $X$-algebra $H(X)$ with equalities
over $W(X)$. $H(X)$ are the  objects  of $H$.

The morphisms in $H$ correspond to morphisms in the category
$\Theta^0$. Every morphism $s_*$ in $H$ has the form
$$
s_* \colon H(X)\to H(Y),
$$
where  $s \colon W(X) \to W(Y) $ is a morphism in $\Theta^0$.

We assume that

1)  The transitions $W(X) \to H(X)$ and $s \to s_*$ yield a
(covariant) functor $\Theta^0 \to H$.

2) Every $s_* \colon H(X) \to H(Y)$ is a Boolean homomorphism.

3)  The coordination with the quantifiers is as follows:

$\qquad $ 3.1)  $s_1 \exists x a = s_2 \exists x a, \quad a \in
H(X)$, if $s_1 y = s_2y$ for every $y
 \in X, \; y \neq
x$.

$\qquad $ 3.2) $s\exists x a = \exists (sx) (sa) $ if $ sx = y \in
Y$ and $y = sx$ not in the support (see 4.2) of $sx'$, $x' \in X,
\; x' \neq x$.

4)  The following conditions describe coordination with equalities

$\qquad $ 4.1) $s_* (w\equiv w') = (sw \equiv sw')$ for $s\colon
W(X) \to W(Y)$, $w, w' \in W(X)$ are of the same sort.

$\qquad $ 4.2) $s^x_w a \wedge (w \equiv w') < s^x_{w'} a $ for an
arbitrary $a \in H(X), x \in X, w, w'$ of the same sort with $x$
 in $W(X)$, and $s^x_w\colon W(X) \to W(X)$ is defined by the rule:
$ s^x_w (x) = w, sy = y, y \in X,\; \;  y \neq x$.

This completes the definition of the Halmos category for a given
$\Theta$.

 Now we turn to
the definition of multi-sorted Halmos algebras.  Suppose that all
finite $X$ are subsets of some infinite universum $X^0$.

The set of sorts we denote by $\Gamma^0$.  It is the set of all
finite subsets in the fixed $X^0$.

For every $X \in \Gamma^0$ we consider the signature
$$
L_X = \{ \vee, \wedge, \rceil; \exists x, x \in X\}
$$
Let $L$ be the union of all $L_X, X \in \Gamma^0$, $V$ be the set
of all  equalities $w\equiv w',  w, w'$ of the same sort in
 $W(X).$ The equalities are considered as nullary operations. Let $S =
S_\Theta$ be the set of all $s: W(X) \to W(Y)$, $X, Y \in
\Gamma^0$. Denote  $\tilde L = L \cup S_\Theta\cup V$, (cf. 1.2).
The multi-sorted signature $\tilde L$ we consider as a new
signature for the FOL and we deal with the algebras in this
signature, which have the form $H = \{ H(X), X \in \Gamma^0\}$.
Here $H(X)$ is an algebra in the signature $L_X$, equipped by a
unary operations $s\colon H(X) \to H(Y)$ for $s: W(X) \to W(Y)$
and nullary operations from $V$. Let us take the subvariety of all
Halmos algebras in the variety of all $H$ and denote it by
$\Hal_\Theta$.

The identities of $\Hal_\Theta$ are the following:

1) \ The Boolean identities for every $H(X),$

2) \ The identities, making every $\exists x: H(X) \to H(X)$ an
existential quantifier,

3) \ The identities, making every $H(X)$ a quantorian $X$-algebra
with equalities,

 4) \ The identities, making every $s: H(X) \to H(Y)$ a Boolean
 homomorphism,

 5)  The identities, correlating the operations of the type $s$
 with quantifiers and equalities,

 6) \ The identities, making the system of all $H(X)$ and all $s$
 a category and $\Theta^0 \to H$ a functor.

 Algebras $H = \{ H(X), X \in \Gamma^0\}$ with the operations from
  the extended  signature $\tilde L$ and subject to the identities
  above are called  {\it multi-sorted Halmos algebras}.
 They constitute the variety $\Hal_\Theta$ of all multi-sorted Halmos
 algebras in $\Theta$.  It is evident that
if all $X$ are subsets of some universum $X^0$ then
 every Halmos category
 is a multi-sorted Halmos algebra, and every multi-sorted Halmos algebra can be treated
 as a Halmos category.  However, one should take into account
 that in this situation the notions  of subalgebra and subcategory  are
 different. The  same situation holds for homomorphisms.  Speaking about
 homomorphisms we mean the homomorphisms of Halmos algebras.

\subhead 2.3 The category and algebra $\Hal_\Theta (G)$
\endsubhead


We fix an algebra $G$ in the variety $\Theta$.  For an arbitrary
finite set $X = \{x_1, \cdots, x_n\}$ and the free in $\Theta$
algebra $W=W(x_1, \cdots, x_n)$ we consider the set of
homomorphisms $\Hom (W, G)$ as an affine space.  The points in
this space are homomorphisms $\mu: W \to G$.  There is a bijection
$$ \Hom (W, G) \to G^{(n)} $$ given by $$ \mu \to (\mu(x_1),
\cdots, \mu(x_n)). $$ Let now Bool $(W,G)$ be the Boolean algebra
of all subsets in $\Hom (W, G)$
$$
Bool (W,G) = Sub (\Hom (W,G)).
$$
We define quantifiers $\exists x: Bool (W, G) \to Bool (W, G), \;
x \in X$ in this algebra in the following way:
$\mu \in\exists x A \Leftrightarrow \exists \nu \in A$ such that
$\mu (y) = \nu(y)$ for $y \in X, \; y \neq x$.  Here $A$ is an
arbitrary set in $\Hom (W, G)$.  The axioms of the existential
quantifiers and equalities are fulfilled (see, for example, [Pl 1]).
We have
here a quantorian $X$-algebra with equalities $Bool (W, G)$ and we set
$$
\Hal_\Theta (G)(X) = Bool (W(X), G).$$


Consider now morphisms of the category $\Hal_\Theta(G)$.
 Take some $s:W(X) \to W(Y)$
in $\Theta^0$. We have
$$
\tilde s: \Hom (W(Y), G) \to \Hom (W(X), G),
$$
$\tilde s (\nu) = \nu s$ for arbitrary $\nu: W(Y) \to G$.  If $A$
is a subset in $\Hom (W(X), G)$, then $\nu \in s_* A = sA$ iff
$\tilde s (\nu) \in A$.  Hence, we get the mapping
$$
s_*: Bool (W(X), G) \to Bool (W(Y), G),
$$
which is a Boolean homomorphism.  It is easy to check that such
$s_*$ are correctly correlated  with quantifiers and equalities
([Pl1]). Thus, the Halmos category $\Hal_\Theta(G)$ is defined.
If, further, we confine ourselves by the sets $X$ in the given
$X^0$, then we come to the Halmos algebra $\Hal_\Theta (G)$.

\proclaim{Theorem 1} The algebras $\Hal_\Theta(G)$ with different
$G\in \Theta$ generate the variety of Halmos algebras
$\Hal_\Theta$.
\endproclaim
The proof of this theorem  will be given in  4.3.

Note that to every $s_*$ there corresponds the conjugate mapping
$$
s^*: Bool (W(Y). G) \to Bool (W(X), G),
$$
where for every $A \subset \Hom (W(Y),G)$ the set $s^* A$ is the
$\tilde s$-image of the set $A$.  From now on we assume that such
a conjugate mapping $s^*: H(Y) \to H(X)$ always exists in Halmos
algebras and categories.

\head 3. The category and the algebra of formulas
$\Hal_{\Phi\Theta}$
\endhead

\subhead{3.1 The definition of the algebra $\Hal_\Theta(\Phi)$.
Compressed formulas.}
\endsubhead

For some reasons we will use also the notation $\Hal_\Theta
(\Phi)$ for $\Hal_{\Phi\Theta }$. We start from the algebra of
pure formulas and will construct  $\Hal_{\Phi \Theta} =
\Hal_\Theta (\Phi)$ as an algebra of the special compressed
formulas.

 Fix the set $X^0$ and a set $\Phi$ of symbols of
relations.  This $\Phi$ can be also an empty set. For the given
$X$ the atomic formulas of the type (sort) $X$ are equalities $w
\equiv w'$,\break $ w, w' \in W(X)$ and formulas $\varphi(w_1,
\cdots, w_n)$ with $n$-ary $\varphi \in \Phi$ and $w_1, \cdots,
w_n \in W(X)$.  We denote the set of all such formulas by $M_X$.
We have also the multi-sorted set $M$ with $M(X) = M_X, X \in
\Gamma^0$.

Further we proceed from the signature $\tilde L$, defined earlier
for FOL.  Let $H$ be the multi-sorted $(\Gamma^0$-sorted)
absolutely free algebra in $\tilde L$ over the set of all atomic
formulas $M$. It is an algebra of pure formulas of FOL in
$\Theta$.  Every formula $u \in H$ has some special sort $X$.

Factorizing $H$ by the identities of the variety $\Hal_\Theta$ we
get the free in $\Hal_\Theta$ algebra $\tilde H$. Here $\tilde H$
is a free Halmos algebra over the same set of atomic formulas $M$.

Every formula in $H$ has a finite record in the signature $\tilde
L$ through atomic formulas.  The same is valid for the formulas in
$\tilde H$.

We consider relations for the atomic formulas of the form:
$$
s_* (\vp (w_1, \cdots, w_n)) = \vp (sw_1, \cdots, sw_n).
$$

The system of all such relations is denoted by $(*)$.

Now, by the definition, the algebra of compressed formulas
$\Hal_{\Phi\Theta} = \Hal_\Theta (\Phi)$ is the result of
factorization of the free Halmos algebra $\tilde H$ by the
defining relations $(*)$. Its elements are called {\it the
compressed formulas}.

Simultaneously, we defined also the Halmos category
$\Hal_\Theta(\Phi)$.

The algebra $\Hal_\Theta (\Phi)$ is no more a free Halmos algebra,
but it preserves some feature of freedom.

Let us consider (multi-sorted) mappings $\tau : M \to Q$, where
$Q$ is an arbitrary Halmos algebra in $\Theta$ and $M$ is the set
of atomic formulas.  For every $X \in \Gamma^0$ we have
$$
\tau_X : M_X \to Q(X).
$$

We call a mapping $\tau$ {\it correct} if for every $s:W(X)\to
W(Y)$:
$$
\eqalign{
&\tau (sw\equiv sw') = s_* (\tau(w\equiv w')),\cr
&\tau(\vp(sw_1, \cdots, sw_n)) = s_* \tau\vp(w_1, \cdots,
w_n),\cr}
$$
for $w, w', w_1, \cdots, w_n \in W(X).$

\proclaim{Theorem 2}  A mapping $\tau: M\to Q$ is extended up to a
homomorphism $\tau: \Hal_\Theta (\Phi) \to Q$ iff this mapping is
correct.
\endproclaim

Proof. Take a homomorphism $\tilde \tau: \tilde H \to Q$ for the
given $\tau : M \to Q$.  Besides, take the natural homomorphism
$\tilde H \to \Hal_\Theta (\Phi)$.  The problem is  to construct
the homomorphism $\tau': \Hal_\Theta (\Phi) \to Q$ with the
commutative diagrams
$$
\CD
M @>\tau >> Q\\
@. @/SE/// @AA\tau' A\\
@. \Hal_\Theta(\Phi)\\
\endCD
\qquad \CD
 \tilde H @>\tilde \tau >> Q\\
@. @/SE/// @AA\tau' A\\
@. \Hal_\Theta(\Phi)\\
\endCD
\qquad
$$


Such $\tau'$ exists if and only if the relations $(*)$  belong to
$\Ker \tau$.  But this exactly means that the mapping $\tau : M
\to Q$ is correct.

In the sequel we will relate the condition of correctness of the
mapping $\tau$  to the idea of the support of the elements of
Halmos algebra.

\subhead{3.2  The value of a formula}
\endsubhead

We fix  a model $(G, \Phi, f), G \in \Theta$, where $f$ is the
interpretation of the set $\Phi$ in $G$.  For every $n$-ary
relation $\vp \in \Phi$ we have $f(\vp) \subset G^{(n)}$.  For the
given model $f$ we define the canonical mapping $\Val_f$, which is
a homomorphism of Halmos algebras $$ \Val_f : \Hal_\Theta (\Phi)
\to \Hal_\Theta (G). $$ First we define it on atomic formulas.
Take some  $X \in \Gamma^0$  and consider formulas in $M_X$.  Now
define
$$
Val_f (w\equiv w') = \Val (w\equiv w')
$$
to be the equality in $\Hal_\Theta (G) (X)$. Thus,
$$
Val_f (w\equiv w')=\{\mu: W(X)\to G| w^\mu = {w'}^\mu\}.
$$
Analogously,
$$
Val_f (\vp(w_1, \cdots, w_n)) = \{
\mu|(w_1^\mu, \cdots, w^\mu_n) \in f(\vp) \}.
$$
 In particular,
it is clear, that the set $f(\vp)$ one-to-one corresponds to the
set
\newline $\Val_f(\vp(x_1, \cdots, x_n))$, all $x_1, \cdots, x_n$ are
different, $X = \{ x_1, \cdots, x_n\}$.

We have a multi-sorted mapping $$ \Val_f:M\to \Hal_\Theta (G). $$

The Halmos algebra $\tilde H$ is freely generated by the set $M$.
Thus, we have a homomorphism of Halmos algebras
$$
\Val_f:\tilde H \to \Hal_\Theta (G)
$$
Simultaneously, we have a homomorphism of $\tilde L$-algebras
$$
\Val_f:H\to \Hal_\Theta (G).
$$
It is easy to check that the relations $(*)$  hold in every
algebra $\Hal_\Theta(G)$, which means that they belong to the
kernels of the homomorphisms above.  This gives the homomorphism
$$
\Val_f: \Hal_\Theta (\Phi) \to \Hal_\Theta (G)
$$
that we do need.

For every formula $u$, pure or compressed, we have defined its
value on the model $( G, \Phi, f)$ as $\Val_f(u)$.  It is a subset
$A= \Val_f (u) $ in the space $\Hom (W(X), G)$.

The kernel $\Ker (\Val_f)$ which lies in $\Hal_\Theta(\Phi)$ is a
filter which  may be considered as the {\it elementary theory of
the model} $(G, \Phi, f)$.

Two formulas $u$ and $v$ of the same sort $X$ in the algebra of
pure formulas $H$ are called {\it semantically equivalent} if for
every model $(G, \Phi, f)$ the equality $\Val_f(u) = \Val_f (v)$
holds.

 \proclaim{Theorem 3}  The formulas $u$ and $v$ are
semantically equivalent iff they coincide in the algebra
$\Hal_\Theta (\Phi)$.
\endproclaim

Proof.  The proof will be  given in 4.3.

Now let us introduce the {\it logical kernel} of a homomorphism
$\mu: W(X) \to G$ by the rule  $\Log\Ker(\mu) = \{ u \in
\Hal_{\Phi\Theta} (X) | \mu \in \Val_f(u)\}$. It is a Boolean
filter, moreover, an ultrafilter in the algebra
$\Hal_{\Phi\Theta}(X)$.  The usual kernel $\Ker(\mu)$ can be
considered as the set of all equalities $w \equiv w'$ in
$\Log\Ker(\mu)$.

Every equality $w\equiv w'$  can be viewed as an equation in the
algebra $G$. The set of all solutions of such equation is the set
$\Val_f(w\equiv w') = \Val (w\equiv w')$.  Now an arbitrary
formula $u \in H$ can be also considered as an equation, but in
the model $(G, \Phi, f)$.  The set of all solutions of such
``equation" $u$ is the set $\Val_f(u)$.

Note here also that the point $\mu: W\to G$ is a solution of the
equation $w\equiv w'$ if and only if $(w\equiv w') \in \Ker
(\mu)$. The same $\mu$ is a solution of the ``equation" $u$ if and
only if $u \in \Log\Ker(\mu)$.  We use this point of view in the
universal algebraic geometry.

\head 4. Structure of the algebra $\Hal_\Theta(\Phi)$.
\endhead

\subhead{4.1 Elementary formulas}
\endsubhead

The operation $s \in S_\Theta$ can be included in the formulas $u
\in \Hal_{\Phi\Theta} (X)$. Formulas without such $s$ are
considered as elementary formulas. We prove here that every $u$ of
the a type $X$ is equivalent in some natural sense to an
elementary formula $v$ of the type $Y$, $X \subset Y$.

The notion of formulas equivalence generalizes the notion of
semantical equivalence.  The corresponding formulas do not need to
be of the same sort.

Let the set $X$ be a subset of the set $Y$.  The identical
inclusion $X \to Y$ defines the identical homomorphism $s^0 : W(X)
\to W(Y)$, with the corresponding $$ \eqalign{&s^0: H(X) \to
H(Y),\cr &s^0_* : \Hal_{\Phi\Theta} (X) \to \Hal_{\Phi\Theta}
(Y).\cr} $$ Denote by $\ol{u}$ a compressed formula in
$\Hal_{\Phi\Theta} (X)$, corresponding to a formula $u \in H$ of
the type $X$.  Let now $u$ and $v$ be two formulas of the types
$X_1$ and $X_2$, respectively. We call these formulas equivalent
if for some $Y$, containing $X_1$ and $X_2$, the equality $$
\ol{s^0_1 u} = \ol{s^0_2 v} $$ which gives $$ s^0_1\ol{u} =
s^0_2\ol{v} $$, holds. We call a formula $u \in H(X)$  an {\it
elementary} one if it can be expressed by atomic formulas of the
set $M_X$ in the terms of the signature $L_X$.

The following theorem is the theorem on the elimination of the
operations of the type $s \in S_\Theta$.

\proclaim{Theorem 4}  Every formula $u$ in the algebra of pure
formulas $H$ is equivalent to some elementary formula.
\endproclaim

Proof.   Denote by $H^0$ a subset in $H$, determined by the rule:
$u \in H^0(X)$ if the formula $u$ is equivalent to some elementary
formula.  It is clear that the set of all atomic formulas $M$ is a
part of the set $H^0$.  We want to check that $H=H^0$.  It is
enough to verify that $H^0$ is a subalgebra in $H$.

Let us make some remarks on the notion of equivalence of two
formulas.

Let $v$ be an elementary formula of the type $Y$ and the inclusion
$s^0: W(Y) \to W(Y')$ be given.  Then the formula $s^0v$ is
equivalent to the elementary formula of the type $Y'$.  Indeed,
let us take a formula $v'$ in $H(Y')$, whose record coincides with
$v$, but in the set $M_{Y'}$.  This $v'$ is an elementary formula
and it is easy to understand that $\ol{s_0v} = \ol{v'}$.  The last
means that $s^0v$ and $v'$ are equivalent.

Let now $u$ and $v$ be two equivalent formulas of the types $X_1$
and $X_2$ respectively and $Y$ a set, containing $X_1$ and $X_2$
with the inclusions
$$
s^0_1:W(X_1) \to W(Y), \; \; s^0_2: W(X_2) \to W(Y),
$$
such that $s^0_1 \bar u = s^0_2 \bar v$ holds true.  Take now an arbitrary set
$Y'$, containing $Y$, with the inclusion $s^0:W(Y) \to W(Y')$.
Then we have also inclusions
$$\eqalign{
&s^0s^0_1:W(X)_1) \to W(Y'),\cr
 &s^0s^0_2:W(X_2)\to W(Y'),\cr}
$$
and
$$
\eqalign{
&(s^0s^0_1)\ol{u} = s^0 (s^0_1\ol{u}) = s^0(s^0_2\ol{v}) =
(s^0s^0_2)\ol{v},\cr
&\ol{(s^0s^0_1)u} \, = \, \ol{(s^0s^0_2)v}.\cr}
$$
Thus, we see that along with the set $Y_1$ containing $X_1$ and
$X_2$ one can proceed from any $Y'$ containing $Y$ in the
definition of the equivalence of two formulas. Using these remarks
it is easy to check that every set $H^0(X)$ is closed under the
operations of the signature $L_X$. Now we need only to verify that
the set $H^0$ is invariant under the operations of the type $s$.

Take $s:W(X) \to W(Y)$ and let $u$ be a formula in $H^0(X)$.
Check that $s u \in H^0(Y)$.

The formula $u$ is equivalent to an elementary formula.  First,
consider the situation when $u$ is elementary itself.  Apply
induction by the record of the formula $u$ in atomic formulas of
the set $M_X$.  Here $X$ is a fixed set while $Y$ and $s:W(X) \to
W(Y)$ are arbitrary.

If $u$ is an atomic formula, then the formula $su$ is equivalent
to the atomic formula which is an elementary one.

Let now $u_1$ and $u_2$ of the type $X$ be elementary formulas and
$su_1$ and $su_2$ be equivalent to the elementary formulas $v_1$
and $v_2$.  Formulas $su_1$ and $su_2$ have the type $Y$.  Let
$Y_1$ be the type of formula  $v_1$ and $Y_2$ the type of formula
$v_2$. There is a set $Y'$ containing $Y, Y_1, Y_2$ with the
inclusions
$$
\eqalign{ &s^0:W(Y)\to W(Y'),\cr &s^0_1: W(Y_1) \to W(Y'),\cr
&s^0_2: W(Y_2) \to W(Y'),\cr}
$$
such that
$$
\eqalign{
&\ol{(s^0s)u_1} \, = \, \ol{s_1^0v_1},\cr
&\ol{(s^0s)u_2} \, = \, \ol{s^0_2 v_2}\cr}
$$
hold.

Take  $u = u_1 \vee u_2$.  We have
$$
\eqalign{
&
\ol{su} = s\ol{u} = s(\ol{u}_1 \vee \ol{u}_2) = s\ol{u}_1 \vee
s\ol{u_2};
\cr
&(s^0 s) \ol{u} = s^0 (s\ol{u}_1 \vee s\ol{u}_2) = (s^0s)\ol{u}_1 \vee
(s^0s) \ol{u}_2 =\cr
&=\ol{s^0_1v_1} \vee \ol{s_2^0v_2} = \ol{s^0_1v_1\vee s_2^0v_2}.\cr}
$$

We use the first remark on the equivalence of formulas.  We have:
$$
\eqalign{
&\ol{(s^0s)u} = \ol{s^0_1v_1} \vee \ol{s^0_2 v_2} = \ol{v'_1} \vee
\ol{v'_2} = \cr
&\ol{v'_1 \vee v'_2} = \ol{v'};\cr}
$$
The formula $v'=v'_1\vee v'_2$ is an elementary formula of the
type $Y'$,  and $su \in H^0(Y)$.

Similarly one can check that if $u = u_1\wedge u_2 $ or $u =
\neg u_1$, then $su\in H^0(Y)$.

Consider further the case $u= \exists xu_1$.
As earlier, $u$ and $u_1$ are of the type $X$ and are elementary,
and the assumption of the induction holds for $u_1$, i.e., $su_1$
is equivalent to an elementary formula for any $s$.

Let us turn to the formulas $su=s\exists xu_1$ and to $\ol{su} =
s\ol{u} = s\exists x\ol{u}_1$.

The variable $x$ belongs to the set $X$, and $su$ belongs to the
set $H(Y)$.  Extend the set $Y$ up to $Y'$,
adding an arbitrary variable $y' $ to $Y$.  Let $s^0: W(Y) \to
W(Y')$ be the corresponding identical inclusion.  Pass to $s^0s:
W(X) \to W(Y')$ and consider a homomorphism
$s_1\colon W(X) \to W(Y')$ acting as $s^0s$ on all elements from
$X$ except $x$, and $s_1 x = y'$.  We have
$$
(s^0s) \ol{u} = (s^0s) \exists x \ol{u_1} = s_1 \exists x
\ol{u}_1.
$$
The element $y'=s_1x$ does not belong to the record of any $s_1
x_1, x_1 \neq x$ by construction. Then
$$
(s^0 s) \ol{u} = \exists y' s_1 \ol{u}_1.
$$

By the assumption of induction, the element $s_1 u_1$ is
equivalent to the elementary $v$.  Let  $v$ be of the type $Y_1$,
$s_1 u_1$ has the type $Y'$.  Take the set $Y'_1$, containing
$Y_1$ and $Y'$, and let
$s^0_1\colon W(Y_1) \to W(Y'_1)$ and
$s^0_2\colon W(Y') \to W(Y'_1)$
be identical inclusions.  Then $(s^0_2 s_1) \ol{u}_1 = s^0_1
\ol{v}$.  We have also $(s^0_2 s^0s) \ol{u} = s^0_2 s_1 \exists
x\ol{u}_1 = s^0_2 \exists y' s_1 \ol{u}_1 = \exists y' s_2^0 s_1
\ol{u}_1 = \exists y' \ol{s^0_1 v} = \exists y' \ol{v^1}$, where
$v'$ is an elementary formula of the type $Y'_1$.  The same is
true for the formula $\exists y'v'$.  We have:
$$
\ol{(s_2^0s^0)(su)} = \ol{\exists y'v'}, su \in H^0 (Y).
$$
This equality holds for every elementary formula $u$ in $H^0(X)$.

Let us consider  the case when $u$ is not necessarily elementary.
The formula $u$ is equivalent to the elementary formula $u'$, for
example, of the type $X_1$.  Take a set $X'$ with the inclusion
$s^0: W(X) \to W(X')$ and $s^0_1$: $W(X_1) to W(X')$ and with the
condition $\ol{s^0u} = s^0\ol{u} =s^0_1 \ol{u'} = \ol{s^0_1u'}$.
Since $u'$ is an elementary formula, $\ol{s^0_1u'} = \ol{v}$,
where $v$ is an elementary formula of the type $X'$. Thus,
$\ol{s^0u} = \ol{v}$.

Select a set $Y'$ containing $Y$ with the commutative diagram
$$
\CD
W(X) @> s>> W(Y)\\
@V s^0VV@VVs_2^0V\\
W(X')@>s_1>> W(Y') \\
\endCD
$$
where $s^0, s^0_2$ are inclusions, and $s_1s^0=s^0_2s$.  Let us
apply it to $\ol{u}$:
$$
s_1 s^0 \ol{u} = s^0_2 s\ol{u} = \ol{s_1v}.
$$
 Here $v$
 is an elementary formula of the type $Y'$, hence the
 formula $s_1v$ is equivalent to the elementary formula $v_1$.
 Let $v_1 $ be of the type $Y_1$ and let the set $Y'_1$ contain
 $Y_1$ and $Y'$.  We have the inclusions $s^0_3 \colon W(Y') \to
 W(Y'_1)$,  $s^0_4\colon W(Y_1) \to W(Y'_1)$ with the condition
 $s^0_3s_1\ol{v} = s^0_4\ol{v_1} = \ol{s_4^0 v_1} = \ol{v'}$,
 where $v'$ is an elementary formula of the type $Y'_1$.  Now
 $$
\eqalign{
&s^0_3 s_1 \ol{v} = s^0_3s_2^0s\ol{u} = \ol{v'},\cr
&(s^0_3 s^0_2) (\ol{su}) = \ol{v'}.\cr}
$$
Thus, the formula $su$ is equivalent to the elementary formula,
i.e., $su \in H^0(Y)$.

The theorem is proved.

\subhead{4.2. Additional remarks}
\endsubhead

Note first of all that if we specify the type $X$ of the formula
$u$, we simultaneously distinguish  space where the value of the
formula $\Val_f(u)$ lies.  The formula $y^2\equiv 2px$ determines
a parabola in the two dimensional real space if $X=\{x, y\}$.  The
same formula for $X=\{ x, y, z\}$ gives a cylinder in the
three-dimensional space.

We view the formulas $u$ and $v$ as different ones if they are of
different types.  On the other hand, we may consider another
equivalence of such formulas.

Let $u$ and $v$ be of the types $X_1$ and $X_2,$ respectively.  We
consider them as equivalent ones if for some $Y$ containing $X_1,
X_2$ the equality $s^0_1 u = s^0_2 v$, where $s^0_1$ and $s^0_2$
are defined as above, holds.  For $X_1 = X_2$ the formulas $u$ and
$v$ are equivalent if they coincide.

Using this new equivalence we can link one-sorted and multi-sorted
Halmos algebras.  We will return to this later.

\proclaim{Proposition 1} For every homomorphism $\sigma:
\Hal_\Theta (\Phi) \to \Hal_\Theta (G)$ there exists a model $(G,
\Phi, f)$ with the condition: $\sigma = \Val_f$.
\endproclaim


\def\om{\omega}

Proof. For every $n$-ary $\vp \in \Phi$ we fix a set $X_\vp = \{
x_1, \cdots, x_n\}$ and consider an atomic formula $\vp(x_1,
\dots, x_n)$ of the type $X = X_\vp$.  All these $\vp(x_1, \dots,
x_n)$ and equality $x \equiv y$ of the type $\{ x, y\}$ generate
the whole algebra $H_\Theta\Phi$. We denote this set of generators
by $M^0$.

Define
$$
f(\vp) = (\sigma\vp(x_1, \dots, x_n))^\pi.
$$
Here the bijection $\pi: \Hom (W(X), G) \to G^{(n)}$ if defined by
$$
\pi(\mu) = (\mu(x_1), \dots, \mu (x_n)),
$$
for every $\mu: W(X) \to G$.

This determines the model $(G, \Phi, f)$.

It is necessary to check that $\sigma=\Val_f$.  We have:
$$
\Val_f(\vp(x_1, \dots, x_n)) = f(\vp)^{\pi^{-1}} = (\sigma\vp(x_1,
\dots, x_n))^{\pi\pi^{-1}} = \sigma\vp(x_1, \dots, x_n).
$$
We used here the definition of the function $\Val_f$ on atomic
formulas.  As for equalities, note that by the definition of a
homomorphism, every  $\sigma$ has to be correlated  with
equalities. This means, in particular, that $\sigma (x\equiv y)$
has to be the equality in $\Hom(W(x,y), G)$.  By the definition we
have $\sigma (x\equiv y) = \Val_f(x\equiv y)$. Thus, the equality
$\sigma = \Val_f$ holds on the generators and, hence, it is always
true.

$\blacksquare$

 Consider in more details the kernels of
homomorphisms and filters in multi-sorted Halmos algebras.

Let $\tau: H \to H'$ be a homomorphism of multi-sorted Halmos
algebras.  For every $X \in \Gamma^0$ there is a homomorphism of
quantorian $X$-algebras with equalities
$$
\tau_X=
\tau(X) : H(X) \to H'(X).
$$
Its kernel $U(X) = \Ker(\tau_X)$ is a full coimage of the unit. It
is a Boolean filter in $H(X)$, compatible with quantifiers. The
condition of compatibility with quantifiers means that $\forall x
u \in U(X)$ holds for every $x \in X$ and every $u \in U(X)$.

Here, as usual, $\forall u = \neg (\exists x \neg u )$. We say that
the filter $U(X)$ is invariant in respect to universal quantifiers
$\forall x , x \in X$.

Note also the conditions of compatibility with the operations of
the type $s \in S_\Theta$:  for every $s: W(X) \to W(Y)$ and every
$u \in U(X)$ we have $su\in U(Y)$.

The kernel $\Ker (\tau) = U$ is a set of all $U(X) = \Ker
(\tau_X)$ for all $X \in \Gamma^0$.  All these $U(X)$ are
compatible with universal quantifiers and operations of the type
$s$. We defined a filter of a multi-sorted Halmos algebra.
Straightforward check shows that such filters and homomorphisms of
Halmos algebras are tied correctly. Also ideals of multi-sorted
Halmos algebras are defined.

The filter $U$ is called {\it trivial}, if it always consists only
of the unit of the algebra $H$.  If $U(X) = H(X)$ holds for every
$X$, then $U$ is called an {\it unproper} filter.  Such filters
are always present in every $H$. If there are no other filters in
$H$, then the algebra $H$ is called {\it a simple} one.

Repeating arguments from [Pl1], one can prove that every algebra
$\Hal_\Theta (G)$ and all its subalgebras are simple Halmos
algebras and these are the only simple Halmos algebras (see 4.3).

Besides, we prove (see Lemma 2) that every Halmos algebra is
semisimple: it is approximated by simple ones.  The proof of
Theorem 1 (see 4.3) uses all above.

Let us consider the notion of the support of an element of Halmos
algebra. Let $H$  be a Halmos algebra and $u$ its element of the
$X$ type. Its {\it support} $\Delta u$ is defined by the rule:
$$
\Delta u = \Delta_X u = \{ x \in X | \exists x u \neq u\}.
$$
This condition means that the given $x$  essentially
participates in the ``record" of the element $u$.
However, speaking about the record of the element, one should take
care.  For example, an element $su=v$ for $s:W(X) \to W(Y)$ and
$u\in H(X)$, has the type $Y$ and its support is calculated in
$Y$.  Then, everything that is included in the record of $u$ in
$H(X)$ does not participate in the record of the element $v$.

Thus, the variables which do not belong to the set $X$ may
participate in the record of the element of the given type $X$, if
operations of the type $s$ are included in this record.  These
observations are  formalized in the algebra of formulas
$\Hal_\Theta (\Phi)$. For the elementary formulas the support of
the formula and the variables participating in the record are
perfectly coordinated. As a matter of fact, the notion of the
support of an element in an arbitrary multi-sorted Halmos algebra
requires a special  consideration.  For example, what can be said
about the relations between the supports $\Delta_Y(su)$ and
$\Delta_X(u)$ if $u$ is of $X$ type and $s: W(X)\to W(Y)$ is
given?  We may suppose that such relations take into account
transitions from the variables $x \in X$ to the support of the
elements $sx\in W(Y)$ (see also [Pl1]).

\proclaim{Proposition 2}  If $\tau: H \to H'$ is a homomorphism of
 Halmos algebras, then $\Delta_X (\tau u) \subset \Delta_X (u)$
for every $u \in H$ of the type $X$.
\endproclaim

Proof.  Let $x \in X$ and $\exists x \tau(u) = \tau (\exists x u)
\neq \tau u$.  Then $\exists x u \neq u, \; \; x \in \Delta_X(u)$.

We refer to such  feature as a coordination between homomorphisms
and supports of the elements.

We have defined the set of generators $M^0$ earlier.  This set can
be viewed as a multi-sorted one.  The supports of the formulas
from $M^0$ are $\Delta \vp (x_1, \cdots, x_n) = X_\vp = \{ x_1,
\dots, x_n\}$, $ \Delta (x\equiv y) = \{ x, y\}$.  Let $H'$ be a
an arbitrary Halmos algebra. A multi-sorted mapping $\tau : M^0
\to H'$ is  called  coordinated with the supports of the elements,
if for every $\vp \in \Phi$ and every $X \in \G^0$ with $X_\vp
\cup \{x, y\} \subset X$ the formulas
$$
\Delta_X(\tau_\vp(x_1, \dots, x_n)) \subset X_\vp \; \; \hbox{\rm
and } \; \; \Delta_X \tau (x \equiv y))) \subset \{ x, y\}.
$$
hold. Every such $\tau$ is uniquely extended up to $\tau: M \to
H'$.

Take $w, w', w_1, \dots, w_n$ of the type $Y$ and consider $s:
W(X) \to W(Y)$ with $s(x) = w, s(y) = w'$, $s(x_1) = w_1, \dots,
s(x_n) = w_n$. We set
$$
\eqalign{ &\tau (\vp (w_1, \dots, w_n )) = s\tau (\vp (x_1, \dots,
x_n)),\cr
&\tau(w\equiv w') = s\tau (x \equiv y).\cr}
$$
Check that this definition does not depend on the choice of $s$.

Consider a set of variables
$X_0 = \{ x, y, x_1, \dots, x_n\} \subset X$ and let $s'$
coincide with $s$ on the set $X_0$.  Take $z \in X\setminus X_0$.
We have:
$$
\eqalign{
&s\tau (\vp (x_1, \dots, x_n)) = s\exists z \tau \vp (x_1, \dots,
x_n) = \cr
&= s'\exists z \tau\vp(x_1, \dots, x_n) = s'\tau \vp (x_1, \dots,
x_n).\cr}
$$
Here $\tau \vp (x_1, \dots, x_n) = \exists z \tau \vp(x_1, \dots,
x_n)$ by the definition of $\tau$ and $z$ is not in  the support
for $\tau \vp(x_1, \dots, x_n)$.  The transition from $s$ to $s'$
is fulfilled  according to the rules of Halmos algebra. Similarly,
$s\tau (x \equiv y) = s'\tau (x \equiv y)$.  Thus, $\tau: M \to
H'$.  We will check that the mapping $\tau$ is a correct one.
Hence, it is uniquely extended up to the homomorphism $\tau:
\Hal_\Theta (\Phi) \to H'$.

Take some $s: W(X) \to W(Y)$ and let $w, w', w_1, \dots, w_n$ be
of the type $X$.  Take also a set $X_0$ containing variables $x,
y, x_1, \dots, x_n$ and let $s_0: W(X_0) \to W(X)$ be defined as
$s$ above.  Then $\tau \vp (sw_1, \dots, sw_n) = \tau\vp (ss_0
x_1, \dots, ss_0 x_n) = s s_0 \tau \vp (x_1, \dots, x_n) = s\tau
\vp (s_0 x_1, \dots, s_0 x_n) = s \tau \vp (w_, \dots, w_n)$.
Similarly, $\tau (s w\equiv s w') = s \tau (w \equiv w')$.  Hence,
the mapping $\tau$ is correct.  This property has been mentioned
in Theorem 2. The following Proposition is also related to the
important properties of the algebra $\Hal_\Theta(\Phi)$.

\proclaim{Proposition 3} Let $\alpha: H\to H'$ be a surjective
homomorphism of two (multi-sorted) Halmos algebras. Let a
homomorphism $\beta:Hal_{\Phi\Theta}\to H'$ be given. Then, there
exists $:\gamma: Hal_{\Phi\Theta}\to H$ such that
$\alpha\gamma=\beta$.
\endproclaim

Proof. Let us use, first, an auxiliary remark. Let $h$ be an
element in $H$ of the type $X$. Then
$\Delta_X(\alpha(h))\subset\Delta_X(h)$. Show that there exists an
element $h_1$ in $H$ of the type $X$, such that
$\alpha(h)=\alpha(h_1)$ and $\Delta_X(h_1)=\Delta_X\alpha(h_1)$.

Let us take
$$
Y=\Delta_X(h)\backslash \Delta_X(\alpha(h)),
$$
and let $h_1=\exists(Y)h$. Since $Y$ does not intersect the
support of the element $\alpha(h)$ then
$$
\alpha(h_1)=\alpha(\exists Yh)=\exists(Y)\alpha(h)=\alpha(h).
$$

Besides that
$$
\Delta_X(h_1)=\Delta_X(h)\backslash Y=\Delta_X\alpha(h)=
\Delta_X\alpha(h_1),
$$
and $\Delta_X(h_1)=\Delta_X\alpha(h_1)$.

Now we use this remark. For every basic element $u$ of the type
$X$ in $\Hal_{\Phi\Theta}$ we take $\beta(u)$. This is an element
of the type $X$ in $H'$. Take an element $h\in H$ with $\alpha(h)=
\beta(u)$. Let us take $h_1$ such that
$\alpha(h_1)=\alpha(h)=\beta(u)$ and
$$
\Delta_X(h_1)=\Delta_X\alpha(h_1)=\Delta_X\beta(u)\subset
\Delta_X(u).
$$
Denote $\gamma(u)=h_1$. Thus, for every $u$ we have:
$$
\Delta_X(\gamma(u)=\Delta_X(h_1)\subset\Delta_X(u).
$$
This means that the homomorphism $\gamma: \Hal_{\Phi\Theta}\to H$
is defined.

For every basic $u$ we have
$(\alpha\gamma)(u)=\alpha(h_1)=\alpha(h)=\beta(u)$. Then
$(\alpha\gamma)(u)=\beta(u)$ for every $u$, i.e.,
$\alpha\gamma=\beta$.

Now we pass to the proof of Theorem 3.

Proof. This commutative diagram follows from the definitions:


$$
\CD
H @[2]>\Val_f >> \Hal_\Theta (G)
\\
@[2]/SE/// @.@. \; @/NE//\Val_f / \\
@.  \Hal_\Theta (\Phi)
\endCD
$$
This means that $\Val_f (u) = \Val_f (\ol u)$ holds for every
formula $u \in H$.  If, now, $\ol u = \ol v$, then $\Val_f (u)=
\Val_f (\ol u ) = \Val_v (\ol v) = \Val_f (v)$.  It is true for
every model $(G, \Phi, f)$. Thus, $\ol u = \ol v $ implies
semantical equivalence of the formulas $u$ and $v$.

Conversely, let $u$ and $v$ be semantically equivalent. Apply the
fact that  every Halmos algebra, is semisimple. In particular, the
algebra $\Hal_\Theta (\Phi)$ is semisimple. We have a system of
homomorphisms $\sigma_\a: \Hal_\Theta (\Phi) \to H_\Theta (G_\a),
\a \in I$, such that if $\sigma_\a (u) = \sigma_\a (v)$ for all
$\a$, then $u = v$.  Every $\sigma_\a$ here can be represented as
$\Val_{f_\a}$ by a model $(G_\a, \Phi, f_\a)$.  If $u$ and $v$ are
semantically equivalent formulas in $\Hal_\Theta (\Phi)$, then
$\Val_{f \a}(u) = \Val_{f\a}(v), \sigma_\a (u) = \sigma_\a(v)$ for
all $\a \in I$. This gives $u = v$.

Note that similar reasoning does not work, for example, for the
algebra $\tilde H$.  Here $\sigma: \tilde H \to \Hal_\Theta (G)$
can be represented as $\sigma = \Val_f $ only if $\sigma$ is
coordinated with the relations $(*)$.

\subhead{4.3  Some  proofs}
\endsubhead

Here we give proof of the theorem 1.

 Besides filters, consider
ideals. An ideal $U$ of the Halmos algebra $H$ selects in $H(X)$ a
Boolean ideal $U(X)$ for every $X\in \G^0$,  which is invariant
under quantifiers $\exists x, x \in X$, and for every $s: W(X) \to
W(Y)$ we have: $u \in U(X)$ implies $su\in U(Y)$. Ideals and
filters are  dual.   The zero ideal is an ideal $U$ for which all
$U(X)$ are zeroes.  The ideal $U$ is trivial, if $U(X) = H(X)$ for
all $X$. The algebra $H$ is simple if and only if it has no other
ideals.

\proclaim{Lemma 1} Every $\Hal_\Theta (G)$ and all its subalgebras
are simple.
\endproclaim

Proof. Take a nonempty subset $A$ in $\Hom (W(X), G)$.  Apply a
quantifier $\exists (X) = \exists x_1\dots \exists x_n, X = \{
x_1, \dots, x_n\}$.  The homomorphism $\mu: W(X) \to G$ belongs to
the set $\exists (X) A $ if $\mu$ coincides with some $\nu: W(X)
\to G$ outside $X$.  Since there are no variables outside $X$,
then every $\mu$ belongs to $\exists (X) A;$ $  \exists X A = \Hom
(W(X), G)$.  It is the unit of the algebra $\Bool (W(X), G)$.

Let now $H$ be an arbitrary Halmos algebra such that $\exists Xa =
1$ holds for every nonzero element $a \in H (X)$ for every $X \in
\G^0$.

Check that $H$ is simple. Let $U$ be a nonzero ideal in $H$.  For
some $X$ we have $U(X) \neq 0$ and in $U(X)$ there is a nonzero
element $a$.  We have $\exists x a = 1$ that is $U(X)$ contains a
unit. Then $U(X) = H(X)$.
 For any $s: W(X) \to W(Y)$ we have $s u \in U(Y)$ for $u \in
 U(X)$.  If $u$ is a unit in $H(X)$, then $su$ is a unit in
 $H(Y)$, $U(Y)$ contains a unit, and $U(Y) = H(Y)$.  The ideal $U$
 is trivial, that is the algebra $\Hal_\Theta (G)$
 and all its subalgebras are simple.

 We will prove later that every simple Halmos algebra is a
 subalgebra of some $\Hal_\Theta (G)$.

 We use the following definition.

 Let $H$ be a multi-sorted Halmos algebra.  We say that this
 algebra is {\it one-sorted representable} if there exists a one-sorted
 Halmos algebra $H(X^0) = H^0$ such that

 1.  \ To each inclusion $s^0: W(X) \to W(X^0)$ there corresponds a
 Boolean inclusion $s^0: H(X) \to H(X^0)$, naturally correlated
 with quantifiers and equalities.

 2. \ To each commutative diagram
 $$
 \CD
 W(X) @>s>> W(Y) \\
 @V s^0_X VV  @VVs^0_Y V \\
 W(X^0) @> s' >> W(X^0)
 \endCD
 $$
there corresponds a diagram
$$
\CD
H(X) @>s>> H(Y) \\
@V s^0_X VV @VVs^0_Y V\\
H^0 @> s' >> H^0
\endCD
$$
As in [Pl1], we can prove that every $H$ is representable.  It is
easy to check directly  for the algebras $\Hal_\Theta (\Phi)$ and
$\Hal_\Theta(G)$.  This notion gives ties between multi-sorted and
one-sorted Halmos algebras.  We  use this fact of representability
in the following theorem.

\proclaim{Theorem 5} Every Halmos algebra is semisimple.
\endproclaim

Proof.  Let the representation $H\to H^0$ be given. For every
element $a \in H^0$ there is a maximal Boolean ideal $V$ not
containing the element $a$. Let $V_*$ be the ideal of Halmos
algebra $H^0$ defined by the rule: an element $u \in H^0$ belongs
to $V_*$ if $\exists (X^0) u \in V$.
 We may think that the algebra $H^0$ is locally finite.  Then,
$\exists X^0 u = \exists (X) u$, where $X \in \G^0$.  It is proved
that $V_*$ is a maximal ideal of the Halmos algebra $H^0$ and $V_*
\subset V$ (see [Pl 1]).

Relying on $V_*$ we build an ideal $U$ in the multi-sorted algebra
$H$. For every finite $X \subset X^0$ denote by $U(X)$ a set of
elements $u \in H(X)$ for which $s^0 u \in V_*$.  Here, $s^0: H(X)
\to H(X^0) = H^0$ is the natural embedding.

It is easy to check that $U(X)$ is a Boolean ideal in $H(X)$,
preserving the quantifier $\exists (X)$. Check also that if $u \in
U(X)$ then $s u \in U(Y)$ for every for every $s: W(X) \to W(Y)$.
Consider a commutative diagram
$$
\CD
H(X) @> s>> H(Y) \\
@Vs^0_X VV @VVs^0_Y V \\
H^0 @>s' >> H^0
\endCD
$$
We have:  $s^0_Y s = s's_X^0$.  We need $s^0_Y (su) \in V_*$.  Now
$s^0_Y su = s's^0_X u$.  By the condition, $s^0_X u \in V_*$.
Then,  $s's^0_X u \in V_*$, $s^0_Y (su) \in V_*$.  Hence, $U$ is
an ideal of the algebra $H$.  Let us show that $U$ is a maximal
ideal.  Let $U_1$ be a greater one.  An ideal $U_*$ in $H^0$,
determined by $U_1$, corresponds to it, and $U_* > V_*$.  Since
$V_*$ is maximal, then $U_* = V_*$.  It should be $U_1(X) > U(X)$
for some finite $X$.  Let $u \in U_1 (X) \setminus U(X)$.  We have
$s^0 u \in U_* = V_*$, and then $u \in U(X)$.

The arising contradiction means that the ideal $U$ is maximal.

Take, further, an element $a \in H(X)$ and $s^0 a \in H^0$. The
element $s^0a$ determines the maximal Boolean ideal $V$ and,
correspondingly, the ideal $V_*$ of the Halmos algebra. Let $U$ be
the filter of the Halmos algebra $H$ corresponding to $V_*$.
Denote $U= U_a$.  Take these $U_a$ for all nonzero elements $a \in
H$ of the different types $X$.


All $U_a$ are maximal ideals in $H$ and their intersection is zero
ideal.  This means that the algebra $H$ is semisimple, since it is
approximated by simple algebras.

\proclaim{Lemma 2} Every simple algebra $H$ is isomorphic to some
subalgebra of the algebra of the type $\Hal_\Theta(G)$.
\endproclaim

Proof. Let $H$ be simple and $H_0$ some set of generating elements
in $H$. For every $h \in H_0$ take a symbol of relation $\vp$.  If
$h$ is of the type $X = \{ x_1, \dots, x_n\}$ then $\vp$ is
$n$-ary. Associate $\vp(x_1, \dots, x_n) \to h$. Recall that we
consider every $H$ as an algebra with equalities. Then the formula
$w \equiv w'$ is associated with the corresponding equality in
$H$.

Collect  all $\vp$ into a set $\Phi$ and take the corresponding
free algebra $\tilde H$ for  $\Phi$ and $\Theta$.   This gives a
homomorphism of the algebra $\tilde H$ on the algebra $H$. Let the
ideal $U$ be a kernel of such homomorphism. We use here one-sorted
representation for the algebras $\tilde H$ and $H$, and let an
ideal $U_0$ correspond to the ideal $U$ in the algebra $\tilde
H^0$.  The algebra $\tilde H^0$ is a free one-sorted Halmos
algebra.  Take a filter $T_0$ corresponding to the ideal $U_0$.
Since the algebra $H$ is simple, the filter $T_0 $ is a maximal
one and $U$ is a maximal ideal.  Filter $T$ in the algebra $H$,
associated with the ideal $U$, corresponds to the filter $T_0$.

The filter $T_0$ determines some model $(G, \Phi, f)$, in which
the set of formulas $T_0$ and  also all the formulas of the filter
$T$ hold.  Take $\Val_f\colon \tilde H \to \Hal_\Theta (G)$ by
$f$. Here, $T=\Ker(\Val_f)$.  From this follows that the algebra
$H$ is isomorphic to some subalgebra in $\Hal_\Theta (G)$.

Lemmas 1 and 2 together with Theorem 5  imply Theorem 2 with the
help of characterization of the varieties by Birkhoff type theorem
for Halmos algebras.

\head {5. Algebraic geometry in first order logic}
\endhead

\subhead{5.1 Sets of formulas and algebraic sets}
\endsubhead

Fix a finite set $X$. Let $W=W(X)$ be the free algebra over $X$ in
the given variety $\Theta$.  The set $\Phi$ of symbols of
relations and a model $(G, \Phi, f), G  \in \Theta$ are also
fixed.  We view the set of homomorphisms $\Hom (W,G)$ as an affine
space (see also 1.3) .

Now, consider the sets $A$ of points $\mu\colon W\to G$ in this
space and the sets $T$ of formulas in the
algebra
$\Hal_{\Phi\Theta} (X)$.  Establish the following Galois
correspondence for the given model $(G, \Phi, f)$:
$$
\cases
T^f=A=\bigcap\limits_{u\in T} Val_f(u) = \{ \mu \big| T
\subset \Log\Ker (\mu) \}\\
A^f=T=\{ u \big| A \subset Val_f (u) \} = \bigcap\limits_{\mu
\in A} \Log \Ker (\mu)\endcases
$$
In each row we have three equalities, first two of which are the
definition, while the third one is an easily checked equality.
These equalities actually hold  and this is a Galois
correspondence. It generalizes the standard correspondence in the
classical algebraic geometry and, also, the Galois correspondence
for the equational geometry in the variety $\Theta$ for the
algebra $G \in\Theta$.

The set $A$ of the type $A = T^f$ for some $T$ we call an {\it
algebraic set, or closed set (algebraic variety)}, over the model
$(G, \Phi, f)$. It is also called an elementary set.

The set $T$ of the type $T = A^f$ for some $A$ is a Boolean filter
in the Boolean algebra $\Hal_{\Phi\Theta} (X)$. We call it an
$f$-{\it closed Boolean filter}. If $A$ is an algebraic set, then
we consider the filter $T=A^f$ as a theory of the set $A$ in the
algebra $\Hal_{\Phi\Theta} (X)$.

Now we consider the Boolean algebra $\Hal_{\Phi\Theta} (X) \big/
A^f$.  It is the coordinate algebra for the given $A$. We have the
following embedding: $$ \Hal_{\Phi\Theta} (X)\big/A^f \to
\prod\limits_{\mu \in A} \Hal_{\Phi\Theta} (X) \big/\Log \Ker
(\mu).
$$

The right side is the cartesian product of two-element algebras.

We can consider the closure $A^{ff}$ for every set $A$ and
$T^{ff}$ for every $T$.  The next proposition gives a
straightforward relation between  $T$ and $T^{ff}$.

\proclaim{Proposition 4}  $v\in T^{ff}$ if and only if the
(infinitary) formula $(\bigwedge\limits_{u\in T} u) \to v$ holds
in the model $(G, \Phi, f)$.
\endproclaim

Proof.  The given formula holds in $(G, \Phi, f)$ if and only if
for arbitrary point $\mu: W \to G$ the inclusion $\mu\in\Val_f(u)$
follows from that of $\mu \in \Val_f(u)$ for all $u \in T$. Let
the formula  hold in the model.  Prove that $v \in T^{ff}$.  By
definition,
$$
T^{ff} = \bigcap\limits_{\mu\in T^f} \Log \Ker (\mu).
$$
We need to prove that if $\mu \in T^f$, then $v \in \Log \Ker
(\mu)$.
We have
$$
T^f=\bigcap\limits_{u \in T} \Val_f (u).
$$
Hence, $\mu \in T^f$ gives $\mu \in \Val_f(u)$ for all $u \in T$.
Besides, $\mu \in \Val_f(v), v \in \Log \Ker (\mu)$ by the condition
of the proposition.
Thus,
$
v \in T^{ff}.
$

Let now $v \in T^{ff}$.  Prove that the formula holds in $(G,
\Phi, f)$.

Take some $\mu: W \to G$ and let $\mu \in \Val_f(u)$ for all $u
\in T$.  We need to prove that $\mu \in \Val_f(v)$.  We have:
$$
\mu\in\bigcap\limits_{u\in T} Val_f (u) = T^f, \quad T^{ff}
\subset \Log \Ker (\mu).
$$
Here $v \in T^{ff}$ means that $v \in \Log \Ker (\mu), \mu \in
Val_f(v)$.

The proposition is proved.

\proclaim{Corollary}
  If $u$ is a formula,
then $v \in u^{ff}$ if and only if the formula $u \to v$ holds in
the model.
\endproclaim

This  is a version of the Hilbert's Nullstellensats in the
algebraic geometry.

\subhead{5.2 The Galois correspondence and morphisms}
\endsubhead

Further the Galois correspondence will be associated with the
morphisms of the categories.  In the case of $\Hal_\Theta (G)$ the
homomorphisms $s\colon W(X) \to W(Y)$ yield the mappings $s_*$ and
$s^*$ for Booleans.  Now, let us define the action of $s_*$ and
$s^*$ on the set of formulas $T$.  If $T \subset \Hal_\Theta
(\Phi)(Y)$, then $s_* T$ is the set of formulas in $\Hal_\Theta
(\Phi) (X)$, determined by the rule:
$$
u \in s_* T \Leftrightarrow s_* u \in T.
$$
If $T\subset \Hal_\Theta (\Phi)(X)$, then $s^*T\subset \Hal_\Theta
(\Phi)(Y)$,
and
$$ s^* T = \{ s_* u | u \in T\}.
$$
\proclaim{Theorem 6}  The rules of compatibility are of the form:

1. If $T \subset \Hal_\Theta (\Phi)(X)$, then
$$
(s^*T)^f = s_* T^f  = sT^f.
$$

2. \ If $A \subset \Hom(W(Y), G)$, then
$$
(s^*A)^f=s_*A^f.
$$
\endproclaim

Proof.  We use the definitions:
$$
T^f=\bigcap\limits_{u\in T} \Val_f(u);\ \ A^f = \{ u \big|
A\subset \Val_f(u)\}.
$$
Let us prove the first rule.  Due to $s^*T \subset
\Hal_\Theta(\Phi) (Y)$, we have $(s^*T)^f, \ s_*T^f$ \break
$\subset \Hom (W(Y), G)$. Consider the point $\nu: W(Y)\to G$. Let
$\nu \in (s^*T)^f = \bigcap\limits_{u\in s^*T} \Val_f (u),$
 $ u = sv, v \in T$.

  For every
$v\in T$ we have:
$$
\nu \in\Val_f(sv) = s\Val_f(v) = s_*\Val_f(v).
$$
Hence, for every $v \in T$ we have
$\nu s\in\Val_f (v)$, and $ \nu s \in T^f = \bigcap\limits_{v\in T}
\Val_f(v), \nu \in s T^f$.

Conversely, let $\nu \in sT^f, \nu s \in T^f=\bigcap_{v\in T}
\Val_f(v)$.  Then $\nu s \in \Val_f (v)$, $\nu \in s \Val_f(v) =
\Val_f (sv)$ for every $v \in T$.
Therefore
$$
\nu \in\bigcap \limits_{v \in T}\Val_f(sv) =
\bigcap\limits_{u \in s^* T} \Val_f (u) = (s^* T)^f.
$$
Consider the second rule.  Here, $s^* A \subset \Hom (W(X), G)$
and $(s^*A)^f$, $s_* A^f \subset \Hal_\Theta (\Phi)(X)$. Let $u
\in (s^*A)^f$, which implies $s^*A\subset \Val_f(u)$. Hence, $\nu
s \in \Val_f(u)$, $\nu \in s \Val_f(u) = \Val_f(su)$, $A \subset
\Val_f (su), su \in A^f, u \in s_*A^f$ for every $\nu \in A$.

Conversely, let $u \in s_* A^f, su\in A^f, A \subset\Val_f(su) = s
\Val_f(u)$.  Then $\nu s \in \Val_f(u)$,  $ s^* A \subset
\Val_f(u), u \in (s^*A)^f$ for every $\nu \in A$.

The rules are proved.

  The given rules induce a good
compatibility between the Galois correspondence in logic and the
morphisms in the categories $\Hal_\Theta ( \Phi)$ and
$\Hal_\Theta(G)$.

\subhead{5.3  Lattices and topology}
\endsubhead

 We introduce
here the functions $Cl_f$ and $Alv_f$, defined for an arbitrary
$W=W(X)$.

The set $Alv_f(W)$ is the set of all algebraic sets for the given
model $(G,\Phi,f)$ in the affine space $\Hom (W(X),G)$.

The set $Cl_f(W)$ is the set of all $f$-closed filters $T$ in the
Boolean algebra $\Hal_\Theta (\Phi)(X)$.

\proclaim{Proposition 5} The sets $Alv_f(W)$ and $Cl_f(W)$ are
dually isomorphic lattices.
\endproclaim

{Proof.}  Let $A, B \in Alv_f(W),$ $ A = T_1^f$, $B = T^f_2$. Then
$A\cap B = (T_1 \cup T_2)^f, A \cap B \in Alv_f(W)$.

1.  \ If $A = T^f$ is an algebraic set, then $sA$ is also the
algebraic set.

2. \ If $T = A^f$ is $f$-closed, then $sT = s_*T$ is also
$f$-closed.

Denote the set of all $u\wedge v, u \in T_1$, $ v \in T_2$ by
$T_1\wedge T_2$.  Then  $A\cup B = (T_1\cap T_2)^f$, $A\cup B \in
Alv_f(W)$.  So, $Alv_f(W)$ is a lattice, which is a sublattice in
the lattice $\Bool (W(X), G)$, which is a distributive lattice.

Let, now $T_1, T_2 \in Cl_f(W), T_1=A^f$, $T_2=B^f$.  Then $T_1\cap
T_2 = A^f\cap B^f = (A\cup B)^f$, $T_1 \cap T_2 \in Cl_f(W)$.

However, we can not proclaim that $T_1 \cup T_2 $ also belongs to
$Cl_f(W)$, since the union of filters can be not a filter.

Thus, we introduce a new operation:
$$
T_1 {\ol\cup}  T_2 = (T_1 \cup T_2)^{ff}.
$$
As a result, we have a lattice $Cl_f(W)$.  Check, that the
lattices $Alv_f(W)$ and $Cl_f(W)$ are dually isomorphic.  This
will imply that the lattice $Cl_f$ is also distributive.

The transition $f$ determines the bijection between algebraic sets
in $\Hom (W,G)$ and $f$-closed filters in $\Hal_\Theta(\Phi)(X)$.
In other words, we have bijections $$ \eqalign{ &f\colon Alv_f(W)
\to Cl_f(W),\cr &f\colon Cl_f(W) \to Alv_f(W).\cr} $$ Let $A, B
\in Alv_f(W)$.  Then $$ (A\cap B)^f = (T_1^f \cap
T_2^f)^f=(T_1\cup T_2)^{ff} = T_1\ol \cup T_2, $$ where $T_1=A^f,
T_2 = B^f$.
Thus,
$$
(A\cap B)^f= A^f {\ol \cup} B^f.
$$
Now, we need to prove
$$
(A\cup B)^f = A^f\cap B^f.
$$
Since the set $A\cup B$ is closed, we need to check that $(A^f\cap
B^f)^f = A\cup B$.
We have:
$$
A^f\cap B^f = (A\cup B)^f, (A^f\cap B^f)^f = (A\cup B)^{ff} =
A\cup B.
$$
The proposition is proved.

We can consider the function $Cl_f\colon\Theta^0 \to Set $ as a
functor.  However, there is no compatibility with the lattice
$Cl_f(W)$ in this case.  On the other hand, the function
$Alv_f\colon \Theta^0 \to Set$ is compatible with the lattice and
one can look at the functor
$$
Alv_f:\Theta^0 \to Lat,
$$
where $Lat $ is the category of lattices.

Consider the function $Alv_f$ also as a multi-sorted set. We have
$Alv_f(X) = Alv_f(W(X))$ for every $X \in \Gamma^0$. Now, $Alv_f$
is a subset in the algebra $\Hal_\Theta(G)$ which can be treated
as a multi-sorted lattice, invariant under the operations of $s
\in S_\Theta$ type.

Let us return to the homomorphism
$\Val_f\colon\Hal_\Theta(\Phi)\to \Hal_\Theta(G)$.  Let $R_f$ be
the image of this homomorphism. It s a subalgebra in $\Hal_\Theta
(G)$.  Every $A \in R_f $ has the form $A=\Val_f(u)$.  It is an
algebraic set, determined by one-element set $T$, consisting of an
element $u$.  We call such algebraic sets simple ones.

All simple algebraic sets for the given model $(G,\Phi,f)$ form a
subalgebra in the Halmos algebra $\Hal_\Theta (G)$.

As we know, every formula $u$ is equivalent to some elementary
formula $v$.  In the sequel we will relate the varieties
$\Val_f(u)$ and $\Val_f(v)$.

Note, further, that the main attention in groups and other
structures in classical algebraic geometry and equational geometry
is paid to the properties of individual algebraic sets.  The most
important goal is to obtain some description of the system of all
solutions of equations.

We pay main attention here  to lattices and categories of
algebraic sets. The properties of individual algebraic sets is a
separate topic.

Let us give some remarks on topologies in the spaces $\Hom (W,G)$,
associated with algebraic sets, namely, generalized Zarisski
topologies.  Closed sets here are algebraic sets.  It is natural,
however, to restrict ourselves with the sets, determined by a
collection of positive formulas, recorded without negation.  The
good topology assumes also that there are no quantifiers (it can
be explained).  We speak of algebraic sets, determined by the
collections of universal positive formulas, and consider such
Zariski topology.

\demo{Example}  Proceed from the classical variety $\Theta$ of
commutative and associative algebras with the unit over the field
of real numbers. Take a unique relation (order relation) as
$\Phi$.  The model $(G, \Phi, f)$ is a field of real numbers, in
which the order relation is naturally realized.  Take $X = \{ x,
y\}$.  The corresponding algebra $W$ is the algebra of real
polynomials of two variables $x$ and $y$.  The space $\Hom (W, G)$
is realized as the real space.

It is easy to understand that in this case the corresponding
generalized Zariski topology coincides with the natural topology,
while the usual Zariski topology does not coincide with this
natural one.  Take, in particular, a disk $x^2 + y^2\le a^2$. It
is closed in the natural topology, but not closed in the usual
Zariski topology.  The same disk can be given by the formula
$\exists z (x^2+y^2+z^2=a^2)$.  The set $\Phi$ is empty here, but
there is a quantifier and the set $X$ consists of three variables.
\enddemo

\subhead{5.4 The categories $K_{\Phi\Theta} (f)$ and
$C_{\Phi\Theta}(f)$}
\endsubhead

Besides the variety of algebras $\Theta$, let us fix also the set
of symbols of relations $\Phi$ and a model $(G, \Phi, f),
G\in\Theta$.  These data determine the category
$K_{\Phi\Theta}(f)$ of algebraic sets for the given variety $f$.
Its objects have the form
$$(X,A),
$$
where $A$ is an algebraic set in the given logic for the given
model.  The set $A$ lies in the affine space $\Hom(W(X), G), A =
T^f$ for some set of formulas $T$, but this $T$ is not fixed.

The morphisms  have the form
$$
(X,A) \to (Y,B).
$$
We start from $s\colon W(Y) \to W(X) $ in $ \Theta^0$.
We have:
$$
\tilde s\colon \Hom (W(X), G) \to \Hom (W(Y), G).
$$
We say that $s$ is admissible for $A$ and $B$ if $\tilde s (\nu) =
\nu s \in B$ for $\nu \in A$.  For every such $s$ we have a
mapping
$$
[s]\colon A\to B.
$$
Now we consider weak and exact categories $K_{\Phi\Theta}(f)$.

In the weak category morphisms are of the form
$$
s\colon (X, A) \to (Y, B),
$$
where $s$ is admissible for $A$ and $B$.

In the exact category morphisms are
$$
[s]\colon (X, A) \to (Y, B).
$$

Consider the notion of admissible triple $(s, A, B)$ in more
detail.
 For $s\colon W(Y) \to W(X)$ we have
 $$
 s_* \colon \Hal_{\Phi\Theta} (Y) \to \Hal_{\Phi\Theta}(X).
 $$
 Let $T_2$ and $T_1$ be some sets of formulas in
 $\Hal_{\Phi\Theta}(Y)$ and $\Hal_{\Phi\Theta}(X)$, respectively.
 We say that $s_*$ is admissible for $T_2$ and $T_1$ if $s_*u \in
 T_1$ for $u \in T_2$.

 \proclaim{Proposition 6} The homomorphism $s\colon W(Y) \to W(X)$
 is admissible for algebraic sets $A$ and $B$ if and only if $s_*$ is
 admissible for $T_2 = B^f$ and $T_1=A^f$.
 \endproclaim

 {Proof.}  Note first, that the inclusion $\nu s \in B$ for
 all $\nu \in A$ means that $A\subset sB=s_* B$.
 This holds in the category $\Hal_\Theta(G)$.
 Applying the transition $f$, we get $(s_* B)^f \subset A^f$.  We
 claim further: $s^* B^f \subset (s_* B)^f$.
 This rule was not among the previously mentioned ones.

 Take $u \in s^* B^f$, $u = s_* v, v \in B^f$, $B\subset
 \Val_f(v)$.
 From this follows $s_* B \subset s_* \Val_f(v) = \Val_f(s_*
 v)=\Val_f (u)$.  Hence, $s_* B\subset \Val_f(u),  u \in (s_*B)^f$.
 The last inclusion holds for ever $u \in s^* B^f$.  Thus, we have
 $s^* B^f\subset (s_*B)^f$.

Let now $s$ be admissible for $A$ and $B$.  We have $A\subset sB$
and $(s_*B)^f\subset A^f = T_1$.  Using $s^*B^f\subset (s_*B)^f$,
we get $s^*B^f=s^*T_2\subset T_1$.  This means that $s_* u \in
T_1$ for $u \in T_2$ and $s_*$ is admissible for $T_2 $ and $T_1$.

Conversely, let now $s^*T_2\subset T_1, s^*B^f\subset A^f$. Again
applying $f$, we get
$$
T_1^f = A\subset (s^*T_2)^f = s_* T^f_2 = s_*B=sB
$$
and $s$ is admissible for $A$ and $B$.  The proposition is proved.

 In this case we have a commutative diagram of Boolean
 homomorphisms
$$ \CD \Hal_{\Phi\Theta}(Y) @>s_*>>
\Hal_{\Phi\Theta}(X)\\ @V\mu_Y VV @VV\mu_X V\\
\Hal_{\Phi\Theta}(Y)\big/B^f @>\ol{s_*}>> \Hal_{\Phi\Theta}
(X)\big/A^f\\
\endCD
$$
 where $\mu_Y$ and $\mu_X$ are natural homomorphisms.

\proclaim{Proposition 7}  Let $s_1$ and $s_2$ be admissible for
$A$ and $B$.  Then $\ol{s_{1*}} = \ol{s_{2*}}$ follows from
$[s_1]=[s_2]$.
\endproclaim

{Proof.}  Let $[s_1]=[s_2]$.  Then $\nu s_1=\nu s_2$ for every
$\nu \in A$.  Consider the following two diagrams:

$$
\CD
\Hal_{\Phi\Theta}(Y) @>s_{1*}>> \Hal_{\Phi\Theta}(X)\\
@V\mu_Y VV  @VV\mu_X V\\
\Hal_{\Phi\Theta}(Y)/B^f @>\ol{s_{1*}}>> \Hal_{\Phi\Theta}(X)/A^f
\endCD
$$

\noindent
and

$$
\CD
\Hal_{\Phi\Theta}(Y) @>s_{2*}>> \Hal_{\Phi\Theta} (X)\\
@V\mu_Y VV @VV\mu_X V\\
\Hal_{\Phi\Theta}(Y)/B^f @>\ol{s_{2*}}>> \Hal_{\Phi\Theta} (X)/A^f
\endCD
$$

We claim  $\ol{s_{1*}} = \ol{s_2*}$, that is, $\ol{s_{1*}} \mu_Y (u)
= \ol{s_{2*}} \mu_Y(u)$ for all $u \in \Hal_{\Phi\Theta}(Y)$.
Hence, $\mu_Xs_{1*} (u) = \mu_X s_{2*} (u)$.

The idea is to prove that $s_{1*}(u)$ and $s_{2*}(u)$ coincide
modulo the filter $A^f$.  This means that $(\neg s_{1*} (u) \vee
s_{2*}(u))\wedge (s_{1*}(u) \vee \neg s_{2*} (u)) \in A^f$.
Denote the left part by $v$.  To prove $v\in A^f$ is the same as
to prove the inclusion $A\subset \Val_f(v)$.  So, we get $\nu
\in\Val_f(v), \forall \nu \in A$. We have
$$
\eqalign{ &\Val_f (v) = (\neg s_{1*}\Val_f (u)\vee s_{2*} \Val_f
(u))\wedge\cr
&(s_{1*}\Val_f(u)\vee \neg s_{2*} \Val_f(u)).\cr}
$$
Check that $\nu$ is included in the first parenthesis.  Suppose
this is not true.

Then $\nu \in s_{1*} \Val_f(u), \nu s_1 \in\Val_f(u)$.
The equality $\nu s_1 = \nu s_2 $ and the inclusions $\nu s_2 \in
\Val_f(u), \nu \in s_{2*} \Val_f (u)$ lead to a contradiction.
Hence, $\nu$ is included in the first parenthesis, and,
analogously, in the second one.  Therefore, $A \subset \Val_f(v)$,
$v\in A^f$.  Thus, $s_{1*} (u)$ and $s_{2*}(u)$ coincide modulo
the filter $A^f$.  This leads to $\mu_Xs_{1*} (u) = \mu_Xs_{2*}
(u), \ol{s_{1*}} \mu_Y (u) = \ol{s_{2*}} \mu_Y(u), \ol{s_{1*}} =
\ol{s_{2*}}$.
This completes the proof of the proposition.

Let us discuss the converse statement.  Let $\ol{s_{1*}} =
\ol{s_{2*}}$.
Then $\ol{s_{1*}}\mu_Y(u) = \ol{s_{2*}}\mu_Y(u), \forall u \in
\Hal_{\Phi\Theta} (Y)$.  Therefore, $\mu_X s_{1*}(u) =
\mu_Xs_{2*}(u)$.
Hence, $s_{1*}(u)$ and $s_{2*} (u)$ are equivalent modulo $A^f$.

Consider the same $v = (\neg s_{1*} (u) \vee s_{2*}(u)) \wedge
(s_{1*}(u) \vee \neg s_{2*} (u)) \in A^f, A\subset \Val_f(v)$.
As before, $\nu \in A$ implies $\nu \in \Val_f(v)$ for every $u$.

Now we are interested in the interpretation of the last inclusion.
Let $\nu \in \neg s_{1*}\Val_f(u)$, $ \nu s_1 \notin \Val_f(u)$.
Using the second parenthesis in the record of $v$, we come to $\nu
s_{2}\notin\Val_f(u)$.

If $\nu s_1 \in\Val_f(u)$, then the inclusion $\nu s_2
\in\Val_f(u)$ follows from the first parenthesis. Does this lead
to $\nu s_1 = \nu s_2$? The problem arising here is that if
$$
\exists \nu\in A, \; \nu s_1 \neq \nu s_2, \; \exists u \in
\Hal_{\Phi\Theta} (Y),
$$
with the properties $\nu s_1\in \Val_f (u), \nu s_2\notin
\Val_f(u)$.  In other words, if there exists $u$ which separates
two different points $\nu s_1$ and $\nu s_2$? It is not true in
general. However, in the case when in the corresponding Zarisski
topology every point is $f$-closed this property holds.  In this
case we have the statement, opposite to Proposition 7. Recall also
that the opposite statement is true in the equational geometry.


The objects of the category $\bbc_{\Phi\Theta}(f)$ are of the form
$\Hal_{\Phi\Theta} (X)/T$, where $T$ is $f$-closed Boolean filter
in $\Hal_{\Phi\Theta} (X)$.  The morphisms are homomorphisms of
such algebras, induced by homomorphisms $s\colon W(X) \to W(Y)$.
As before, they are of the type $\ol{s}_*$, i.e., not arbitrary
homomorphism of the given Boolean algebras. The transition
$$
(X,A) \to \Hal_{\Phi\Theta} (X)/A^f
$$
determines a contravariant functor
$$
K_{\Phi\Theta} (f) \to \bbc_{\Phi\Theta} (f)
$$
for exact $K_{\Phi\Theta} (f)$, which follows from the proposition
6.  This functor is not a duality in general.  Consider now a
special case when it is a duality.

Let $\Theta$ be a variety of algebras and let $G$ be an algebra in
$\Theta$.  Consider the variety $\Theta(G)$ of algebras in
$\Theta$ with the constants from $G$.  For example, if $\Theta$ is
the variety of commutative and associative rings with the unit and
$P$ is a field, then $\Theta(P) = Var - P$ is the variety of
algebras over $P$.

In the general case we consider the $G$-algebra $G$ in $\Theta(G)$
and the equations in it.  Here, the free algebra is of the type
$W(X) = G* W_0 (X)$, where $W_0 (X)$ is free in $\Theta$, $*$ is
the free product in $\Theta$.

Let $\mu: W(X) \to G$ be a point, $X = \{ x_1, \dots, x_n\}$, and
 $\mu (x_1) = a_1, \dots, \mu (x_n)= a_n$. The constants $a_1,
\dots, a_n$ are considered as elements in $W(X)$.  Thus, we can
speak of a formula $u$ of the type
$$
(x_1 \equiv a_1)  \wedge \dots \wedge (x_n = a_n).
$$
In this case $\Val_f (u) $ consists only of the element $\mu$, and
the point $\mu$ is closed in the Zariski topology. Here we have
the duality of categories.

Note also that the point $\mu$ is a simple variety, determined by
elementary formula $u$.  Let us make some remarks, related to this
observation.

Let $u$ be a formula of the type $X$ and $v$ be an elementary
formula of the type $Y \supset X$, equivalent to $v$.
Let us associate algebraic varieties $\Val_f(u)$ and $\Val_f(v)$.
Proceed from the identical inclusion $s^0\colon W(X) \to W(Y)$ and
$v=s^0_* u$. We have:
$$
\Val_f (v) = \Val_f(s^0_* u) = s^0_* \Val_f(u).
$$
Here $\Val_f(v)$ is a cylinder in $\Hom(W(Y), G)$, over the
initial variety $\Val_f(u)$.  This new cylinder is determined by
the elementary formula.

Consider further one special characteristic of the objects in the
category $\bbc_{\Phi\Theta} (X)$.  These are Boolean algebras of
the form $\Hal_{\Phi\Theta} (X)/A^f$.

We have a canonical homomorphism
$$
\Val^X_f\colon \Hal_{\Phi\Theta} (X) \to \Bool (W(X), G)
$$
for every $X\in \Gamma^0$.  We have also an isomorphism of Boolean
algebras.
$$
\chi:\Bool(W(X), G) \to Fun(\Hom (W(X), G), Z).
$$
Here $Z=\{0, 1\}$ is a two-element Boolean algebra, $ Fun(\Hom
(W(X), G), Z)$ is a Boolean algebra of binary functions, and
$\chi(A)$ is the characteristic function of the set $A$
in $\Hom (W(X), G)$.

Consider further the composition of homomorphisms
$$
\sim = (\chi\Val_f^X)\colon\Hal_{\Phi\Theta}(X) \to Fun
(\Hom(W(X), G), Z).
$$
Here $\tilde u$ is a characteristic function of the algebraic set
$\Val_f(u)$, corresponding to a formula $u \in \Hal_{\Phi\Theta}
(X)$.

Take an algebraic set $A=T^f$ with $T=A^f$ for the given model
$(G, \Phi, F)$, and consider a mapping
$$
\psi_A\colon Fun(\Hom(W(X), G), Z) \to Fun (A, Z)
$$
which associates a restriction of the function on $A$ to every
function from the left side.  It is a homomorphism of Boolean
algebras.

Let $\tilde u_A =\psi_A(\tilde u)$.  The transition $u\to\tilde u$
is a homomorphism of Boolean algebras.  Call its image {\it an
algebra of regular functions on the set} $A$.

\proclaim{Proposition 8}  The Boolean algebra $\Hal_{\Phi\Theta}
(X)/T$ is isomorphic to the algebra of regular functions defined
on the variety $A$.
\endproclaim

{Proof.}  For every $u$ and $\mu\colon W\to G$ we have $\tilde
u(\mu) = 1$ if and only if $\mu \in \Val_f(u)$ or, the same, $u
\in \Log\Ker (\mu)$.  Now, $u \in T=A^f = \bigcap_{\mu \in A}
\Log\Ker (\mu)$ if and only if $\tilde\mu =1$  for every $\mu\in
A$.  This means that $u \in T$ if and only if $\tilde u_A$ is a
unit in the considered algebra of regular functions on the set
$A$.  This implies the proposition.

Further we need the following definition: \proclaim{Definition}
Two models $(G_1, \Phi, f_1)$ and $(G_2, \Phi, f_2)$ are
geometrically equivalent if and only if for the arbitrary set of
formulas $T$ of some type $X$ we have
$$
T^{f_1f_1} = T^{f_2f_2}.
$$
\endproclaim

\proclaim{Proposition 9}  If the models $(G_1, \Phi, f_1)$ and
$(G_2, \Phi, f_2)$ are geometrically equivalent then the (weak)
categories $K_{\Phi\Theta} (f_1)$ and $K_{\Phi\Theta}(f_2)$ are
isomorphic.
\endproclaim
{Proof.}  This proposition for the exact categories in the
equational theory  follows directly from the corresponding
duality. In our case there is no duality and we need a proof.

Let the given models $(G_1, \Phi, f_1)$ and $(G_2, \Phi, f_2)$  be
geometrically equivalent.  Build an isomorphism
$$
F\colon K_{\Phi\Theta} (f_1) \to K_{\Phi\Theta} (f_2).
$$
Define $F(A) = A^{f_1f_2} = B$ for every object $A$ in the
category $K_{\Phi\Theta} (f_1)$.  If $A$ is of the type $X$, then
$B$ is of the same type. Similarly, $B^{f_2f_1}$ for $B$.  Here,
$A^{f_1f_2f_2f_1} = A^{f_1f_1f_1f_1} = A^{f_1f_1} = A$. This means
that $F$ is a bijection on the objects. Let now $s\colon W(Y)\to
W(X)$ be given. Check that this $s$ is admissible for $A_1$ and
$A_2$ in $K_{\Phi\Theta} (F_1)$ if and only if this $s$ is
admissible for the corresponding $B_1 $ and $B_2$ in
$K_{\Phi\Theta}(f_2)$. Let $A_1\subset sA_2$ be given. Then
$A^{f_1}_1\supset (sA_2)^{f_1}$ and $A_1^{f_1f_2} \subset
(sA_2)^{f_1f_2}$.  We have:
$$
\eqalign{
&(sA_2)^{f_1}\supset s^*A^{f_1}_2,\cr
&(sA_2)^{f_1f_2}\subset (s^*A^{f_1}_2)^{f_2} = sA_2^{f_1f_2}.\cr}
$$
Hence,
$$
B_1\subset (sA_2)^{f_1f_2}\subset sB_2.
$$
This means that $s$ is admissible for $B_1$ and $B_2$.

The opposite direction can be checked similarly.
The proposition is proved.

Let us make a remark on the notion of geometrical equivalence of
models.
\proclaim{Proposition 10}  If two models are geometrically
equivalent, then they are elementary equivalent.
\endproclaim
{Proof.}  Let the models $(G_1, \Phi, f_1)$ and $(G_2, \Phi, f_2)$
be geometrically equivalent, and let the formula $u$ of some type
$X$ belongs to the elementary theory of the first model.  We have
$u \in\Hal_{\Phi\Theta}(X)$ and $\Val_{f_1}(u) = \Hom(W(X), G_1)$.

Let $A  = \Hom (W(X), G_1)$ and $T=A^{f_1}$.  Obviously, $T=T(X)$
is the elementary theory of the type $X$ of the first model and
$u\in T$.  The set $T$ is $f_1$-closed.  Moreover, it is the
minimal one with this property in $\Hal_{\Phi\Theta} (X)$.  By
assumption, $T$ is also $f_2$-closed and is the minimal one with
this property.  Hence, if $B=\Hom(W(X), G_2)$, then $B^{f_2} = T$.
 Therefore, $T$ is an elementary theory of the type $X$ of the
 second model.

 Since $u \in T$ is the formula, it belongs to the elementary
 theory of the second model.  Analogously, if the formula $u$
 belongs to the elementary theory of the second model, then $u$
 holds in the first model.

\proclaim{Problem 1}   Is the opposite true?, i.e., is it true
that elementary equivalent models are geometrically equivalent?
\endproclaim

Note in the conclusion that the category $K_\Theta(G)$ in the
equational theory is a full subcategory in $K_{\Phi\Theta}(f)$ for
the given model $(G, \Phi, f)$.  There is no such a relation for
$C_\Theta(G)$ and $C_{\Phi\Theta}(f)$.  The objects of the first
one are algebras in $\Theta$, while the objects of the second one
are Boolean algebras.

\subhead {5.5. The categories $K_{\Phi\Theta}$ and
$C_{\Phi\Theta}$}
\endsubhead


The objects of $K_{\Phi\Theta}$ are of the form
$$
(X, A, (G, \Phi, f)).
$$
Here, the model $(G, \Phi, f)$ is not fixed in the category,
and $A=T^f$ for some $T$ in $\Hal_{\Phi\Theta} (X)$.

The objects of $C_{\Phi\Theta}$  are of the form
$$
(\Hal_{\Phi\Theta} (X)\big/ T,  (G, \Phi, f)),
$$
where $T$ is $f$-closed Boolean filter in $\Hal_{\Phi\Theta} (X)$.

Let us pass to morphisms.  Given a homomorphism $\delta: G_1\to
G_2$ in $\Theta$, we have $\tilde \delta:\Hom (W(X), G_1) \to \Hom
(W(X), G_2)$ by the rule $\tilde\delta(\nu) = \delta\nu$.
According to the commutative diagram
$$
\CD
W(Y) @>s>> W(X)\\
@V\nu' VV @VV\nu V\\
G_2 @>\delta >> G_1
\endCD
$$
we can write $\nu'=\delta\nu s = \tilde \delta (\nu s) = \tilde
\delta \tilde s(\nu) = \tilde s \tilde \delta (\nu) = (s, \delta)
(\nu)$.

Let now two objects $(X, A, (G_1, \Phi, f_1))$ and $(Y, B, (G_2,
\Phi, f_2))$ be given.  We say that the pair $(s,\delta)$ is
admissible for $A$ and $B$ if $(s, \delta) (\nu) \in B$ for every
$\nu \in A$.  Such admissible $(s, \delta)$ determines morphisms
in the weak category $K_{\Phi\Theta}$.

For every set $B$ of the type $Y$ let us consider a set $(s,
\delta)B$ of the type $X$, determined by the rule:
$$
\nu: W(X) \to G_1 \in (s, \delta)B \; \; \hbox{\rm iff} \; \; (s,
\delta) (\nu) \in B.
$$
We have:  $\tilde s (\tilde \delta(\nu))\in B$ and $\tilde
\delta(\nu) \in sB$, $ \nu \in \delta s B$.  Here $sB$ is of the
type $X$.  If $C$ is a set of a type $Y$, then $\delta C$ is of
the same type, and $\nu \in \delta C$ if $\tilde \delta (\nu) =
\delta \nu \in G$.  Hence, $(s, \delta) B = \delta s B = s \delta
B$.  Besides, the pair $(s, \delta)$ is admissible  for $A$ and
$B$ if and only if $A \subset s \delta B$.  We write also $\delta
B = \delta_*B$ and consider $\delta^*$ determined by $\delta^* B =
\{ \delta \nu | \nu \in B\}$.

We can define also exact morphisms of the type $([s], \delta)$.
Given $\delta : G_1 \to G_2$, we have a mapping $[s]: A \to B$ by
the rule $[s] (\nu) = \delta \nu s, \nu \in A$.
As a morphism we take a pair $([s], \delta)$.  Morphisms
$$
\eqalign{
& (\ol s, \delta): (\Hal_{\Phi\Theta}(Y) /T_2, \; (G_2, \Phi,
f_2)) \rightarrow \cr
& (\Hal_{\Phi\Theta} (X) /T_1, (G_1, \Phi, f_1))\cr}
$$
are determined by the suitable pairs $(s, \delta)$.  Like $([s],
\delta)$, the homomorphism $\ol s$ of Boolean algebras depends on
$\delta: G_1 \to G_2$.  We will return to it further.

Consider now some details.

Given $A \subset \Hom (W(X), G_1)$, take $(\delta_* A)^{f_2} $ in
$\Hal_{\Phi\Theta} (X)$.
\proclaim{Proposition 11}  The pair $(s,
\delta)$ is admissible for the given $A$ and $B$ if and only if $s
u \in ( \delta^*A)^{f_2}$ for every $u \in B^{f_2}$.
\endproclaim
{Proof.}  Note first that for $(s, \delta)$ we have a homomorphism
$\ol s: \Hal_{\Phi\Theta} (Y)/B^{f_2} \to \Hal_{\Phi\Theta} (X)
/(s^* A)^{f_2}$ .  It is a morphism in the category
$C_{\Phi\Theta} (f_2)$. Let $(s,\delta)$ be admissible for $A$ and
$B$.  We claim that if $u\in T_2 = B^{f_2}$ then $s u \in
(\delta^* A)^{f_2}$.

The inclusion $u \in T_2 = B^{f_2}$ means that $B\subset
\Val_{f_2}(u)$.  Take an arbitrary $\nu $ in $A$.  Then $\nu' =
\delta \nu s \in B$ and $\delta \nu s \in \Val_{f_2} (u)$, $\delta
\nu \in \Val_{f_2} (su)$.  So, $\delta^* A\subset \Val_{f_2} (su),
s u \in (\delta^* A)^{f_2}$.

Let now $s u \in (\delta^* A)^{f_2}$ for every $u \in T_2$.  Then
$\nu' = \delta \nu s \in B$ for every $\nu \in A$.  This means
also that $\delta \nu \in sB = sT_2^{f_2} = (s^* T_2)^{fs} =
\mathop\cap\limits_{v\in s^* T_2} \Val_{f_2}(v)$. We should check that
$\delta \nu \in \Val_{f_2} (v)$ holds for every $v \in s^* T_2$.
We have $v = su, u \in T_2$.  Hence, $v = su \in (s^*A)^{f_2}$ and
$\delta^* A\subset \Val_{f_2} (v)$.  Thus, $\delta\nu \in
\delta^*A, \delta \nu \in \Val_{f_2} (v)$.  The proposition is
proved.

Consider once more the pairs $(s, \delta)$, admissible for the
fixed $A$ and $B$.  Fix also $\delta: G_1 \to G_2$ and vary the
component $s$.  For every such $s$ we have $[s]: A \to B$ and $\ol
s:\Hal_{\Phi\Theta}(Y)/B^{f_2} \to \Hal_{\Phi\Theta} (X)/(\delta^*
A)^{f_2}$.
\proclaim{Proposition 12}  For every given $\delta$ we
have $[s_1]=[s_2] \to \ol{s_1} = \ol{s_2}$.
\endproclaim
{Proof.} The proof of this proposition repeats one of the
proposition 6 and we omit it.

Further we return to morphisms in the category $C_{\Phi\Theta}$.
As we have seen, they have the form
$$
\eqalign{ &(\ol s, \delta):
(\Hal_{\Phi\Theta} (Y)/T_2, (G_2, \Phi, f_2))\rightarrow
\cr
&\rightarrow (\Hal_{\Phi\Theta} (X) /T_1, (G_1, \Phi, f_1)).\cr}
$$
The pairs $s$ and $\delta$ are correlated in a special way.  Here,
$\delta: G_1 \to G_2$ is a homomorphism of algebras in $\Theta$.
Take $$ (\delta^* T_1^{f1})^{f_2} = T^\delta_1 $$ and assume that
$s: W(Y) \to W(X)$ induces a homomorphism of Boolean algebras
$$
\ol s:
\Hal_{\Phi\Theta} (Y) /T_2 \to \Hal_{\Phi\Theta} (X) /T^\delta_1
$$ which is a morphism in the category $C_{\Phi\Theta} f_2$.
These  conditions mean correlation between $s$ and $\delta$. This,
in turn, implies that the pair $(s, \delta)$ is admissible for $A
= T^{f_1}_1$ and $B=T^{f_2}_2$.

Pay attention to the fact that the homomorphism $\ol s$, is not in
general a homomorphism of the initial Boolean algebras.

Let us consider multiplication of morphisms in $K_{\Phi\Theta}$
and $C_{\Phi\Theta}$ in more detail using the arguments above. Let
$$([s], \delta): (X, A; (G_1, \Phi, f_1)) \to (Y, B; (G_2, \Phi, f_2))
$$
and
$$
([s'], \delta'):(Y, B; (G_2, \Phi, f_2)) \to (Z, C; (G_3, \Phi, f_3))
$$
be given in the exact category $K_{\Phi\Theta}$.
Define the multiplication:
$$([s'], \delta')([s], \delta) = ([ss'], \delta' \delta):
(X, A; (G_1, \Phi, f_1)) \to (Z, C; (G_3, \Phi, f_3)).
$$
One should check that the homomorphisms $ss'$ and $\delta'\delta$
are correctly coordinated.  This means that $\delta'\delta\nu
ss'\in C$ holds for every $\nu \in A$.  We have $\delta \nu s \in
B$ and, further, $\delta'(\delta \nu s) s'\in C$.  Hence, $ss'(u)
\in T^{(\delta'\delta)}_1$ for every $u \in T_3$.  This gives the
homomorphism
$$
\ol{ss'} : \Hal_{\Phi\Theta} (Z)/T_3 \to\Hal_{\Phi\Theta} (X) /T_1^
{(\delta'\delta)}.
$$
Define the multiplication of morphisms in $C_{\Phi\Theta}$ as
$$
(\ol s', \delta') (\ol s, \delta) = (\ol{ss'}, \delta'\delta).
$$
We see that the transition $(X, A; (G, \Phi, f))\to
(\Hal_{\Phi\Theta}(X)/A^f, (G, \Phi,f))$ determines a
contravariant functor $K_{\Phi\Theta} \to C_{\Phi\Theta} $ for
weak  and exact category $K_{\Phi\Theta}$.

Prove  that the homomorphism $\ol{ss'}$ may be  represented as a
product of two morphisms in the category $C_{\Phi\Theta} (f_3)$
related to $s$ and $s'$.

We have a diagram for $s'$:
$$
\CD
Hal_{\Phi\Theta} (Z) @>s'>> \Hal_{\Phi\Theta}(Y)\\
@V\mu_Z VV @VV\mu_Y V\\
\Hal_{\Phi\Theta} (Z)/T_3 @>\ol{s'} >>
\Hal_{\Phi\Theta}(Y)/T^{\delta'}_2
\endCD
$$

\noindent
and another one for $s$

$$
\CD
\Hal_{\Phi\Theta}(Y) @>s>> \Hal_{\Phi\Theta}(X)\\
@V\mu_Y VV @VV\mu_X V\\
\Hal_{\Phi\Theta}(Y)/T_2 @>\ol{s} >>
\Hal_{\Phi\Theta}(X)/T^\delta_1
\endCD
$$

\noindent
The third diagram is built on the base of the second one:

$$
\CD
\Hal_{\Phi\Theta}(Y) @> s>> \Hal_{\Phi\Theta} (X)\\
@[1]V\mu'_Y VV @[1]V\mu'_X VV\\
\Hal_{\Phi\Theta} (Y) / T_2^{\delta'} @>\ol{s}^{\delta'}>>
\Hal_{\Phi\Theta} (X)/T_1^{(\delta'\delta)}
\endCD
$$

\noindent
with $\ol{ss'} = \ol s^{\delta'} \ol{s'}$.   Here
$\ol{s}^{\delta'}$ and $\ol{s'}$ are morphisms in $C_{\Phi\Theta}
(f_3)$.

We need to show that $s u \in T_1^{(\delta'\delta)}$ holds for
every $u \in T^{\delta'}_2$.  The filter $T^{\delta'}_2$ is
constructed with the use of the commutative diagram

$$ \CD
@.\!\! W(Y)\\
@[2]/SW/\nu'// @.@. @/SE//\mu/ \\
G_3 @[2]<\delta'<< \! G_2
\endCD
$$
Here $T^{\delta'}_2 = (\delta^{'*} B)^{f_3} \subset
\Hal_{\Phi\Theta} (Y)$.

Let us pass to $(T^{\delta'}_2)^{f_3} = (\delta^{'*} B)^{f_3f_3} =
C_1$.

Consider  the diagram
$$ \CD
@.@.W(X)\!\!\!\!\!\!\!\!\! \\
@.@[2]/SW/ \nu^{''}// @.@/SW// \nu'/ @VV\nu V \\
G_3 @<\delta'<< G_2
@<\delta << G_1
\endCD
$$
It determines $T^{(\delta'\delta)}_1 \subset \Hal_{\Phi\Theta}
(X)$. The inclusion $s u \in T^{(\delta'\delta)}_1$ for every $u
\in T^{\delta'}_2$ means that the pair $(s, \delta' \delta)$ is
admissible for $A$ and $C_1$.  Check the last statement.  Let $\nu
\in A$.  Then $\delta \nu s \in B$ and, hence, $\delta'\delta\nu s
\in \delta^{'*}B\subset (\delta^{'*}B)^{f_3f_3} = C_1$.

It is left to check the equality $\ol{ss'} = \ol s^{\delta'}
\ol{s'}$. Using commutative diagrams
$$
\CD
\Hal_{\Phi\Theta}(Z)
@>s'>> \Hal_{\Phi\Theta}(Y) @>s>> \Hal_{\Phi\Theta} (X) \\
@V\mu_ZVV @V\mu'_Y VV @VV\mu'_XV\\
\Hal_{\Phi\Theta}(Z)/T_3 @>\ol{s'} >>
 \Hal_{\Phi\Theta}(Y)/T^{\delta'}_2 @>\ol{s}^{\delta'}
 >> \Hal_{\Phi\Theta}(X)/T^{(\delta'\delta)}_1
 \endCD
 $$
and
$$
\CD
\Hal_{\Phi\Theta}(Z)
@>ss'>> \Hal_{\Phi\Theta}(X) \\
@V\mu_Z VV @V\mu'_X VV \\
\Hal_{\Phi\Theta}(Z)/T_3 @>\ol{ss'} >>
 \Hal_{\Phi\Theta}(X)/T^{\delta'\delta}_1
 \endCD
 $$
we can rewrite the product of morphisms as
$$
(s^{\ol{'}}, \delta') (\ol s, \delta) = (\ol s^{\delta'} \ol{s'},
\delta' \delta).
$$
Let us return to the definition of the morphisms in
$K_{\Phi\Theta}$.  Given the morphism
$$
(s, \delta): (X, A; (G_1, \Phi, f_2)) \to (Y, B; (G_2, \Phi,
f_2)),
$$
we have $\delta \nu s \in B$ for every $\nu \in A$. Pass to the
set $\delta^* A$, i.e., all $\delta \nu, \nu \in A$. Here,
$\delta^* A\subset sB$.  The set $sB $ is $f_2$-closed. Hence, if
$A_1 = (\delta^* A)^{f_2f_2}$, then $A_1 \subset sB$ and the
mapping $s:W(Y)\to W(X)$ is admissible for $f_2$-closed sets $A_1$
and $B$.  We have here a morphism
$$
[s]: (X,A_1) \to (Y,B)
$$
in the category $K_{\Phi\Theta} (f_2)$.  This morphism is
coordinated with the morphism
$$
\ol{s}:\Hal_{\Phi\Theta} (Y)/T_2 \to \Hal_{\Phi\Theta} (X)
/T^\delta_1
$$
in the category $C_{\Phi\Theta} (f_2)$. The definitions of
morphisms in the categories $C_{\Phi\Theta}$ and $K_{\Phi\Theta}$
are now  coordinated.

Note further that every category $K_{\Phi\Theta} (f)$ is a
subcategory in $K_{\Phi\Theta}$ and every $C_{\Phi\Theta}(f)$ is
also a subcategory in $C_{\Phi\Theta}$.  Besides, we have
subcategories $K_{\Phi\Theta}(G)$ in the category $K_{\Phi\Theta}$
for every given algebra $G \in \Theta$ and with arbitrary
interpretations $f$ of the set $\Phi$ in $G$.  There are
subcategories $C_{\Phi\Theta}(G) $ in the category
$C_{\Phi\Theta}$.

Some more details on these categories.  Objects of the category
$K_{\Phi\Theta}(G)$ can be represented in the form
$$
(X, A, f),
$$
since $G$ and $\Theta$ are fixed for this category.  Homomorphisms
$\delta: G \to G$ are identical homomorphisms.

Morphisms in $K_{\Phi\Theta}(G)$ are
$$
[s]: (X, A, f_1) \to (Y, B, f_2),
$$
where $s$ and the unit are coordinated  in accordance with the
general definition of the category $K_{\Phi\Theta}$.  Here, $A =
A^{f_1}_1, B = T_2^{f_2}$ for some $T_1$ in $\Hal_{\Phi\Theta}
(X)$ and $T_2=\Hal_{\Phi\Theta}(Y)$.  We may also assume that $T_1
=A^{f_1}$ and $T_2 = B^{f_2}$.  The corresponding set $A_1$ is now
$A^{f_2f_2}$ and $T^\delta_1$ is $T^{f_1f_2}_1 = A^{f_2}$.

Homomorphism $s:W(Y)\to W(X)$ induces a homomorphism of Boolean
algebras
$$
\ol{s}: \Hal_{\Phi\Theta} (Y)/B^{f_2} \to \Hal_{\Phi\Theta} (X)
/A^{f_2}.
$$
Every such $s$ determines a mapping $[s]:A\to B$ and,
simultaneously, a mapping $A_1 \to B$.

 As it was pointed out, the
subcategory $R_f$ of simple varieties (i.e., varieties, determined
by one-element sets $T$) is selected in the category
$K_{\Phi\Theta}(f)$ for a given model $(G, \Phi, f)$. This  $R_f$
is also a Halmos algebra (subalgebra in $\Hal_\Theta(G))$
coinciding with the image of the homomorphism $$ \Val_f:
\Hal_\Theta (\Phi) \to \Hal_\Theta(G).$$ We can take the
homomorphism
$$
\Val_f: \Hal_\Theta (\Phi) \to R_f
$$
and consider inclusions on the level of objects:
 $$ \Val_f: \Hal_\Theta (\Phi) \to
K_{\Phi\Theta}(f). $$ Every such inclusion determines transitions
from the sets of formulas $T$ to the corresponding objects of the
category $K_{\Phi\Theta}(f)$.

Consider a new category with the objects
$\Val_f:\Hal_{\Phi\Theta}(f) \to K_{\Phi\Theta}(f)$ and
commutative diagrams
$$
\CD
\Hal_\Theta(\Phi) @[2]>\Val_{f_1}>>  K_{\Phi\Theta}(f_1)\\
@[3]/SE/\Val_{f_2}// @.@.\!\!\!\!\! @VV\gamma V\\
@.@. K_{\Phi\Theta}(f_2)
\endCD
$$
as morphisms.

Here, $\gamma  $ is a functor of categories, and the commutativity
of the diagram  is understood on the level of objects.

Besides, $f_1$ is associated with the model $(G_1, \Phi, f_1)$,
and $f_2$ with $(G_2, \Phi, f_2)$. Thus, also the algebras $G$
vary.

The corresponding functor $\gamma$, subject to the condition of
commutativity of the diagram, we call a {\it special functor}.

The functor $\gamma $ is called {\it strong}, if $\gamma$ induces
the algebra homomorphism $R_{f_1} \to R_{f_2}$.

Every special $\gamma $ is strong.

If $(G_1, \Phi, f_1)$ and $(G_2, \Phi, f_2)$ are geometrically
equivalent, then the functor $\gamma=\gamma_{f_1}, f_2 = f$,
defined above, is a special isomorphism $K_{\Phi\Theta} (f_1) \to
K_{\Phi\Theta} (f_2)$.

Define now the semimorphisms in the considered category.  They are
of the form
$$
\CD
\Hal_{\Phi\Theta} @[2]>\Val_{f_1}>> K_{\Phi\Theta}(f_1)\\
@[1]V\sigma VV @. @[1]VV\gamma V \\
\Hal_{\Phi\Theta} @[2]>\Val_{f_2} >> K_{\Phi\Theta} (f_2)
\endCD
$$
where $\sigma$ is an automorphism of the algebra
$\Hal_{\Phi\Theta}  = \Hal_\Theta (\Phi)$.

In this case $\gamma$ induces some semihomomorphism of algebras
$R_{f_1} \to R_{f_2}$. The main problem is to find conditions on
the models which provide an isomorphism of the corresponding
categories of algebraic sets. We have already seen that if the
models are geometrically equivalent then this is the case.  We may
also claim that in this situation the models are elementary
equivalent and, in particular, if they are finite, then they are
isomorphic.

More general information we get from the Galois theory, to which
we pass right now.

\head 6. The Galois-Krasner theory in algebraic logic and
algebraic geometry
\endhead

\subhead{6.1. The group $Aut(\Hal_\Theta (G))$}
\endsubhead

The theory presented in  this and  next sections goes back to the
work of M. Krasner [Kr]. Later there were the papers [Be]; [Ba],
[Pl1], and others.  Here we follow the schemes from [Be1] and
[Pl1]. Some new details are related to peculiarities of Halmos
algebra $\Hal_\Theta (G)$.

Automorphism $\sigma$ of the Halmos algebra $\Hal_\Theta (G)$
determines  an automorphism $\sigma_X$ of the quantorian
$X$-algebra $\Bool (W(X), G)$ for every finite set $X \subset
X^0$.  All these automorphisms $\sigma_X $ are naturally
coordinated with the definition of $s\in S_\Theta$.

Show that some  automorphism $\delta_* =\sigma\in Aut
(\Hal_\Theta(G))$ corresponds to an automorphism $\delta \in Aut
G$.

Start from the construction used earlier.  Let a homomorphism
$\delta: G_1 \to G_2$ be given.
Then
$$
\tilde\delta : \Hom(W(X), G_1) \to \Hom(W(X), G_2).
$$
If now $B\subset \Hom(W(X), G_2)$, then
$$A=\delta_*B=\tilde\delta^{-1} (B)\subset \Hom(W(X), G_1),$$
that is $ \mu \in A = \delta_*B$ if and only if $\delta\mu \in B$.
Thus, we have
$$
\delta_*: \Hal_\Theta(G_2) \to\Hal_\Theta(G_1).
$$
This $\delta_*$ is  compatible with the Boolean structure and with
the morphisms:  $\delta_* (sB) = s\delta_* (B)$ for every $s:
W(Y)\to W(X)$.  Indeed, let $\nu $ lie in $\Hom(W(Y), G_1)$. Then
$\nu \in \delta_* (sB)$ if and only if $\delta\nu\in sB$, $(\delta
\nu) s \in B$, $(\delta \nu)s = \delta (\nu s) \in B$.
Since the multiplication of morphisms is associative, this takes
place if and only if $\nu s \in \delta_* B$ and $\nu \in s
\delta_* B$.  This gives also that if $s: W(Y) \to W(X)$ is
admissible for $A$ and $B$, then the same $s$ is admissible for
$\delta_* (A)  $ and $\delta_* (B)$. Indeed, if $A \subset s B$,
then $\delta_* (A) \subset \delta_* (sB) = s\delta_* (B)$.  But
$\delta_*$ is not compatible with quantifiers in general.
$\delta_*$ is an isomorphism of Halmos algebras if and only if
$\delta: G_1 \to G_2$ is an isomorphism in $\Theta$.  In this case
$s$ is admissible for $A$ and $B$ if and only if this $s$ is
admissible for $\delta_*(A)$ and $\delta_*(B)$.  Hence if $\delta
\in Aut (G)$, then $\delta_* = \sigma\in Aut (\Hal_\Theta (G))$.

We have the representation of groups
$$
Aut (G) \to Aut (\Hal_\Theta(G)).
$$
The main result now is the following:
\proclaim{Theorem 7}  The representation $Aut(G) \to
Aut(\Hal_\Theta (G))$ is an isomorphism of groups.
\endproclaim

{Proof.} We start the proof from the following statement.  Let $M$
be an arbitrary non-empty set, ${\frak M} = \sub  M$ is the
Boolean algebra of all subsets in $M$.  Then every automorphism of
the Boolean algebra $\frak M$ determines a substitution on the set
$M$, which induces this automorphism.  See, for example, [Pl1], p.
338.

Here, if $\tau$ is a substitution on the set $M$, then the
corresponding automorphism $\tau_*$ is defined by the rule
$$
a \in\tau_* A = \tau A \; \; \hbox{\rm iff} \; \; \tau a \in A, \;
A \subset M.
$$

Fix further a finite set $X$ and proceed from $M = M_X = \Hom
(W(X), G)$.  Every automorphism $\sigma$ of the algebra
$\Hal_\Theta (G)$ determines an automorphism $\sigma_X $ of the
Boolean algebra $\Bool (W(X), G)$, i.e., of the algebra $\frak M
$$_X$.
 This $\sigma_X $ is given by a substitution $\tau_X$ on the set
 $M_X = \Hom (W(X), G)$ depending on $\sigma$ and $X$.

 Denote by $\Sigma_X$ the group of all substitutions of the set $M_X$,
 by $\Sigma_G$ the group of all substitutions of the set $G$ and by
 $\tilde \Sigma_G = \Sigma^X_G$ its cartesian power.  Consider an
 inclusion
 $$
 {}^{\sim:} \tilde\Sigma_G \to \Sigma_X.
 $$
 This inclusion is defined as follows:

 if $\zeta \in \tilde\Sigma_G$ then a substitution $\tilde
 \zeta
 \in \Sigma_X $ is defined by the rule: for each $\mu:W(X)\to G$, there correspond
  a homomorphism
 $\tilde \zeta(\mu): W(X) \to G$  defined by: $\tilde \zeta(\mu) (X) =
\zeta(X) (\mu(X))$ for every $x \in X$.

 We would like to characterize the group $\tilde \Sigma_G$ as a
 subgroup in $\Sigma_X$.

Consider an auxiliary proposition.

 \proclaim{Definition}  The substitution $\tau\in \Sigma_X$ is called
 correct if
 $$
 \mu(x) = \nu(x) \Leftrightarrow (\tau\mu)(x) = (\tau\nu)(x)
 $$
 for every $\mu, \nu: W(X) \to G$ and every $x \in X$
\endproclaim

\proclaim{Proposition 13}  A substitution $\tau$ is correct if and
only if $\tau = \tilde \zeta$ for some $\zeta=\tilde\Sigma_G$.
\endproclaim
{Proof.}  Let $\tau = \tilde\zeta$ for some $\zeta=\tilde
\Sigma_G$.  Check the correctness of $\tau$. Take $x \in X$, $\mu,
\nu: W(X) \to G$.  Then $\tilde\zeta(\mu)(x) = \zeta (x) \mu (x);
\tilde \zeta(\nu)(x) = \zeta (x) \nu (x)$. It is clear that the
equality $\tilde \zeta (\mu) (x) = \tilde \zeta (\nu) (x)$ holds
if and only if $\mu (x) = \nu (x)$.

Let not $\tau$ be an arbitrary correct substitution on the set
$M_X = \Hom (W(X), G)$.  Build $\zeta \in \tilde \Sigma_G$ with
the condition $\tau =\tilde \zeta$.

Fix some $x \in X$ and define a substitution $\zeta(x)$ on the set
$G$.  Take an arbitrary $g \in G$, and $\mu: W(X) \to G$ with $\mu
(x) = g$.  Set:
$$ \zeta(x)g = (\tau\mu)(x).
$$
Since $\tau $ is correct, this definition does not depend on the
choice of $\mu$.  If $\nu (x) = \mu (x)$, then $(\tau\nu)(x) =
(\tau\mu) (x) = \zeta(x)(g)$.
Check that the mapping $\zeta(x): G\to G$ is a substitution.

Let us note first of all that if $\tau$ is a correct substitution,
then the same is valid for $\tau^{-1}$.  Indeed, take $x \in X,
\mu, \nu: W(X) \to G$.  Denote $\mu_1 = \tau^{-1}\mu,
\nu_1 = \tau^{-1}\nu, \mu = \tau\mu_1, \nu =
\tau\nu_1$.
We have: $(\tau\mu_1) (x) = (\tau\nu_1)(x) \Leftrightarrow
(\mu_1(x) = \nu_1(x)$.  This gives $\mu(x) \nu(x) \Leftrightarrow
(\tau^{-1} \mu)(x) = (\tau^{-1}\nu)(x)$.  The substitution
$\tau^{-1} $ is correct.

Take now an arbitrary element $g_1\in G$ and show that $\zeta(x)
(g) = g_1$ for some $g \in G$.  Let $\mu_1(x) = g_1, (\tau^{-1}
\mu_1)(x) = g = \mu (x), \mu = \tau^{-1} \mu_1$. Then $(\tau\mu)
(x) = \zeta(x)(g) = (\tau\tau^{-1} \mu_1) (x) = \mu_1 (x) = g_1$.

It remains to verify that $\zeta (x) (g_1) = \zeta (x) (g_2)$ implies
$g_1 = g_2$. Let $\mu_1(x) = g_1$,  $\mu_2 (x) = g_2$.
We have $\zeta (x) (g_1) = (\tau\mu_1)(x) = \zeta(x) (g_2) = (\tau
\mu_2)(x)$.  The last gives $\mu_1 (x) = \mu_2(x) $ and $g_1 =
g_2$.  We checked that $\zeta (x)$ is a substitution.

Take $\zeta \in \tilde \Sigma_G = \Sigma^X_G$.  According to the
rule: $\zeta (x)$ is determined by $\tau$.  Verify that $\tau =
\tilde \zeta$.  Take an arbitrary $\mu: W(X) \to G$ and check that
$\tau \mu = \tilde \zeta\mu$.  Take an arbitrary $x \in X$.  Then
$(\tau \mu) (x) = \zeta(x) \mu(x) = (\tilde \zeta \mu)(x)$.
This finishes the proof of the proposition.

\proclaim{Proposition 14} A substitution $\tau$ commutes with the
quantifiers in $\Bool (W(X),G)$ if and only if $\tau$ is correct.
\endproclaim
{Proof.}  We mean commutativity of the form
$$
\tau\exists x A = \exists x \tau A
$$
for $x \in X$ and $A\subset \Hom(W(X),G)$.

Let $\tau$ be correct, $\tau = \tilde \zeta, \zeta \in\tilde
\Sigma_G$. Take $\mu \in\tau\exists x A = \tilde \zeta \exists x
A$.  We have $\tilde \zeta \mu \in\exists x A$.  Choose $\nu \in
A$ with $(\tilde \zeta \mu) (x') = \nu (x')$ for every $x'\in X,
x'\neq x$.  Let $\nu = \tilde \zeta \nu_1$.  Since $\nu \in A$,
then $\nu_1 \in \tilde\zeta A$.  The equality $\tilde \zeta \mu
(x') = \tilde\zeta \nu_1 (x')$ gives $\mu(x')= \nu_1(x')$ for
every $x' \neq x$, and thus $ \mu \in \exists x \tilde\zeta A$.

Let now $\mu \in \exists x \tilde \zeta A$.  Take $\nu_1$ in
$\tilde \zeta A$ with $\mu (x') = \nu_1 (x')$ for $x' \neq x$.  We
have $(\tilde \zeta\mu)(x') = (\tilde\zeta \nu_1) (x')$.  Here,
$\nu = \tilde \zeta \nu_1 \in A$.  Hence, $\tilde \zeta \mu \in
\exists x A$, $ \mu \in \tilde \zeta \exists x A$.  We have
checked commutativity for $\tau = \tilde \zeta$.

Let us check the opposite direction, i.e., the commutativity
fulfills and we verify the correctness of the substitution $\tau$.

Let $\mu$ and $\nu$ be two elements in $\Hom (W(X), G)$, and
$\mu(x) = \nu(x)$ for some $x \in X$.  We need to check that
$(\tau\mu)(x) = (\tau\nu)(x)$ holds true.  Denote by $Y$ the set
$X$ without $x $ (i.e., $Y=X\setminus x)$.
Denote the product of all $\exists y, y \in Y$ by $\exists (Y)$.
Commutativity of $\tau$ with each $\exists y$ implies
commutativity with $\exists (Y)$.  Take, further, the set $A$
consisting of one element $\nu$, $A = \{
\nu \}$.  Then $\mu \in \exists (Y) \{ \nu\}$.  Now,
$$
\tau \mu \in \tau^{-1} \exists (Y) \{ \nu\} = \exists (Y)
\tau^{-1} \{ \nu\} = \exists (Y) \{ \tau\nu\}.
$$
Hence, $(\tau\mu)(x) = (\tau\nu)(x)$.

Here we used commutativity of $\tau^{-1}$ with quantifiers.
Similarly we derive that \break $(\tau\mu) (x) = (\tau\nu)(x)$ implies
$\mu(x) = \nu(x)$.  Thus, $\tau$ is a correct substitution.

Hence, one can  assume that every automorphism of a quantifier
$X$-algebra $\Bool (W(X), G)$ is realized by a correct
substitution of the set $\Hom (W(X), G)$, depending on the given
automorphism.

Consider further coordination with  the operations $s\in S_\Theta$.

Pick out a subgroup $\Sigma_0$ in $\tilde\Sigma_G = \Sigma^X_G$,
consisting of constant $\zeta$, i.e., $\zeta (x) = \delta \in
\Sigma_G$ for all $x \in X$.  A constant $\zeta$ corresponds to  a
$\delta \in \Sigma_G$.  This determines the inclusion $\Sigma_G
\to \tilde\Sigma_G = \Sigma^X_G$ whose image is $\Sigma_0$.
Simultaneously we have an inclusion $\Aut(G) \to \Sigma^X_G$ for
each $X \subset X^0$.

Denote the pointed inclusion for every $X$ by $^{\sim} X$.  If
$\delta \in \Aut (G)$ or $\delta \in \Sigma_G$, then
$\delta^{^\sim X} \in \tilde \Sigma_G$.  Further we vary the
finite set $X \subset X^0$ and pass to the multi-sorted variant of
sets and algebras.  In particular, we consider multi-sorted set
$M$ with the components $M_X = \Hom (W(X), G)$.

The substitution $\tau$ of $M$ is a function, determining a
substitution $\tau^X$ of the set $M_X$ for every $X$.  Such $\tau$
is correct if all $\tau^X$ are correct, $\tau^X=\zeta^{^\sim X},
\tau = \tilde\zeta$.  Here $\zeta$ is also a function, determining
an element $\zeta^X$ in $\Sigma^X_G$ for every $X$.  To every
$\delta\in\Sigma_G$ corresponds $\delta^{^\sim X}$ in
$\Sigma^X_G$, which determines a substitution $\tilde\delta$ of
the multi-sorted set $M$.

Consider homomorphisms $s\colon W(X) \to W(Y)$.  An operation $s:
\Bool (W(X), G)$ $\to $ $\Bool (W(Y), G)$ in the algebra
$\Hal_\Theta(G)$ corresponds to a homomorphism $s$.  The
substitution $\tau$ of a multi-sorted set $M$ commutes with all
such $s$ if for every $A$ in $\Hom(W(X), G)$ holds $\tau^Y
sA=s\tau^X A$.

Let $\tau=\tilde \zeta$ be an arbitrary substitution.
\proclaim{Proposition 15}  The substitution $\tau$ commutes with
all $s \in S_\Theta$ if and only if $\tau = \tilde\delta$ where
$\delta$ is an automorphism of the algebra $G$.
\endproclaim
{Proof.}  We already know that every $\delta \in \Aut (G)$ acts in
the algebra $\Hal_\Theta (G)$ as an automorphism.  In particular,
this means that the function $\tilde \delta = \tau $ for such
$\delta$ has the needed properties.

Let now $\tau = \tilde\zeta$ commute with operations ot the type $s$.
We want to
check if $\tau = \tilde \delta$, where $\delta \in \Aut (G)$.

Show first that every $\zeta^X$ is a constant, and then notice
that all $\zeta^X$ for different $X \subset X^0$ get the same
value $\delta \in \Sigma_G$.  We will do it simultaneously.  Take
$x \in X$ and $y \in Y$ and assume that $\zeta^X (x) \neq
\zeta^Y(y)$.  Select an element $a \in G$ with $\zeta^X (x) a\neq
\zeta^Y(y) a$, and take $\mu: W(Y) \to G$ with $\mu(y) = a$.  Let
$s: W(X) \to W(Y)$ transfer $x$ to $y$.  Then
$$ \eqalign{
&(\tau^Y\mu)(sx) = (\tau^Y\mu)(y) = \zeta^Y(y) \mu(y) =\cr
&=\zeta^Y(y) (a),\cr &\zeta^X(x) (a) = \zeta^X(x) \mu(y) =
\zeta^X(x) (\mu s) (x) =\cr &=(\zeta^X(x) (\mu s)) (x),\cr
&(\zeta^X(x)(\mu s)) (x) = \tau^X (\mu s)(x).\cr}
$$
Hence,
$$
\eqalign{
&((\tau^Y\mu)s)(x)\neq \tau^X (\mu s)(x),\cr
&(\tau^Y \mu) (s) \neq \tau^X (\mu s).\cr}
$$
We can rewrite it as
$$
\tilde s(\tau^Y \mu) \neq \tau^X(\tilde s (\mu)).
$$

Take now an arbitrary $A \subset \Hom (W(X), G)$, and let $\mu
\in s \tau^X A$.  Thus, $\mu s = \tilde s (\mu) \in \tau^X A$,
$\tau^X (\mu s) \in A$, $\tau^X(\tilde s (\mu)) \in A$.

On the other hand, $\mu \in\tau^Y(sA)$ means that $(\tau^Y)(\mu)
\in s A$, $\tilde s (\tau^Y \mu) \in A$.  Since $\tau^X(\tilde
s(\mu)) \neq \tilde s (\tau^Y\mu)$, it is clear that for some $A$
there is no commutativity.  It suffices to take one-element set
$A$.  Thus, if $\tau$ is not constant, then there is no
commutativity.  Hence, $\tau$ is constant and $\tau = \tilde
\delta, \delta \in \Sigma_G$.

Prove further that the substitution $\delta$ is an automorphism of
the algebra $G$ if the commutativity takes place.

Let $\omega \in \Omega$ be an $n$-ary operation and $a_1, \dots, a_n$
elements in $G$.  Check that
$$
\delta(a_1\dots a_n \omega) = (\delta a_1) \dots (\delta a_n)\omega.
$$
Take $X=\{ x, x_1, \dots, x_n\}$, $W=W(X)$ and $\mu: W\to G$ with
the condition $\mu(x_1) = a_1, \dots, \mu(x_n) = a_n$.  Take also
$s: W\to W$ with $sx = x_1 \dots x_n \omega$.   The commutativity with
$s$ means for $\tau = \tilde \delta$ that
$$
\tau (\mu s) = \tau(\mu)\cdot s.
$$
For the given $x$ we have
$$
\eqalign{
&\tau(\mu s)(x) = \tilde \delta (\mu s ) ( x) = \delta(\mu s )(x)
= \cr
&=\delta(\mu s (x)) = \delta (\mu (x_1\dots x_n \om)) =
\delta(a_1\dots a_n\om),\cr
&(\tau(\mu)\cdot s)(x) = \tau(\mu)(sx) = \tau(\mu)(x_1\dots
x_n)\om=\cr
&(\tau\mu(x_1)) \dots (\tau\mu(x_n))\om = \delta (a_1) \dots
\delta (a_n)\om.\cr}
$$
This finishes the proof of Proposition 14.  Return to the proof of
Theorem 7.

Let $\sigma$ be an automorphism of the algebra $\Hal_\Theta (G)$.
This $\sigma$ induces an automorphism of every Boolean algebra
$\Bool (W(X), G)$.  Thus, $\sigma$ is determined by a substitution
$\tau$ of the multi-sorted set $M$ with the components $M(X) =
\Hom (W(X), G)$.  The commutativity with quantifiers means that
$\tau = \tilde \zeta$, while commutativity with operations of $s
\in S_\Theta$ type leads to $\tau = \tilde \delta$ where $\delta$
is an automorphism of the algebra $G$.  It is easy to see that
$\delta_*= \sigma$ for such $\delta$.  Theorem 7 is proved.
\subhead{6.2 The main results}
\endsubhead

Let a model $(G, \Phi, f)$ be given.  A group of automorphisms
$\Aut (f)$ of this model is a subgroup in a group $\Aut (G)$.  As
we have seen, $\Aut (G)$ acts in each set $\Hom(W(X), G)$.  A
group $\Aut (f)$ also acts there. \proclaim{Theorem 8}  Every
$f$-algebraic set in $\Hom (W(X), G)$ is invariant under the
action of the group $\Aut (f)$.
\endproclaim
{Proof.}  Check first that for every formula $ u \in
\Hal_{\Phi\Theta}(X)$, every point $\mu: W\to G$ and every
automorphism $\delta \in \Aut  (f)$ we have $$ \mu \in \Val_f(u)
\Leftrightarrow \delta \mu \in \Val_f(u).
$$
Let us call this property of  the formula $u$ the correctness
property. We want to show that every $u$ is correct.

Although this is trivial, we give the precise proof. All atomic
formulas are correct, and it is easy to see that if $u$ and $v$
are correct, then $u \vee  v$, $u \wedge  v$, $\lnot u$ are
correct as well.  Let $u$ be correct.  Take the formula $\exists x
u $ for some $x \in X$. Let $\mu \in \Val_f (\exists x u ) =
\exists x \Val_f(u)$.  We can take a point $\nu: W\to G$ such that
$\nu \in \Val_f (u), \mu (y) = \nu (y)$ for every $y \in X, y \neq
x$.  Take now $\delta \in \Aut (f)$.  We have:  $\sigma \mu (y) =
\delta \nu (y), \delta \nu \in \Val_f(u)$. So, $\delta \mu \in
\exists x \Val_f(u) = \Val_f (\exists x u)$.

Let, conversely, $\delta \mu \in \Val_f(\exists x u)$.
Take some $\nu = \delta \nu_1 \in \Val_f(u)$ with $\nu(y) = \delta
\nu_1 (y) = \delta \mu(y), y \neq x$.  Then $\mu (y) =
\nu_1(y),\nu_1 \in \Val_f(u)$ and $\mu \in \exists x \Val_f(u) =
\Val_f (\exists x u)$.

Finally, let $v=su$, where $u$ of some type $Y$ is correct.  Given
$s: W(Y) \to W(X)$, we get $\Val_f(v) = s\Val (u)$.  Take $\delta
\in Aut (f)$ and let $\mu \in \Val_f(v), \mu \in s \Val_f (u), \mu
s \in \Val_f (u)$.  By the condition on $u$ we have $\delta \mu s
\in \Val_f(u), \delta \mu \in s \Val_f (u) = \Val_f (su) = \Val_f
(v)$.

Let now $\delta \mu \in \Val_f (v) = s \Val_f (u)$ be given.  Then
$\delta \mu s \in \Val_f (u)$ and $\mu s \in \Val_f(u), \mu \in
\Val_f (v)$.

So, all correct formulas  form a subalgebra in $\Hal_\Theta
(\Phi)$.   This subalgebra contains the set of generators $M$, and
thus, all formulas in $\Hal_\Theta (\Phi)$ are correct.

Let now $A = T^f=\mathop\cap_{u\in T} \Val_f (u)$ for some $T$ of
the type $X$, $\mu \in A, \delta \in \Aut (f)$.  Then for every $u
\in T$ we have $\mu \in \Val_f (u)$, $\delta\mu \in \Val_f (u),
\delta \mu \in A$.

The theorem is proved.

Note in addition that if the set $A$ is $f$-closed and $\delta \in
\Aut (f)$, then $\delta_*A = \delta A = A$.  Indeed, we may claim
that if $\delta_* A = A$ for some $\delta \in \Aut (G)$ and $A
\subset \Hom (W(X), G)$, then there is a model $(G, \Phi, f)$ with
some $\Phi$ and $f$, such that $\delta \in \Aut (f) $ and the set
$A$ is $f$-closed.

Let us prove this claim.

Take $X = \{ x_1, \dots, x_n\}$.  Let the set $\Phi $ consist of
one $n$-ary relation $\varphi$.  Consider a formula $\varphi (x_1,
\dots, x_n)$, and realize $\vp$ in $G$.  We proceed from the
standard bijection $\pi: \Hom (W(X), G) \to G^{(n)}$, setting
$f(\vp) = A^\pi \subset G^{(n)}$.  We have $\Val_f (\vp (x_1,
\dots, x_n)) = A$ for such $f$.  Thus, $A$ is an $f$-algebraic
set. The condition $\delta A = A$ means that the automorphism
$\delta$ is coordinated with the interpretation $f$ of the given
set $\Phi$, i.e., $\delta \in \Aut(f)$.

Let us make some remarks on the constructions above.  Given a
model $(G_1, \Phi, f)$ and an algebra isomorphism $\delta: G_1 \to
G_2$, consider another model $(G_2, \Phi, f^\delta)$ with the
commutative diagram
$$
\CD
@.\!\!\Hal_\Theta (\Phi)\\
@[2]/SW/\Val_{f^\delta}//@.@. @/SE//\Val_f/ \\
\Hal_\Theta(G_2) @[2]> \delta_* >> \Hal_\Theta(G_1)
\endCD
$$
We determine $f^\delta$ such that the models are isomorphic and
the diagram is indeed commutative.

Take a formula $u$ of the type $X$ in $\Hal_\Theta (\Phi)$.  It
should be:
$$
\Val_f (u) = \delta_* \Val_{f^\delta} (u).
$$
Let $\vp$ be an $n$-ary relation in $\Phi$.  Take $X = \{ x_1,
\dots, x_n\}$ and $u = u(x_1, \dots, x_n)$. Pass to $\delta^{(n)}:
G^{(n)}_1 \to G^{(n)}_2$ and to bijections $\pi_1: \Hom (W(X),
G_1) \to G^{(n)}_1$ and $\pi_2: \Hom(W(X), G_2) \to G^{(n)}_2$.
The diagram
$$
\CD
\Hom(W(X), G_1) @>\tilde \delta>>\Hom (W(X), G_2)\\
@V\pi_1 VV  @VV\pi_2 V\\
G^{(n)}_1 @>\delta^{(n)} >> G^{(n)}_2
\endCD
$$
is commutative.

We should determine $f^\delta (\vp)$, solving the equation $\Val_f
(u) = \delta_* \Val_{f^\delta} (u)$ for $u = \vp (x_1, \dots,
x_n)$.  Denote $A = \Val_f (u)$ and $B = \Val_{f^\a} (u)$. Then
$A=\delta_*B$.  Here, $f(\vp) = \pi^*_1 A$ and $f^\delta (\vp) =
\pi^*_2 B$.  We have:  $\nu \in A$ if and only if $\delta \nu \in
B$, which gives $(a_1, \dots, a_n) \in f (\vp) $ if and only if
$(\delta a_1, \dots, \delta a_n) \in f^\delta (\vp)$.  The last
determines the model $(G_2, \Phi, f^\delta)$, isomorphic to the
initial model $(G_1, \Phi, f)$ by the isomorphism $\delta$.

It can be checked that the second model realizes the commutative
diagram.

Let us relate the algebraic sets for the models $f$ and
$f^\delta$. Let $A = T^f, B=T^{f^\delta}$ for some $T\subset
\Hal_{\Phi\Theta} (X)$. Then $A=\delta_* B$.  Indeed,
$$
\eqalign{
&A=T^f=\mathop\cap\limits_{u\in T} \Val_f (u) =\cap_{u\in T}
\delta_*\Val_{f^\delta} (u) =\cr
&=\delta_* (\mathop\cap\limits_{u\in T} \Val_{f^\delta} (u)) = \delta_*
T^{f^\delta} = \delta_* B.\cr}
$$
Note that the isomorphic models are geometrically equivalent.
This follows, for example, from the proposition 3.  Hence, the
models $(G_1, \Phi, f)$ and $(G_2, \Phi, f^\delta)$ are
geometrically equivalent.  If $T=A^f$ then the filter $T$ is
simultaneously $f$ and $f^\delta$-closed.  For such $T$ we have
$$
T=A^f = (\delta_*B)^f= B^{f^\a} = T, (\delta_*B)^f = B^{f^\delta}.
$$




The theory we consider here is quite natural to call
Galois-Krasner theory.

The algebra $\Hal_\Th(G)$ is considered together with its group of
automorphisms, represented as $\Aut G.$ For every
$A\in\Hal_\Th(G)$ and $\sigma\in\Aut G$, we write $\sigma A$
instead of $\sigma_*A.$ We give further standard definitions.

Let $H$ be a subset in $\Aut G.$ Then $R=H'$ is a subalgebra in
$\Hal_\Th(G),$ consisting of all elements $A$, for which $\sigma
A=A$ for every $\sigma \in H.$ Here $R=H'$ is actually a
subalgebra of a multi-sorted algebra $\Hal_\Th(G)$ and for every
$X\in\G^0$ the corresponding component $R(X)$ is a subalgebra of
the quantifier $X$-algebra $\Bool(W(X),G).$ The algebra $R(X)$
contains all equalities of the type $X.$

Let now $R$ be a subset (multi-sorted) of $\Hal_\Th(G).$ Then
$R'=A$ is a subgroup in $\Aut(G)$, consisting of all $\sigma\in
\Aut(G)$ for which $\sigma A=A$ for every $A\in R.$ It is actually
a subgroup in $\Aut G.$ We have Galois correspondence and Galois
closure $H''=(H')'$ and $R''=(R')'.$ A subgroup $H$ is closed if
$H''=H.$ Every closed subalgebra $R$ contains all equalities in
$\Hal_\Th(G).$ Recall that the equality in $\Hal_\Th(G)$ is an
algebraic set determined by an equation of the $w\equiv w'$ type.

\proclaim{Theorem 9}
For every model $(G,\Phi,f)$ we have
$$R_f'=\Aut(f).$$\endproclaim

{Proof.} It follows from Theorem $8$ that the inclusion
$\Aut(f)\subset R_f'$ takes place.

Let now $\sigma\in R_f'.$ It means that $\sigma A=A$ for every
$A\in R_f'.$ We need to show that the automorphism $\sigma$ of the
algebra $G$ is compatible with the interpretation $f(\vp)$ of the
arbitrary $\vp\in\Phi.$

Let the relation $\vp$ be $n$-ary. Take $X=\{x_1,\dots,x_n\}$ and
consider the formula $\vp(x_1,\dots,x_n).$ Take
$A=\Val_f(\vp(x_1,\dots,x_n)).$ For the bijection
$\pi:\Hom(W(X),G)$ $\to$ $G^{(n)}$ we have $A=\pi^{-1}(f(\vp)).$
Coordination of $\sigma$ with $f(\vp)$ is equal to the condition
$\sigma A=A.$ This condition holds true, since $A\in R_f.$ Thus,
$\sigma \in\Aut(f)$ and $R_f'=\Aut (f).$

We cannot claim that $(\Aut(f))'=R_f$ holds.
Indeed, in the general case we can consider an algebraic set $A$ which is
not simple, is not
determined by a single formula $u$ and, hence, does not belong to the
algebra $R_f.$
However, $\sigma A=A$ for every $\sigma \in\Aut(f)$ by Theorem 8.
Thus, $R_f\subset\Aut(f)'$ and the inclusion may be proper.

Let us formulate the first main theorem.
\proclaim{Theorem 10}
If the algebra $G$ is finite, then every subgroup in $\Aut(G)$ is closed
and every subalgebra
$R$ in $\Hal_\Th(G),$ containing all equalities, is also closed.\endproclaim

We have here a bijection between such $R$ and subgroups in $\Aut(G).$

We give the proof of this theorem  in the next subsection.

\proclaim{Theorem 11}
If the algebra $G$ is finite, then for every model $(G,\Phi,f)$ we have
$$\Aut(f)'=R_f.$$\endproclaim

{Proof.} From Theorems 9 and 10 follows that
$$R_f''=R_f=(R_f')'=\Aut(f)'.$$

It follows from Theorem 11 that if an algebra $G$ is finite, then
every $f$-algebraic set $A$ for the model $(G,\Phi,f)$ is a simple
algebraic set.
Now we can state that algebraic varieties  possesses the property
to be noetherian . This can also be proved directly if the algebra
$G$ is finite.

Consider the second main Galois theorem.

Take an isomorphism $\dl: G_2\to G_1$ of algebras in $\Th.$ An
isomorphism of Halmos algebras
$$\dl_*: \Hal_\Th(G_1)\to \Hal_\Th(G_2)$$
corresponds to $\dl.$

The groups $\Aut(G_1)$ and $\Aut(G_2)$ are isomorphic.
Their isomorphism is determined by $\sigma\to\dl\sigma\dl^{-1}$,\
$\sigma\in\Aut(G_2).$
This gives $\Aut(G_2)\to\Aut(G_1).$

Let now $R_1$ be a subalgebra in $\Hal_\Th(G).$
We have $\dl_*(R_1)=R_2$ in $\Hal_\Th(G).$
Hence, $R_1$ and $R_2$ are isomorphic.
We check directly
$$R_1'=\dl R_2'\dl^{-1}(\dl_*(R_1))'=\dl^{-1}R_1'\dl.$$
In this sense, the subgroups $R_1'$ in $\Aut(G_1)$ and $R_2'$ in
$\Aut(G_2)$ are conjugated by the isomorphism $\dl: G_2\to G_1.$

\proclaim{Theorem 12} If algebras $G_1$ and $G_2$ are finite, $R_1$ is a
subalgebra in
$\Hal_\Th(G_1)$ and $R_2$ is a subalgebra in $\Hal_\Th(G_2),$ then every
isomorphism $\g:
R_1\to R_2$ is induced by an isomorphism $\dl: G_2\to G_1.$
We have $\g=\dl_*: R_1\to R_2.$\endproclaim

\subhead{6.3. Proofs}
\endsubhead

We assume  that a finite algebra $G$ is considered also as a
subset in $X^0.$ A free in $\Th$     algebra $W(G)$ is an object
of the category $\Th^0.$ Consider an affine space $\Hom(W(G),G)$
and a quantifier $G$-algebra $\Bool(W(G),G)$ as an object of the
Halmos category $\Hal_\Th(G).$

Let us  make remarks on the set $\Hom(W(G),G)$ and the algebra
$\Bool(W(G),G).$

We have a standard bijection
$$\pi: \Hom(W(G),G)\to G^G,$$
where $\pi(\nu): G\to G$ is a restriction of the homomorphism $\nu$ on the
set $G.$

The right part, $G^G,$ is a semigroup of transformations of the
set $G.$ Let us extend the multiplication in this semigroup to the
left part. If $\nu_1,$ $\nu_2$ are two elements in $\Hom(W(G),G)$,
then $\nu_1\nu_2: W(G)\to G$ is a homomorphism, determined by the
condition $\pi(\nu_1\nu_2)=\pi(\nu_1)\cdot\pi(\nu_2).$

Now, $\Hom(W(G),G)$ is a semigroup, isomorphic to the semigroup $G^G$ with
the isomorphism
$\pi.$

The group $\Sigma_G$ is a group of all substitutions of the set $G.$
Now it may be represented as a group of all invertible elements in the
semigroup
$\Hom(W(G),G).$
It is a subset in $\Hom(W(G),G)$ and an element in $\Bool(W(G),G).$

Let $X$ be an arbitrary finite set in $X^0.$ Consider a mapping
$\mu: X\to G.$ This $\mu$ is simultaneously a mapping $\mu': X\to
W(G).$ We have homomorphisms $\overline\mu: W(X)\to G$ and $s:
W(X)\to W(G)$ determined by $\mu$ and $\mu', $ respectively. Since
$s$ is, in fact, determined by the mapping $\mu: X\to G,$  we
write $s=s(\mu).$

Apply these notations to the case when $X$ coincides with $G.$
Consider two mappings $\sigma_1$ and $\sigma_2:G\to G.$
We have $\overline{\sigma_1\sigma_2}=\overline\sigma_1\overline \sigma_2.$
We are interested in the product $\overline\sigma_1\cdot s(\sigma_2).$
It is easy to see that $\overline \sigma_1\cdot s(\sigma_2)=\overline{\sigma_1\sigma_2}.$

Take further an arbitrary subalgebra $R$ in $\Hal_\Th(G),$ containing
equalities.
This means that  $R$ preserves operations of the type $s^*,$
conjugated to the
operations of the $s_*$ type (compare [Pl1]).

We have $R(X)\subset\Bool(W(X),G)$ for every finite $X\subset X^0.$
Since the algebra $R(X)$ is finite, it has atoms.
Every element from $R(X)$ is a sum of atoms and all such atoms generate the
algebra $R.$

Atoms are easily built. Let $\mu: W(X)\to G$ be a point in the
space $\Hom(W(X),$ $G).$ This space is a unit in the algebra
$\Bool(W(X),G)$ and, simultaneously, a unit in $R(X).$ It contains
the point $\mu.$

Denote by $(\mu)$ the intersection of all elements of the algebra $R(X),$
containing the point
$\mu.$ This $(\mu)$ is an atom in $R(X),$ and all atoms in $R(X)$ are built
in such a way.

Consider further an algebra $R(G).$ It is a subalgebra in
$\Bool(W(G),G).$ Take a unit $e$ of the semigroup $G^G.$ A point
$\overline e: W(G)\to G$ which is a unit of the semigroup
$\Hom(W(G),G)$ corresponds to this $e.$ Take an atom $(\overline
e)$ in $R(G)$ determined by $\overline e.$ Denote $(\overline
e)=H=H^R.$

\proclaim{Proposition 16} For every algebra $R$ the subset $H$ in
$\Hom(W(G),G)$ is a subgroup in the group $\Sigma_G.$
\endproclaim
{Proof.} Let us make one general remark. Let a morphism $s:
W(X)\to W(Y)$ be given. If $A\in R(X),$ then $s_*A\subset R(Y)$
for the given algebra $R.$ Besides, we know that if $A\in R(Y),$
then $s^*A\in R(X).$

Take now an arbitrary $\sigma: G\to G.$ We have $\overline\sigma:
W(G)\to G$ and $s=s(\sigma): W(G)\to W(G).$ Take, further,
$s^*:\Bool(W(G),G)\to\Bool(W(G),G).$ By the construction, the atom
$H$ is an element in the algebra $\Bool(W(G),G).$ Take $s^*H.$ It
is also an element in $\Bool(W(G),G),$ consisting of all $\nu s,$
\ $\nu\in H.$ An element $\nu=\overline e$ is contained in $H.$
Hence, $\overline e s(\sigma)=\overline{e\sigma}=\overline
\sigma\in s^*H.$

Note also that both $H$ and $s^*H$ belong to the algebra $R(G).$
Assume that $\overline\sigma\in H.$
Then $\overline\sigma\in H\cap s^*H.$
This intersection is a nonzero element in $R(G).$
Since $H$ is an atom, then $H\cap s^*H=H,$\ $H\subset s^*H.$

Assume that $s^*H$ is strictly greater than $H.$
Then nontrivial decomposition
$$s^*H=H\cup(\neg H\cap s^*H))$$
takes place.
Apply a Boolean endomorphism $s_*.$
Then
$$H\subset s_*(s^*H)=s_*H\cup s_*(\neg H\cap s^*H)).$$
This gives the decomposition
$$H=(H\cap s_*(H))\cup(H\cap s_*(\neg H\cap s^*H)).$$
This decomposition is nontrivial (see [Pl1], p.347). Since the
summands lie in $R(G)$ and $H$ is an atom, the decomposition is
impossible. Hence, $H=s^*H.$ Let now $\overline\si_1$ and
$\overline\si_2$ be two elements in $H.$ Then $\overline
\si_1\cdot \overline\si_2=\overline{\si_1\si_2}=\overline{s_1}
s(\si_2)\in s^*H=H$, and the set $H$ is closed under
multiplication. Besides, $\overline e=\overline\si_1
s(\sigma)=\overline{\si_1\si},$ since $\overline e \in s^*H.$ This
means that for every $\overline\si\in H $ there is the inverse
element $\overline\si_1\in H.$

Hence, the whole set $H$ is the subgroup in the group $\Sigma_G$.

\proclaim{Proposition 17} Let $X$ be an arbitrary finite set in
$X^0$ and $A$ be an atom of the Boolean algebra $R(X).$ Then there
is a mapping $\mu: X\to G$ such that $s^*H=A$ holds for $s=s(\mu):
W(X)\to W(G).$\endproclaim

{Proof.} Take an arbitrary $\mu$ with $\overline\mu\in A.$ Take
$s=s(\mu)$ for such $\mu.$ Check that $s^*H=A.$

Take $\overline e\in H$ and $\overline e s(\mu)=\overline{e\mu}=\overline \mu\in A.$
The same $\mu$ lies also in $s^*H.$
We have $\overline\mu\in A\cap s^*H.$
Both components lie in $R(X),$ and $A\subset s^*H,$ since $A$ is an atom.
As previously, we see that $A=s^*H.$

\proclaim{Proposition 18} The subgroup $H$ in $\Sigma_G$ coincides
with the subgroup $R'$ in $\Aut(G).$
\endproclaim

{Proof.} We intend to prove that if $\si: G\to G$ is a
substitution, then $\overline\si\in H$ if and only if $\si$ is an
automorphism of the algebra $G,$ belonging to the subgroup $R'.$

Let, at the beginning, $\si\in R'.$ Take the point $\overline\si:
W(G)\to G.$ It is easy to understand that
$\si\overline\mu=\overline{\si\mu}$ for every mapping $\mu:X\to
G.$ If $X=G,$ then
$\si\overline{\mu}=\overline{\si\mu}=\overline\si\overline\mu.$ In
particular, $\si H=H$ for the given $\si.$ This means that
$\si\nu\in H\Leftrightarrow \nu\in H.$ Take an element $\overline
e$ in $H.$ Then $\si\overline e=\overline{\si e}=\overline\si\in
H.$ Thus, $\si\in R'$ implies $\overline\si\in H.$

Let us prove the opposite.
For every $A\subset \Hom(W(X),G)$ and $\si\in\Sigma_G$ define $\tilde\si A.$
The set $\tilde\si A$ is the set of all $\mu: W(X)\to G$ for which $\tilde
\si(\mu)=\nu\in A.$
Here $\mu=\tilde \si^{-1}(\nu)\in\widetilde{\si^{-1}}{}^*A.$
If $\si\in\Aut(G),$ then $\tilde\si A=\si A=(\si^{-1})^*A.$

Take now $A=H$ and let $\overline \si_1\in H.$
Check that $\widetilde{\si^{-1}}(\overline\si_1)=\overline{\si^{-1}}\ \overline
{\si_1}=\overline{\si^{-1}\si_1}.$
Take an arbitrary $g\in G.$
Then
$$\widetilde{\si^{-1}}(\overline
\si_1)(g)=\si^{-1}\si_1(g)=\overline{\si^{-1}\cdot\si_1}(g)=(\overline{\si^{-1}}\cdot
\overline\si_1)(g).$$
Hence,
$$\widetilde{\si^{-1}}{}^*H=\overline{\si^{-1}}H.$$
In the right hand part there is a coset by $H$ with the representative
$\overline{\si^{-1}}.$
In particular, if $\overline\si\in H,$ then
$\widetilde{\si^{-1}}{}^*H=H=\tilde\si H.$

Show now that $\tilde \si A=A$ for every atom $A$ in $R(X).$ Take
$\mu: X\to G$ with $\overline\mu\in A$ and let $s=s(\mu).$ Then
$A=s^*H,$\ $s(x)\in G$ for every $x\in X.$ We have $$\tilde
\si(A)=\widetilde{\si^{-1}}{}^*(A)=\widetilde{\si^{-1}}{}^*(s^*H)=s^*(\widetilde
{\si^{-1}}{}^*H) =s^*(H)=A$$ for $\si\in\Sigma_G.$ We used here
the fact that $\widetilde{\si^{-1}}(\nu
s)=\widetilde{\si^{-1}}(\nu)\cdot s$ holds for every $\nu\in H.$
Thus, all atoms are invariant under every $\tilde \si$ with
$\overline\si\in H.$ But then all $A\in R$ are also invariant over
such $\tilde \si.$

Show that this implies $\si\in \Aut(G).$

Let $\om\in\Om$ be an $n$-ary operation, $a_1,\dots,a_n$ elements in $G,$\
$\overline \si\in H.$
We need to check that $(a_1\dots a_n\om)^\si=a_1^\si\dots a_n^\si \om.$
Take variables $\{x_1,\dots,x_n,y\}\in X$ and proceed from the equality
$x_1\dots
x_n\om\equiv y.$
This is the equality of the $X$ type.
Take further $\mu: W(X)\to G$ with the condition
$\mu(x_1)=a_1,\dots,\mu(x_n)=a_n$ and
$\mu(y)=b=a_1\dots a_n\om.$
Then $\mu\in \Val(x_1\dots x_n\om=y)=A.$
By the condition, $A\in R(X),$\ $\tilde\si(A)=A.$
Thus, also $\tilde \si(\mu)\in A.$
We have:
$$\align \tilde\si(\mu)(x_1\dots x_n\om)&=(\tilde
\si(\mu)(x_1))\dots(\tilde\si(\mu)(x_n))\om\\
&=(\si\mu(x_1))\dots(\si\mu(x_n))\om=(\si
a_1)\dots(\si a_n)\om\\ &=\tilde\si(\mu)(y)=\si\mu(y)=\si b=\si(a_1\dots
a_n\om).\endalign$$
Hence, if $\tilde\si\in H,$ then $\si\in\Aut(G)$ and $\tilde
\si(A)=\si(A)=A$ for every $A\in
R,$\ $\si\in R'.$

Let us finish the proof of Theorem 10.

Let $R$ be a subalgebra in $\Hal_\Th(G),$ containing equalities.
We want to show that $R''=R.$
The inclusion $R\subset R''$ is always true.
Take $H=R'$ and let $\overline H$ be a corresponding atom in $R(G).$
This $\overline H$ generates the whole algebra $R.$
We have also $(R'')'=H.$
Then the algebra $R''$ is generated by the same $\overline H.$
Hence, $R''=R.$

Let now $H$ be a subgroup in $\Aut(G).$
A subset $\overline H$ in $\Hom(W(G),G),$ which is an element in $\Bool(W(G),G),
$ corresponds to
$H.$
Let us generate a subalgebra $R$ in $\Hal_\Th(G)$ by this element.
Check that $H=R'=\overline H'.$
This implies that $H'=R''=R$ and $H''=R'=H.$

Let us pass to the proof of Theorem 12.

Let $G_1$ and $G_2$ be two finite algebras in the given variety $\Th$.
Assume that the initial universum $X^0$  contains the sets $G_1$ and $G_2$
and consider
Halmos algebras $\Hal_\Th(G_1)$ and $\Hal_\Th(G_2)$.
We consider them together with the groups of automorphisms $\Aut(G_1)$ and
$\Aut(G_2).$
Take subalgebras $R_1$ and $R_2$ in these algebras, respectively, both with
equalities.
Suppose that there is an isomorphism
$$\g: R_1\to R_2.$$
We need to build an isomorphism $\dl: G_2\to G_1,$ inducing the
given $\g.$ Take a group $H_2=R_2'\subset\Aut(G_2)$ by the given
$R_2.$ As before, take a set $\overline
H_2\subset\Hom(W(G_2),G_2)$ in $R_2(G_2).$ This set is an atom in
$R_2(G_2)$ over the unit $\overline e.$ Apply $\g^{-1}.$ Then
$A=\overline H_2^{\g^{-1}}$ is an atom in $R_1(G_2)=\Bool
(W(G_2),G_1).$ Take an arbitrary $\mu: W(G_2)\to G_2$ in this
nonempty set $A.$ Let $\dl: G_2\to G_1$ be the restriction of
$\mu$ on $G_2.$

We want to show that every such $\dl$ solves the problem.

Consider a group $H_1=R_1'$ in $\Aut(G_1),$ and let $\overline H_1$ be a
corresponding atom in
$R_1(G_1).$
We have also an atom $A$ in $R_1(G_2).$
As we know, $A=s^*\overline H_1=H_2^{\g^{-1}}.$

Take $B=\overline H_1^\g=s_1^*(\overline H_2)$ by $\overline H_1.$
Here $s_1=s(\si)$ with $\si: G_1\to G_2,$ selected in $B$ as $\dl$ in $A.$

 Using ([Pl1], p.351),
we prove that $\dl: G_2\to G_1$ is a bijection, as well as $\si:
G_1\to G_2.$ We also prove that $\dl$ and $\si$ are isomorphisms
of algebras ([Pl1], p.356) and that $H_1=\dl H_2\dl^{-1},$ or,
similarly $R_1'=\dl R_2'\dl^{-1}.$ 
This allows to show that $\dl_*: \Hal_\Th(G_1)\to \Hal_\Th(G_2)$
induces an isomorphism $R_1\to R_2.$ We check further that this
isomorphism coincides with the initial $\g.$ It is sufficient to
verify the last one on the set $A,$ generating the algebra $R_1.$
The equality $\g(A)=\dl_*(A)$ follows from the definition.

\head{7. Geometric properties of models}
\endhead

\subhead{7.1. Isomorphisms of categories}
\endsubhead

We consider a new model $(G_1,\Phi,f^\dl)$ for the given isomorphism $\dl:
G_2\to G_1$ and a
model $(G_1,\Phi,f).$
This new model is isomorphic to the initial one.
We consider also the commutative diagram

$$
\CD
\Hal_\Theta(\Phi) @>\Val_{f}>> \Hal_\Theta(G_1) \\
 @. @/SE/\Val_{f^{\delta}}// @VV\delta_* V\\
 @. \Hal_\Theta(G_2)\
\endCD
$$


If $A=T^f$ is an algebraic set in $\Hal_\Th(G_1),$ then
$\dl_*T^{f^\dl}$ is an algebraic set in $\Hal_\Th(G_2).$ Hence, we
have a bijection between $f$-algebraic sets in $\Hal_\Th(G_1)$ and
$f^\dl$-algebraic sets in $\Hal_\Th(G_2).$ We can  speak now of a
commutative diagram
$$
\CD
\Hal_\Theta(\Phi) @>\Val_{f}>> K_{\Phi\Theta}(f) \\
 @. @/SE/\Val_{f^{\delta}}// @VV\delta_* V\\
 @.  K_{\Phi\Theta}(f^\delta) \
\endCD
$$


Here the bijection $\dl_*$ is well coordinated with the morphisms of the
categories, and
thus, is an isomorphism of categories.
Besides, $\dl_*$ induces an isomorphism of algebras $R_f$ and $R_{f^\dl}.$

We are  interested  in the general problem of isomorphism of two
categories of the type $K_{\Phi\Th}(f).$ Here we define a {\it
special} isomorphism. For the models $(G_1,\Phi,f_1)$ and
$(G_2,\Phi,f_2)$, it is an isomorphism, determined by the diagram
$$
\CD
\Hal_\Theta(\Phi) @>\Val_{f_1}>> K_{\Phi\Theta}(f_1) \\
 @. @/SE/\Val_{f_2}// @VV\gamma V\\
 @.  K_{\Phi\Theta}(f_2) \
\endCD
$$



Here $\g$ is an isomorphism of categories.
It follows from the definition that the same $\g$ induces an isomorphism of
algebras
$R_{f_1}$ and $R_{f_2}.$
Besides, $\g(T^{f_1})=T^{f_2}$ for every  $T\subset \Hal_\Th(\Phi)$ of the
definite type $X.$

\proclaim{Theorem 13}
If the algebras $G_1$ and $G_2$ are finite, then the categories
$K_{\Phi\Th}(f_1)$ and
$K_{\Phi\Th}(f_2)$ are specialy isomorphic if and only if the models are
isomorphic.\endproclaim

{Proof.} We have already seen this in one direction. Now let us
have a special isomorphism $\g: K_{\Phi\Th}(f_1)\to
K_{\Phi\Th}(f_2).$ We show that there exists an isomorphism of
models $\dl: G_2\to G_1,f_1=f_2^\dl,$ such that $\g=\dl_*$ on the
objects of the categories. We have an isomorphism $\g: R_{f_1}\to
R_{f_2}.$ According to Theorem 12, there is an isomorphism of
algebras $\dl: G_2\to G_1,$ such that $\g$ and $\dl_*$ coincide on
the elements from $R_{f_1}.$ From this  follows that $\g$ and
$\dl_*$ coincide on the objects of the categories. It is also easy
to see that the commutative diagram

$$
\CD
\Hal_\Theta(\Phi) @>\Val_{f_1}>>\Hal_\Theta(G_1)  \\
 @. @/SE/\Val_{f_2}// @VV\delta_* V\\
 @. \Hal_\Theta(G_2)   \
\endCD
$$

\noindent
takes place.

Applying the remarks from the end of the Section 7, we conclude that
$f_1=f_2^\dl$ and the
models are isomorphic.

Along with the special isomorphisms we consider {\it strict}
isomorphisms. It is an isomorphism $\g: K_{\Phi_1\Th}(f_1)\to
K_{\Phi_2\Th}(f_2),$ inducing an isomorphism of algebras
$R_{f_1}\to R_{f_2}.$ Here $\Phi_1$ and $\Phi_2$ are different in
general, and the models have the form $(G_1,\Phi_1,f_1)$ and
$(G_2,\Phi_2,f_2).$ We do not assume relations with the algebras
of formulas. We cannot speak here about isomorphism of models and
we use the notion of automorphic isomorphism.
\proclaim{Definition}
The modes $(G_1,\Phi_1,f_1)$ and $(G_2,\Phi_2,f_2)$ are called {\it
automorphically
isomorphic}, if

1. Algebras $G_1$ and $G_2$ are isomorphic.

2. There exists an isomorphism of algebras $\dl: G_2\to G_1$ with
$\Aut(f_1)=\dl\Aut(f_2)\dl^{-1}$.

Here the groups of automorphisms of models are conjugated.
\endproclaim

If the models are, in particular, isomorphic, $\Phi_1=\Phi_2,$
then they are automorphically isomorphic. \proclaim{Theorem 14}
For the finite algebras $G_1$ and $G_2$ the categories
$K_{\Phi_1\Th}(f_1)$ and $K_{\Phi_1\Th}(f_2)$ are strictly
isomorphic if and ony if the models are automorphically
isomorphic.\endproclaim

{Proof.} Let a strict isomorphism
$$\g:K_{\Phi_1\Th}(f_1)\to K_{\Phi_2\Th}(f_2)$$
be given.
By the condition, this $\g$ induces an isomorphism of algebras $\g:
R_{f_1}\to R_{f_2}.$
Once more by Theorem 12, we have an isomorphism $\dl: G_2\to G_1$ such that
$\dl_*$ induces
$\g.$
Here, $R_{f_2}=\dl_*(R_{f_1}).$
This, in its turn, gives $R_{f_1}'=\dl R_{f_2}'\dl^{-1}.$
But we know from Theorem 9 that $R_{f_1}'=\Aut(f_2),$\  $R_{f_2}'=\Aut(f_2).$
Thus, $\Aut(f_1)$ and $\Aut(f_2)$ are conjugated by the isomorphism $\dl:
G_2\to G_1.$
We have verified that the models are automorphically isomorphic.

Let us prove the opposite.
Let the models $(G_1,\Phi_1,f_1)$ and $(G_2,\Phi_2,f_2)$ be automorphically
isomorphic with
the isomorphism of algebras $\dl: G_2\to G_1.$
We have an isomorphism $\dl_*:\Hal_\Th(G_1)\to \Hal_\Th(G_2).$
We have a subalgebra $R_{f_1}$ in $\Hal_\Th(G_1)$ and $R_{f_2}$ in
$\Hal_\Th(G_2).$
Here $R_{f_1}'=\Aut(f_1)$ and  $R_{f_2}'=\Aut(f_2)$.
By the condition we have $R_{f_1}'=\dl R_{f_2}'\dl^{-1}.$
This means that $\dl_*(R_{f_1})=R_{f_2}.$
Indeed, $(\dl_*(R_{f_1}))'=\dl^{-1}R_{f_1}'\dl=R_{f_2}'.$
Hence, $\dl_*(R_{f_1})=R_{f_2}$ and $\dl_*(R)'=\dl^{-1}R'\dl.$
This always takes place.
Then $\dl_*$ induces an isomorphism $R_{f_1}\to R_{f_2}.$

Let now $A\in R_{f_1}.$
Then $A=\Val_{f_1}(u)$ for some $u\in\Hal_\Th(\Phi_1),$\ $\dl_*(A)=B\in
R_{f_2},$ \
$B=\Val_{f_2}(v),$\ $v\in \Hal_\Th(\Phi_2).$
Here $u$ and $v$ are of the same type $X.$
We have $\dl_*\Val_{f_1}(u)=\Val_{f_2}(v).$

Let, further, $A$ be an object of the category $K_{\Phi_1\Th}(f_1),$\
$A=T^{f_1},$ where
$T$ is a collection of formulas in $\Hal_\Th(\Phi)$ and $T$ is of type $X.$
We have
$$A=\bigcap_{u\in T}\Val_{f_1}(u).$$
We have also
$$\dl_*(A)=\bigcap_{u\in T}\dl_*\Val_{f_1}(u)=\bigcap_{v\in T^*}\Val_{f_2}(v).$$
Here $T^*$ is a collection of formulas $v$ in $\Hal_{\Th}(\Phi_2)$ of the
type $X,$ somehow
connected with the collection $T.$
In any case, we may claim that $\dl_*A,$ as well as $A,$
 is an algebraic set.
This means that there exists a bijection on the objects
$K_{\Phi_1\Th}(f_1)\to
K_{\Phi_2\Th}(f_2).$
Let now $s: A_1\to A_2$ be a morphism in $K_{\Phi,\Th}(f_1).$
We have $A_1\subset sA_2.$
But then $\dl_*(A_1)\subset \dl_*(sA_2)=s\dl_*A_2.$
Hence, the same $s$ gives a morphism $s:\dl_*(A_1)\to \dl_*(A_2).$
The opposite is also true.
Thus, we come to the isomorphism of categories $\dl_*: K_{\Phi,\Th}(f_1)\to
K_{\Phi_2\Th}(f_2),$ and this isomorphism is a strict one.
This completes the proof of Theorem 14.

\subhead{7.2. A remark on the categories $K_{\Phi\Th}(G)$}
\endsubhead

Here $G$ is an arbitrary algebra in $\Th.$ Recall that the objects
of this category have the form $(X,A,f),$ where $f$ is a
interpretation of the given set $\Phi$ in the given algebra $G.$

Consider $K_{\Phi\Th}(G_1)$ and $K_{\Phi\Th}(G_2)$, while the
isomorphism $\dl: G_1\to
G_2$ is given.
To the model $(G_1,\Phi,f)$ corresponds a model $(G_2,\Phi,f^\dl).$
Simultaneously, to the object $(X,A,f)$ corresponds an object
$(X,\dl_*(A),f^\dl).$
This gives an isomorphism of categories
$$K_{\Phi\Th}(G_1)\to K_{\Phi \Th}(G_2).$$
If, further, $\dl: G\to G$ is an automorphism of algebras, then an
automorphism of categories
$K_{\Phi\Th}(G)$ corresponds to $\dl.$
This leads to the representation
$$\Aut(G)\to\Aut(K_{\Phi \Th}(G)).$$
This representation is naturally tied with the representation
$$\Aut(G)\to \Aut(\Hal_\Th(G)).$$

\subhead{7.3. Geometrical properties }
\endsubhead

\proclaim{Definition 1} The model $(G,\Phi,f)$ is called {\it
geometrically noetherian} if for every finite set $X$ and every
set of formulas $T$ in $\Hal_{\Phi\Th}(X)$ 
 there is some
finite part $T_0\subset T$ with $T_0^f=T^f.$ This means that
$T_0^{ff}=T^{ff}.$
\endproclaim

\proclaim{Theorem 15} The model $(G,\Phi,f)$ is geometrically
noetherian if and only if the minimality condition holds in the
lattice $\Alv_f(X).$ Correspondingly, in the lattice $\Cl_f(X)$ we
have the maximality condition.
\endproclaim

Now let $T_0=\{u_1,\dots, u_n\}$ and $u$ is $u_1\wedge\dots\wedge u_n$.
Then $T_0^f=\Val_f(u)=\{u\}^f.$
Thus, we may claim that if the model $(G,\Phi,f)$ is geometrically
noetherian, then every
algebraic set over this  model is  a simple algebraic set.

However, the corresponding element $u=u_1\wedge\dots\wedge u_n$
does not necessarily belong to the initial set $T.$ We call a
model $(G,\Phi,f)$ {\it weak geometrically noetherian} if every
algebraic set over it is a simple algebraic set. Weak noetherian
model is not necessarily geometrically noetherian.

However, we may claim that every finite model is geometrically
noetherian (compare Theorem 8 in Galois theory). We may also claim
that any finite cartesian product of geometrically noetherian
models is also a geometrically noetherian model. A submodel of a
geometrically noetherian model is also geometrically noetherian.
 The similar properties are not true in general, in respect to  to cartesian
powers.



Let us pass to the notion of geometrical equivalence of two models.
Let two models $(G_1,\Phi,f_1)$ and $(G_2,\Phi,f_2)$ with the same $\Phi$
be given.

\proclaim{Definition 2}
The models $(G_1,\Phi,f_1)$ and $(G_2,\Phi,f_2)$ are geometrically
equivalent if
$T^{f_1f_1}=T^{f_2f_2}$ holds for every finite $X$ and every $T$ in
$\Hal_{\Phi\Th}(X).$

If the models $(G_1,\Phi,f_1)$ and $(G_2,\Phi,f_2)$ are geometrically
equivalent, then

1. The lattices $\Cl_{f_1}(X)$ and $\Cl_{f_2}(X)$ coincide,
while the lattices
$\Alv_{f_1}(X)$ and $\Alv_{f_2}(X)$ are isomorphic.

2. The categories $C_{\Phi\Th}(f_1)$ and $C_{\Phi\Th}(f_2)$
coincide, while the
categories $K_{\Phi\Th}(f_1)$ and $K_{\Phi\Th}(f_2)$ are isomorphic.

3. These models are elementary equivalent.

It follows from the first claim that if the models are geometrically
equivalent and one of
them is geometrically noetherian, then the second one is also geometrically
noetherian.
\endproclaim

\proclaim{Theorem 16} Let the model $(G,\Phi,f)$ be geometrically
noetherian. Then each of its utrapower is also geometrically
noetherian, and all these ultrapowers are geometrically equivalent
to the initial $(G,\Phi,f).$\endproclaim

All the described notions are naturally tied with the logic of
generalized (infinitary) formulas of the kind $(\wedge_{u\in
T})\to v,$\ or, what is the same $T\to v.$ For the geometrically
noetherian models it is sufficient to proceed from the usual
finite formulas.


Elementary equivalence of the models does not generally imply
their geometrical equivalence. However, we may claim the following

\proclaim{Theorem 17} If two models are elementary equivalent and
one of them is geometrically noetherian, then the second one is
also geometrically noetherian and these models are geometrically
equivalent.
\endproclaim

In concern with the notion of geometrically noetherian model let us
return to the notion of the logical kernel of a homomorphism of
the form $\mu: W(X) \to G$. It is easy to see that if $\Log_f\Ker(\mu)$
is such a kernel and $\Val_f(\Log_f\Ker(\mu))$ is its image in the algebra
$R_f(X)$ then this image is a principal ultrafilter if the model
$(G,\Phi, f)$ is geometrically noetherian. This ultrafilter is
generated by the algebraic set ${\{\mu\}}^{ff}$.

\head{8. Applications to the knowledge science}
\endhead

\subhead{8.1 Introduction}
\endsubhead

Knowledge theory and knowledge bases provide an important example
of the field where applications of universal algebra and algebraic
logic are very natural, and their interacting with quite practical
problems arising in computer science is very productive. Another
examples of such interaction are given by relational database
theory, constraint satisfaction problem (\cite{BJ},\cite{JCP}),
theory of complexity, and by others.

One can speak about knowledge and a system of knowledge. As a
rule, a domain of knowledge or of a system of knowledge is fixed.
 In our approach only knowledge
that allows a formalization in some logic is considered. The logic
may be different. It is often oriented towards the corresponding
field of knowledge cf. \cite{G},\cite{L},\cite{S}.

In this paper we focus on the special situation of elementary
knowledge.

Elementary knowledge is considered to be a
 first order knowledge, i.e., the knowledge
that can be represented by the means of the First Order Logic
(FOL). The corresponding applied field (field of knowledge) is
grounded on some variety of algebras $\Theta$, which is arbitrary
but fixed. This variety $\Theta$ is considered as a knowledge
type. Its counterpart in database theory is the notion of datatype
$\Theta$.

 We also fix a set of symbols of relations $\Phi$.
The subject of knowledge is a triple $(G, \Phi, f)$, where $G$ is
an algebra in $\Theta$ and $f$ is a interpretation of the set
$\Phi$ in $G$. It is a model in the ordinary mathematical sense.
As a rule, we use  shorthand and write $f$ instead of $(G, \Phi,
f)$.  For the given $\Phi $ we denote the corresponding applied
field by $\Phi \Theta$.

FOL is also oriented towards the variety $\Theta$.

We assume that every knowledge under consideration is represented
by three components: \pmf 1) {\it The description of  knowledge.}
It is a syntactical part of knowledge, written out in the language
of the given logic. The description reflects, what do we want to
know. \pmf 2) {\it The subject of  knowledge} which is an object
in the given applied field, i.e., an object for which we determine
knowledge. \pmf 3) {\it The content of  knowledge} (its
semantics).

The first two components are relatively independent, while the
third one is uniquely determined by the previous two. In the
theory under consideration, this third component has a geometrical
nature.  In some sense it is an algebraic set in an affine space.
If $T$ is a description  of knowledge and $(G, \Phi, f)$ is a
subject, then $T^f$ denotes the content of knowledge. We would
like to equip the content with its own structure, algebraic or
geometric, and to consider some aspects of such structure.

 We want to underline that there are three aspects in our approach
 to knowledge representation:  logical (for knowledge
 description), algebraic (for the subject of knowledge) and
 geometric (in the content of knowledge).
This geometry is of algebraic nature. However, the involved
algebra inherits some geometric intuition.


 We consider categories of elementary knowledge.  The language of
 categories in knowledge theory is a good way to organize and
 systematize primary elementary knowledge.  Morphisms in a knowledge
 category  give links between knowledge.  In particular, one
 can speak of isomorphic knowledge.  The categorical approach also
 allows us to use ideas of monada and comonada [ML].  It turns out that
 this leads to some general views on enrichment and
 computation of knowledge.  Enrichment of a structure can be
 associated with a suitable monada over a category, while the
 corresponding computation is organized by comonada.

 Let us make one more remark.  In every well described
 field of knowledge, one can study a category of elementary
 knowledge belonging to this field.  Consideration of such
 categories might be of special interest.
\subhead{8.2 The category of knowledge}
\endsubhead

Define first the category  $\Know_{\Phi\Theta}(f)$.

Fix a model (subject of knowledge) $(G, \Phi, f)$. Let us define a
category of knowledge for this model and denote it by Know$_{\Phi
\Theta}(f)$. This is the knowledge category for the given subject
of knowledge. Since the model is fixed, the objects of the
category \knowf \ have to have the form $(X, T, A)$.  We do not
fix the subject of knowledge in the notation of the object,
since it is fixed in the notation of the category.

  The set $X$ is
multi-sorted.  It marks the ``place" where the knowledge is
situated. The set $X$ points also  the ``place of the knowledge'',
i.e., the space of the knowledge $Hom (W(X), G)$, while the
subject of the knowledge $(G, \Phi, f)$ is given. The set $T$ is
the description of the knowledge in the algebra $Hal_{\empty
\Theta}(X)$, and $A=T^f$ is the content of knowledge, depending on
$T$ and $f$. The set $T^{ff} = A^f$ is the full description of the
knowledge $(X, T, A)$ which is a Boolean filter in Hal$_{\Phi
\Theta}(X)$.

Now about morphisms $(X, T_1, A) \to (Y, T_2, B)$. Take $s\colon
W(Y) \to W(X)$.  We have also $s\colon \Hal_{\Phi \Theta}(Y)\to
\Hal_{\Phi \Theta}(X)$ (see 2.2). This is a homomorphism of
Boolean algebras.  The homomorphism $s$ gives rise  to
$$
\tilde s \colon Hom (W(X), G) \to Hom (W(Y), G).$$ As above, the
first $s$ is admissible for $A$ and $B$ if $\tilde s (\nu) =\nu s
\in B$ for every point $\nu\colon W(X)\to G$ in $A$.

As we know, $s$ is admissible for $A$ and $B$ if and only if $su
\in A^f$ for every $u \in B^f$.  This holds for  $s_\ast$, for
which we have also a homomorphism $\ol s:\Hal\em (Y)/B^f\to
\Hal\em (X)/A^f$.  It is easy to prove that $s$ is admissible for
$A$ and $B$ if and only if $s u \in A^f$ holds for every $u \in
T_2$. We consider admissible $s$ as a morphism
$$
s\colon (X, T_1, A) \to (Y, T_2, B),
$$
in the weak category \knowf.

  We have $\tilde s (\nu) = \nu s \in
B$ if $\nu \in A$, and $s$ induces a mapping $[s]\colon A \to B$.
Simultaneously, there is a mapping $s\colon T_2 \to A^f$ and a
homomorphism
$$
\ol s\colon \Hal\em (Y)/B^f \to \Hal\em(X)/A^f.
$$
We have already mentioned (Proposition 2) that $\ol s_1 = \ol s_2$
follows from $[s_1]=[s_2]$. Thus, we can take the morphisms of the
form
$$
[s]\colon (X, T_1, A) \to (Y, T_2, B),
$$
for the morphisms of the exact category \knowf. The canonical
functors \knowf$\to K\em(f)$ for weak and exact categories are
given by the transition $(X, T, A) \to (X, A)$. In this transition
we ``forget" to fix the description of knowledge $T$.

Now we define the  category $\Know_{\Phi\Theta}$.
Let us define the category of elementary knowledge for the whole
applied field $\Phi\Theta$; the subject of the knowledge $(G,
\Phi, f)$ is not fixed.  As earlier, we proceed from the category
$\Phi \Theta$ whose morphisms are homomorphisms in $\Theta$. They
ignore the relations from $\Phi$.

An object of the knowledge category \know \   has the form
$$
(X, T, A; (G, \Phi, f)),
$$
and we write $(X, T, A; G,  f)$, because $\Phi$ is fixed for the
category. Here $X$ marks the place of knowledge. The components $A
= T^f$, $G$ and $f$ may change.

Consider morphisms:
$$
(X, T_1, A; G_1, f_1) \to (Y, T_2, B; G_2, f_2).
$$
We apply the same approach as in Section 3.3 with some
modifications.

Start from $s: W(Y) \to W(X)$ and $\delta: G_1 \to G_2$. These $s$
and $\delta$ should correlate. Let us explain the correlation
condition.  Take a set $A_1 = \{ \delta \nu, \nu \in A\}=\delta^*
A$ and
 take further $T^\delta_1 = A_1^{f_2}$.
 Correlation of $s$ and $\delta$ means that $su\in T^\delta_1$ holds for
 any $u \in T_2$.  The same holds for every $u \in B^{f_2}$.
 The last also says that there is a homomorphism
 $$
 \ol s\colon \Hal\em (Y) /B^{f_2} \to \Hal\em (X)/A_1^{f_2}.
 $$
 The first of the two mappings $(s, \delta) \colon A \to B$ and
 $s\colon T_2 \to T^\delta_1$ transforms the content of knowledge,
 while the second one acts on the description.  Here $T_2$ and
 $T^\delta_1$ describe knowledge associated with the same subject
 $(G_2, \Phi, f_2)$.

 With the fixed $\delta$ there is also an exact mapping $([s], \delta): A \to
 B$.  This brings us to  weak and exact categories \know.  The morphisms of
 the first one are $(s, \delta)$ and in the second one they are of the form $([s],
 \delta)$ for $(X, T_1, A; G_2, f_1) \to (Y, T_2, B;
 G_2, f_2)$.
 The canonical functors \know$ \to K\em$ are defined by the
 transition
 $$
 (X, T, A; G, f) \to (X, A; G, f).
 $$
 As above, we remove the description of knowledge from the
 notations.

For given algebras $G$ in $\Theta$ we consider the  categories
$K\em(G)$ and $\know_{\Phi\Theta}(G)$.

An algebra $G \in \Theta$ is fixed in the categories $K\em (G)$
and $\know_{\Phi\Theta}(G).$  A set of symbols of relations $\Phi
$ is fixed as usual, but interpretations $f$ of $\Phi$ in $G$ may
change. Thus, $K\em(G)$ is a subcategory in $K\em$ and
$\know_{\Phi\Theta}(G)$ is a subcategory in $\know_{\Phi\Theta}.$
Here the corresponding $\delta: G \to G$ are identical
homomorphisms. Objects of the category $K\em(G)$ have the form
$(X, A, f)$, and those of the category $\know_{\Phi\Theta}(G)$ are
written  as $(X, T, A, f)$. There is a canonical functor
$\know_{\Phi\Theta}(G) \to K\em(G)$. As for morphisms
$$
\eqalign{ &(X, A, f_1) \to (Y, B, f_2) \; \; \hbox{\rm and} \cr
&(X, T_1, A, f_1) \to (Y, T_2, B, f_2),\cr}
$$
we  note that $A=A_1, A_1^{f_2} = T^\delta_1$ and
$A^{f_2}=T_1^{f_1f_2}$. Hence, the corresponding admissible $s:
W(Y) \to W(X)$ transfers each $u \in T_2$ into $s u \in T_1^{f_1
f_2}$ and it induces a homomorphism
$$
\ol s: \Hal\em (Y) /B^{f_2} \to \Hal\em (X) /A^{f_2}.
$$
Every $s$ yields a mapping $[s]: A \to B$. This provides a
morphism $(X,A,f_1)\to (Y,B,f_2)$.

\subheading{8.3 Category of knowledge description (the category
$L_\Theta(\Phi)$)}.

Denote the category of knowledge description by $L_{\Phi\Theta}$
or $L_\Theta(\Phi)$.

Its objects are of the form $(X,T)$, where $X$ is a finite set and
$T$ is a set of formulas of $\Hal_{\Phi\Theta}(X)$. Define
morphisms $(X,T_1)\to(X,T_2)$. According to description of the
category
 $\Hal_{\Theta}(\Phi)$ proceed from the functor $\Theta^0\to
\Hal_{\Theta}(\Phi)$ which assigns a mapping $s_\ast:
\Hal_{\Phi\Theta}(X) \to \Hal_{\Phi\Theta}(Y)$ to every
homomorphism $s: W(X)\to W(Y)$. We say that $s$ is{\it admissible
in respect to $T_1$ and $T_2$} if $s_\ast(u)\in T_2$ for every
$u\in T_1$. For such admissible $s$ we have a mapping $s_\ast:
T_1\to T_2$ which determines
$$
s_\ast: (X, T_1)\to (X,T_2)
$$

\subheading{8.4 Functor of transition from knowledge description
to knowledge content (the functor $\Ct_f$)}

Proceed from the model $(G,\Phi, f)$ and consider a functor
$$
\Ct_f:L_{\Phi\Theta}\to K_{\Phi\Theta}(f).
$$
Here, $K_{\Phi\Theta}(f)$ is the corresponding category  of
algebraic (elementary) sets over the given model and $\Ct$ stands
for "contents". The functor $\Ct_f$ is a contravariant one. To
every object $(X,T)$ of the category $L_{\Phi\Theta}$ it assigns
the corresponding content $(X, T^f)=(X,A)$ which is an object of
the category $K_{\Phi\Theta}(f)$.

Now one has to define the functor $\Ct_f$ on morphisms. Let  a
morphism
$$
s_\ast: (Y, T_2)\to (X,T_1)
$$
be given for $s: W(Y)\to W(X)$. Show that $s$ induces a morphism
$$
\widetilde {s_\ast}: (X, A)\to (Y,B),
$$
where $A=T_1^f$, and $B=T_2^f$.

We proceed from $\widetilde s: \Hom(W(X),G)\to \Hom (W(Y),G)$.

Let us define a transition $s\to \widetilde s.$ 

Check first that if $s$ is admissible for $T_2$ and $T_1$ then
this $ s$ is admissible for $A=T_1^f$ and $B=T_2^f$. The last
means that $\widetilde s(\nu)\in B$ if $\nu\in A$. The inclusion
$\nu\in A$ says that $\nu\in \Val_f(v)$ for every $v\in T_1$. We
need to verify that $\nu s\in B,$ that is $\nu s \in \Val_f(u)$
for every $u\in T_2$.

Take an arbitrary $u\in T_2$. We have: $v= {s_\ast}(u)\in T_1$;
$\nu\in \Val_f(v)=\Val_f(s_\ast u)= s \Val_f(u)$. This gives $\nu
s  \in \Val_f(u)$. We used that $s$ and $\Val_f$ commute, since $
\Val_f$ is a homomorphism of algebras.

 The mapping $[s]: A\to B$
corresponds to the homomorphism $s:W(Y)\to W(X)$. This mapping is
considered simultaneously as a morphism in the category
$K_{\Phi\Theta}(f)$ (see 3.2)
$$
[s]: (X,A)\to (Y,B).
$$
We define: $\Ct_f(s_\ast)=\widetilde s_\ast=[s]$.

 Check now compatibility of the definition of
$\Ct_f$ with the multiplication of morphisms. Given $s_1: W(X)\to
W(Y)$ and $s_2: W(Y)\to W(Z)$ we have $s_2s_1: W(X)\to W(Z)$.
Using the fact that the transition $\Theta^0\to \Hal_\Theta(\Phi)$
is a functor, we get $(s_2s_1)_\ast= s_{2\ast}s_{1\ast}$. Here, we
have
$$
s_{1\ast}:\Hal_{\Phi\Theta}(X)\to \Hal_{\Phi\Theta}(Y),
$$
$$
s_{2\ast}:\Hal_{\Phi\Theta}(Y)\to \Hal_{\Phi\Theta}(Z),
$$
and
$$
(s_2s_1)_\ast:\Hal_{\Phi\Theta}(X)\to \Hal_{\Phi\Theta}(Z).
$$
Let $(X,T_1)$, $(Y,T_2)$ and $(Z,T_3)$ be objects in
$L_\Theta(\Phi)$ , and $s_1,s_2$ admissible in respect to $T_1$,
$T_2$ and, correspondingly, for $T_2$, $T_3$. In this case there
are morphisms
$$
s_{1\ast}:(X, T_1)\to (Y,T_2),
$$
$$
s_{2\ast}:(Y, T_2)\to (Z,T_3),
$$
and
 $$ s_{2\ast}s_{1\ast}= (s_2s_1)_\ast: (X,T_1)\to (Z,T_3).$$

Take $T_1^f=A$, $T_2^f=B$, $T_3^f=C$. We have
$$
\widetilde {s_{1\ast}}: (Y,B)\to (X,A),
$$
$$
\widetilde {s_{2\ast}}: (Z,C)\to (Y,B),
$$
and
$$
\widetilde{{s_2s_1}_\ast}=\widetilde{s_{1\ast}}\widetilde{s_{2\ast}}:
(Z,C)\to (X,A).
$$
This gives compatibility of the functor $\Ct_f$ with the
multiplication of morphisms. Compatibility with the unity morphism
is evident. This finishes the definition of the contravariant
functor $\Ct_f: L_{\Phi\Theta}\to K_{\Phi\Theta}(f)$.



\subheading{ 8.5 Homomorphisms of Halmos algebras
$\Hal_\Theta(\Phi)$ and functors of the categories
$L_\Theta(\Phi)$}

Given a homomorphism $\beta:\Hal_\Theta(\Phi_1)\to
\Hal_\Theta(\Phi_2)$,  define the corresponding functor
$\widetilde\beta:L_\Theta(\Phi_1)\to L_\Theta(\Phi_2)$. For every
set of formulas $T\subset \Hal_{\Phi_1\Theta}(X),$ denote by
$T^\beta$ the set $T^\beta=\{u^\beta, u\in T\}$. If $(X,T)$ is an
object in $L_\Theta(\Phi_1)$, then, setting
$$
\widetilde\beta(X,T)=(X,T^\beta),
$$
we get an object in $L_\Theta(\Phi_2)$.

In order to define the functor $\widetilde\beta$ on morphisms let
us make a remark. Proceed from the functors
$\Theta^0\to\Hal_\Theta(\Phi_1)$ and
$\Theta^0\to\Hal_\Theta(\Phi_2)$. The morphisms
$$
s_{\ast}^1:\Hal_{\Phi_1\Theta}(X)\to \Hal_{\Phi_1\Theta}(Y),
$$
$$
s_{\ast}^2:\Hal_{\Phi_2\Theta}(X)\to \Hal_{\Phi_2\Theta}(Y)
$$
correspond to every $s: W(X)\to W(Y)$. We have also
$$
\beta=(\beta_X, X\in \Gamma^0):\Hal_{\Theta}(\Phi_1)\to
\Hal_{\Theta}(\Phi_2).
$$
The fact that the homomorphism $\beta$ is compatible with the
operation $s$ is represented by the commutative diagram

$$
\CD
 \Hal_{\Phi_1\Theta}(X)@>s_\ast^1>> \Hal_{\Phi_1\Theta}(Y)\\
@V\beta_X VV @VV\beta_Y V\\
\Hal_{\Phi_2\Theta}(X)@>s_\ast^2>> \Hal_{\Phi_2\Theta}(Y)\\
\endCD
$$
So, for a homomorphism $s: W(X)\to W(Y)$ we have the equality
$\beta_Ys_\ast^1(u)=s_\ast^2\beta_X(u)$ for every $u\in
\Hal_{\Phi_1\Theta}(X).$

Now we are able to define an action of the functor
$\widetilde\beta$ on morphisms. Let a morphism $s_\ast^1:
(X,T_1)\to(Y,T_2)$ in the category $L_{\Phi_1\Theta}$ be given and
$s_\ast^1(u)\in T_2$ if $u\in T_1$. Then, we have $s_\ast^2(v)\in
T_2^\beta$ if $v\in T_1^\beta$.

Indeed, let $v=\beta_X(u)$, $u\in T_1$, $v\in T_1^{\beta_X}$. We
have:
$$
s_\ast^2\beta_X(u)=s_\ast^2(v)=\beta_Ys_\ast^1(u)\in
T_2^{\beta_Y},
$$
since $s_\ast^1(u)\in T_2$. Hence, $s_\ast^2(v)\in T_2^{\beta_Y}$
for every $v=\beta_X(u)\in T_1^{\beta_X}$.

 We set
$s_\ast^2=\widetilde\beta(s_\ast^1): T_1^{\beta_X}\to
T_2^{\beta_Y}$.
 A morphism
 $$
 s_\ast^2=\widetilde\beta(s_\ast^1):
(X,T_1^{\beta_X})\to (Y,T_2^{\beta_Y})
$$
 corresponds to $s_\ast^1:
(X,T_1)\to (Y,T_2)$.

Check now compatibility of the transition $s_\ast^1\to s_\ast^2$
with the multiplication of morphisms. Given $s_1: W(X)\to W(Y)$
and $s_2: W(Y)\to W(Z)$, we have $s_2s_1: W(X)\to W(Z)$. Using
once more the fact that the transition $\Theta^0\to
\Hal_{\Theta}(\Phi)$ is a functor, we get
$$
(s_2^1s_1^1)_\ast= s_{2\ast}^1s_{1\ast}^1,
$$
$$
(s_2^2s_1^2)_\ast= s_{2\ast}^2s_{1\ast}^2,
$$

Apply $\widetilde\beta$. We need to verify that
$\widetilde\beta(s_{2\ast}^1s_{1\ast}^1)=
\widetilde\beta(s_{2\ast}^1)\widetilde\beta(s_{1\ast}^1). $ We
have
$$
\widetilde\beta(s_{2\ast}^1s_{1\ast}^1)=\widetilde\beta(s_{2}^1s_{1}^1)_\ast=
(s_{2}^2s_{1}^2)_\ast=
s_{2\ast}^2s_{1\ast}^2=\widetilde\beta(s_{2\ast}^1)\widetilde\beta(s_{1\ast}^1).
$$
This gives compatibility with the multiplication as well as with
the unit. Hence, we have the functor $\widetilde\beta:
L_\Theta(\Phi_1)\to L_\Theta(\Phi_2)$.

\subhead{ 8.6 Knowledge bases}
\endsubhead

We proceed from a multi-model $(G, \Phi, F)$. A multi-model  $(G,
\Phi, F)$ defines a system of models  $(G, \Phi, f,)$ where $f$
runs the set $F$. Here $G$ is an algebra in $\Theta$, and $\Phi$
is a set of relations. Recall that both the algebra $G\in \Theta$
and a relation $f\in F$ are multi-sorted.  The set $F$ is a set of
instances $f$, where $f$ is a interpretation of the set $\Phi$ in
$G$.

To every such multi-model corresponds  a knowledge base $KB =
KB(G,\Phi,F)$. The definition slightly differs from that of
\cite{PTP}.

\proclaim{Definition} A knowledge base $KB = KB(G,\Phi,F)$
consists of two categories. The first one is the category of
knowledge description $L_\Theta(\Phi)$, and the second one is the
category of knowledge content $K_{\Phi\Theta}(f)$. These two
categories are related by the functor
$$
\Ct_f : L_{\Theta}(\Phi) \to K_{\Phi\Theta}(f).
$$
\endproclaim

This functor $\Ct_f$ transforms knowledge description to content
of knowledge.
We do not assume that between different $f_1$ and $f_2$ in $F$
there are any ties: instances are independent. On the other hand,
between some $f_1$ and $f_2$ there may be relations that we will
try to take into account.

A content of knowledge $\Ct_f(X,T)=(X,T^f)$ corresponds to an
object $(X,T)$ of the category $L_\Theta(\Phi)$, which is a
description of knowledge. We view the description $T$ as a {\it
query} to a knowledge base, and $A=T^f$ as a {\it reply to this
query}.

Besides, if there is a relation $s_\ast$ between $(X,T_1)$ and
$(Y,T_2)$, then  there will be a relation $\widetilde s
=\widetilde {s_\ast}$ between $(X,A)$ and $(Y,B)$, where $A=T_1^f,
\ B=T_2^f$  .

This peculiarity of the definition naturally reflects geometrical
essence of knowledge.

In fact, in this definition of a knowledge base the category of
knowledge is decomposed to two categories: the category of
description of knowledge and the category of content of knowledge,
tied by the functor of transition from description to content.



\head{ 9 Equivalence of  knowledge bases}
\endhead

\subhead{ 9.1 Definition}
\endsubhead

Let knowledge bases $KB_1=KB(G_1,\Phi_1,F_1)$  and
$KB_2=KB(G_2,\Phi_2,F_2)$ correspond to the given multimodels
$(G_1,\Phi_1,F_1)$ and $(G_2,\Phi_2,F_2)$.

\proclaim{Definition 1} Knowledge bases $KB_1$ and $KB_2$ are
called informationally  equivalent, if there exists a bijection
$\alpha: F_1\to F_2$ such that for every $f\in F_1$ there exist
homomorphisms
$$
\beta_f: \Hal_\Theta(\Phi_1)\to \Hal_\Theta(\Phi_2)
$$
$$
\beta_f': \Hal_\Theta(\Phi_2)\to \Hal_\Theta(\Phi_1)
$$
and an isomorphism of categories
$$
\widetilde\gamma_f: K_{\Phi_1\Theta}(f)\to
K_{\Phi_2\Theta}(f^\alpha)
$$
such that the commutative diagrams of functors of categories hold:
\endproclaim
$$
\CD
 L_\Theta(\Phi_1)@>\widetilde\beta_f>> L_\Theta(\Phi_2)\\
@V\Ct_f VV @VV \Ct_{f^\alpha}V\\
K_{\Phi_1\Theta}(f)@>\widetilde\gamma_f>> K_{\Phi_2\Theta}(f^\alpha)\\
\endCD
$$
and
$$
\CD
 L_\Theta(\Phi_1)@<\widetilde\beta'_f<< L_\Theta(\Phi_2)\\
@V\Ct_f VV @VV \Ct_{f^\alpha}V\\
K_{\Phi_1\Theta}(f)@<(\widetilde\gamma_f)^{-1}<< K_{\Phi_2\Theta}(f^\alpha)\\
\endCD
$$

Denote these diagrams by $\ast$ and $\ast\ast$, respectively.
Rewrite commutative diagrams for the object $(X, T)$ of the
category $L_\Theta (\Phi_1)$ in the form $(X, T^f)^{\widetilde
\gamma _f} = (X, T^{\beta_f f^\alpha})$ and for the object $(X,
T)$ of the category $L_\Theta (\Phi_2)$ in the form $(X,
T^{f^\alpha})^{\widetilde \gamma _f ^{-1}} = (X, T^{\beta' _f
f})$.

From this follows
$$
(X, T^f) = (X, T^{\beta_f
f^\alpha})^{\widetilde{(\gamma_f)}^{-1}},
$$
$$
(X, T^{f^\alpha}) = (X, T^{\beta'_f f})^{\widetilde \gamma_f },
$$

The last means that everything which can be known from $KB_1$ can
be also known from $KB_2$ and vice versa. Similar property holds
for morphisms, i.e. connections between objects. Equivalence of
knowledge bases we consider as a triple $(\alpha,\ast,\ast\ast)$,
where $\alpha: F_1\to F_2$ is a bijection, while $\ast$ and
$\ast$$\ast$ define the corresponding diagrams for every $f\in
F_1$.

The next proposition deals with the transition from  knowledge
bases to data\-bases. Let $R_f$ be the image of the homomorphism
$\Val_f: \Hal_\Theta(\Phi)\to \Hal_\Theta(G).$

\proclaim{Proposition 19} If a bijection $\alpha: F_1\to F_2$
determines equivalence of the bases $KB_1$ and $KB_2$ then for
every $f\in F_1$ we have an isomorphism of Halmos algebras
$\gamma_f:R_f\to R_{f^\alpha}$.
\endproclaim

Proof.

Proceed from the corresponding diagrams $\ast$ and $\ast\ast$.
Given  set $X$, take a set $T$ of one element
$u\in\Hal_{\Phi_1\Theta}(X)$. In this case $T^f=\Val_f(u)$. We
have $Ct_f(X,T)=(X,\Val_f(u))$,
$$
(X,\Val_f(u))^{\widetilde\gamma_f}=Ct_{f^\alpha}(X,u^\beta)=(X,(u^\beta)^{f^\alpha})=
(X,\Val_{f^\alpha}(u^\beta)).
$$

Hence, $\widetilde \gamma_f$ transfers $\Val_f(u)$ to
$\Val_{f^\alpha}(u^\beta)$ for every $u$, which means that
$\widetilde\gamma_f$ induces a mapping $\gamma_f: R_f\to
R_{f^\alpha}$. It is a homomorphism since $\Val_f$ and $\beta$ are
homomorphisms of algebras, and it is an injection since every
$R_f$ is a simple algebra \cite{Pl1}.

Let now $u_1$ be an arbitrary element of $\Hal_{\Phi_2\Theta}(X).$
Then the second diagram gives
$$
(X,\Val_{f^\alpha}(u_1))^{\widetilde\gamma_f^{-1}}=
(X,\Val_f(u_1^{\beta'_f})),
$$
and
$$
(X,\Val_{f^\alpha}(u_1))=(X,\Val_f(u_1^{\beta'_f}))^{\widetilde\gamma_f}=
(X,\Val_{f^\alpha}(u))^{\widetilde\gamma_f},
$$
where $u=u_1^{\beta'_f}$. This implies that $\gamma_f:R_f\to
R_{f^\alpha}$ is a surjection. Hence, we have an isomorphism
$\gamma_f:R_f\to R_{f^\alpha}$.

Proceed now from the isomorphism of Halmos algebras
$\gamma_f:R_f\to R_{f^\alpha}.$

\subhead{9.2 The case of finite models}
\endsubhead

First of all it is clear that for finite models $(G,\Phi, F)$ the
corresponding KB remains, in general, infinite.

We prove the following main

\proclaim{Theorem 18} Let the given models be finite. Then the
knowledge bases $KB_1$ and $KB_2$ are equivalent if and only if
there exists a bijection $\alpha: F_1\to F_2$ such that for every
$f\in F_1$ there is an isomorphism $\gamma_f: R_f\to
R_{f^\alpha}$.
\endproclaim

Proof.

In one direction the statement is always true. Let  now $\gamma_f:
R_f\to R_{f^\alpha}$ be an isomorphism for every $f\in F_1$.
According to the proposition 3  (see also \cite{PT}) we have the
homomorphisms $\beta_f:\Hal_{\Theta}(\Phi_1)\to
\Hal_{\Theta}(\Phi_2)$ and $\beta'_f:\Hal_{\Theta}(\Phi_2)\to
\Hal_{\Theta}(\Phi_1)$ such that the diagrams

$$
\CD
 \Hal_{\Theta}(\Phi_1)@>\beta_f>> \Hal_{\Theta}(\Phi_2)\\
@V\Val_f VV @VV\Val_{f^\alpha} V\\
R_f@>\gamma^f>> R_{f^\alpha}\\
\endCD
$$

$$
\CD
 \Hal_{\Theta}(\Phi_1)@<\beta'_f<< \Hal_{\Theta}(\Phi_2)\\
@V\Val_f VV @VV\Val_{f^\alpha} V\\
R_f@<\gamma_f^{-1}<< R_{f^\alpha}\\
\endCD
$$

\noindent
 are commutative.

 Simultaneously, there are functors
 $$
\widetilde\beta_f: L_{\Theta}(\Phi_1)\to L_{\Theta}(\Phi_2),
$$
$$
\widetilde\beta'_f: L_{\Theta}(\Phi_2)\to L_{\Theta}(\Phi_1).
$$

 It is left to define the isomorphism of categories $
\widetilde{\gamma_f}:K_{\Phi_1\Theta}(f)\to
K_{\Phi_2\Theta}(f^\alpha)$ such that the diagrams of the types
$\ast$ and $\ast\ast$ be commutative.

First we define $\widetilde\gamma_f$ on objects and then on
morphisms. Take an object (X,T) of the category $L_\Theta(\Phi_1)$
for an arbitrary object $(X,A)$ of the category
$K_{\Phi_1\Theta}(f)$ with $T^f=A$. We have
$Ct_f(X,T)=(X,T^f)=(X,A)$. Set
$$
(X,A)^{\widetilde\gamma_f}=(X,T^f)^{\widetilde\gamma_f}= (X,
\bigcap_{u\in T}\gamma_f\Val_f(u))= $$
$$(X, \bigcap_{u\in
T}\Val_{f^\alpha}(u^{\beta_f}))=(X,T^{\beta_ff^\alpha}).
$$

We want to show that this definition does not depend on the choice
of the set $T$ with $T^f=A$. Consider first the case when
$T_1^f=T_2^f=A$ and the sets $T_1$ and $T_2$ are finite. We have:
$(X,T_1^f)^{\widetilde\gamma_f}= (X,T_1^{\beta_ff^\alpha})$ and
$(X,T_2^f)^{\widetilde\gamma_f}= (X,T_2^{\beta_ff^\alpha}).$

We need to check that $T_1^{\beta_ff^\alpha}=
T_2^{\beta_ff^\alpha}. $ Indeed,
$$
T_1^{\beta_ff^\alpha}=\bigcap_{u_1\in
T_1}\Val_{f^\alpha}(\beta_fu_1))=\bigcap_{u_1\in
T_1}\gamma_f\Val_f(u_1)).
$$
Since $\gamma_f:R_f\to R_{f^\alpha}$ is an isomorphism of algebras
and $T_1$, $T_2$ are finite sets, we can rewrite the expression in
the form
$$
T_1^{\beta_ff^\alpha}=\gamma_f(\bigcap_{u_1\in
T_1}\Val_{f}(u_1))=\gamma_f(\bigcap_{u_2\in T_2}\Val_f(u_2))=
\bigcap_{u_2\in T_2}\gamma_f\Val_{f}(u_2))=T_2^{\beta_ff^\alpha}.
$$

Passing to the general case we proceed from finite models. Every
finite model is geometrically noetherian, i.e., if
$A=T_1^f=T_2^f$, then in $T_1$ and $T_2$ one can find finite
subsets $T_{01}$ and $T_{02}$ with $T{_{01}}^f=T_{02}^f=A$. Here,
$T_{01}^{\beta_ff_\alpha}=T_{02}^{\beta_ff_\alpha}$. We have to
verify that $ T_1^{\beta_ff^\alpha}=T_2^{\beta_ff^\alpha} $ and
$T_{01}^{\beta_ff_\alpha}=\bigcap_{u_1\in
T_1}\Val_{f^\alpha}(\beta_fu_1)$. We can take a finite subset
$T_{10}$ in $T_1$ such that $ T_1^{\beta_ff^\alpha}=
T_{10}^{\beta_ff^\alpha}$. Take the union of sets $T_{10}$ and
$T_{01}$ and denote it by $T_{001}$. Then $T_{001}^f=A=T_1^f$,
$T_1^{\beta_ff^\alpha}=T_{001}^{\beta_ff^\alpha}$. Analogously,
for $T_2$ take $T_{002}$ and $A=T_{001}^f=T_{002}^f$. Besides
that,
$$
T_1^{\beta_ff^\alpha}=T_{001}^{\beta_ff^\alpha}=
T_2^{\beta_ff^\alpha}.
$$
The equality $ T_1^{\beta_ff^\alpha}= T_2^{\beta_ff^\alpha}$ gives
commutativity of the diagram for objects.

Similarly, we build $\widetilde\gamma_f^{-1}$ having
$\gamma_f^{-1}$ and the equality
$\widetilde\gamma_f^{-1}=\widetilde{\gamma_f^{-1}}$ holds.

Now let us pass to morphisms. Remind first of all that to every
homomorphism $s:W(Y)\to W(X) $ there correspond
$$
s_\ast^1:\Hal_{\Phi_1\Theta}(Y) \to \Hal_{\Phi_1\Theta}(X)
$$
$$
s_\ast^2: \Hal_{\Phi_2\Theta}(Y) \to \Hal_{\Phi_2\Theta}(X).
$$

Let the objects $(Y,T_2)$ and $(X,T_1)$ be given in
$L_\Theta(\Phi_1)$. Recall that $s$ is admissible in respect to
$T_2$ and $T_1$ if $ s_\ast^1(u)\in T_1$ for every $u\in T_2$.
Here $s_\ast^1:(Y,T_2)\to (X,T_1)$ is a morphism. Proceed further
from an arbitrary homomorphism
$\beta:\Hal_\Theta(\Phi_1)\to\Hal_\Theta(\Phi_2)$. It had been
proved that if $s$ is admissible in respect to $T_2$ and $T_1$
then the same $s$ is admissible in respect to $T_2^\beta$ and
$T_1^\beta$ as well, i.e., $ s_\ast^1(u)\in T_1^\beta$ for every
$u\in T_2^\beta$. Hence, we have a morphism
$$
\widetilde\beta(s_\ast^1)=s_\ast^2:(Y,T_2^\beta)\to(X,T_1^\beta).
$$
Take now $\beta=\beta_f$ and apply $\Ct_{f^\alpha}$:
$$
\Ct_{f^\alpha}(s^2_\ast):\Ct_{f^\alpha}(X,T_1^{\beta_X})\to
\Ct_{f^\alpha}(Y, T_2^{\beta_X}).
$$
It can be rewritten as
$$
\Ct_{f^\alpha}(s_\ast^2):(X,T_1^{\beta_Xf^\alpha})\to (X,
T_2^{\beta_Yf^\alpha})
$$
or $$ \Ct_{f^\alpha}(s_\ast^2):(X,T_1^f)^{\gamma_f}\to (Y,
T_2^{f})^{\gamma_f}.
$$
Let now $T_1^f=A$, $T_2^f=B$. For $s_\ast^1:(Y,T_2)\to (X,T_1)$ we
have
$$
\Ct_{f}(s_\ast^1):(X,T_1^f)\to (Y, T_2^f)
$$
and  a related morphism
$$
\Ct_{f^\alpha}(s_\ast^2):(X,T_1^f)^{\gamma_f}\to (Y,
T_2^{f_2})^{\gamma_f}
$$

Commutativity of the diagram on morphisms means that
$$
\widetilde\gamma_f\Ct_f(s_\ast^1)=\Ct_{f^\alpha}(\widetilde\beta_f(s_\ast^1))
$$
for every $s^1_\ast: (Y,T_2)\to (X,T_1). $

Continuing consideration of finite models, proceed from the
isomorphism $\gamma_f:R_f\to R_{f^\alpha}$ and the corresponding
functor $\widetilde\gamma_f: K_{\Phi_1\Theta}(f)\to
K_{\Phi_2\Theta}(f^\alpha)$. This functor had been defined on the
objects, and now we are going to define it on morphisms.

Let $\tau: (X,A)\to (Y,B)$ be a morphism in $K_{\Phi_1\Theta}(f)$.
This $\tau$ appears as follows. A morphism
$$
s_\ast^1:\Hal_{\Phi_1\Theta}(Y) \to \Hal_{\Phi_1\Theta}(X)
$$
corresponds to $s: W(Y)\to W(X)$. If now $A=T_1^f$, $B=T_2^f$ and
$s^1_\ast$ is admissible in respect to $T_2$ and $T_1$ then we
have $\widetilde s_\ast^1: (X,A)\to (Y,B). $ We may say that
 $\tau=\widetilde s_\ast^1$ for some $s_\ast^1$.

 Define
 $$
 \widetilde\gamma_f(\widetilde s_\ast^1)
 =\widetilde s_\ast^2: (X,
 T_1^f)^{\widetilde\gamma_f}\to (Y,T_2^f)^{\widetilde\gamma_f}.
 $$
Here,
$$
(X, T_1^f)^{\widetilde\gamma_f}=(X,T_1^{\beta_ff^\alpha}),
 $$
  $$
 (Y, T_2^f)^{\widetilde\gamma_f}=(Y,T_2^{\beta_ff^\alpha}).
$$
do not depend on the choice of $T_1$ and $T_2$ with $T_1^f=A$ and
$T_2^f=B$. Check further that $\widetilde\gamma_f:
K_{\Phi_1\Theta}(f)\to K_{\Phi_2\Theta}(f^\alpha)$ determined in
such a way is in fact a functor and this functor provides
commutativity of the diagram on morphisms.

Note first of all that the definition of $\widetilde\gamma_f$ on
morphisms can be rewritten as
$$
\widetilde\gamma_f(\Ct_f(s_\ast^1))=(\Ct_{f^\alpha}(s_\ast^2)).
$$
Take two morphisms $\widetilde{s_{1\ast}^1}=\Ct_f(s_{1\ast}^1)$
and $\widetilde{s_{2\ast}^1}=\Ct_f(s_{2\ast}^1)$ and consider the
product
$$
\widetilde{s_{1\ast}^1}\widetilde{s_{2\ast}^1}=\Ct_f(s_{1\ast}^1)\Ct_f(s_{2\ast}^1)=
\Ct_f(s_{2\ast}^1s_{1\ast}^1)=\widetilde{s_{2\ast}^1s_{1\ast}^1}=
\widetilde{(s_2s_1)_\ast^1}.
$$

Apply $\widetilde\gamma_f$:
$$
\widetilde\gamma_f(\widetilde{(s_2s_1})_\ast^1)=(\widetilde{(s_2s_1)_\ast^2})=
\widetilde{s_{2\ast}^2s_{1\ast}^2}=\widetilde{s_{1\ast}^2}\widetilde{s_{2\ast}^2}=
\widetilde\gamma_f(\widetilde{s_{1\ast}^1})\widetilde\gamma_f(\widetilde{s_{2\ast}^1})
$$
Now check the commutativity of the diagram
$$
\CD
 L_{\Theta}(\Phi_1)@>\widetilde\beta_X>> L_{\Theta}(\Phi_2)\\
@V\Ct_f VV @VV\Ct_{f^\alpha} V\\
K_{\Phi_1\Theta}(f)@>\widetilde\gamma_f>> K_{\Phi_1\Theta}(f^\alpha)\\
\endCD
$$

Take a morphism $s_\ast^1:(Y,T_2)\to(X,T_1)$ in
$L_\Theta(\Phi_1)$. We have
$$
\widetilde{\beta_X}
(s_\ast^1):(Y,T_2^{\beta_X})\to(X,T_1^{\beta_X}),
$$
and
$$
\Ct_{f^\alpha}\widetilde{\beta_X}
(s_\ast^1):(X,T_1^{\beta_Xf^\alpha})\to(Y,T_2^{\beta_Xf^\alpha}),
$$
Rewrite it as
$$
\Ct_{f^\alpha}\widetilde{\beta_X}
(s_\ast^1):(X,T_1^f)^{\widetilde\gamma_f}\to(Y,T_2^{f})^{\widetilde\gamma_f},
$$
Further,
$$
\Ct_{f} (s_\ast^1):(X,T_1^f)\to(Y,T_2^{f}),
$$
$$
\widetilde\gamma_f\Ct_{f}
(s_\ast^1):(X,T_1^f)^{\widetilde\gamma_f}\to(Y,T_2^{f})^{\widetilde\gamma_f}.
$$
Check now the equality
$$
\widetilde\gamma_f\Ct_{f}
(s_\ast^1)=\Ct_{f^\alpha}\widetilde\beta_X (s_\ast^1)
$$
for every $ s_\ast^1$. We have
$$
\widetilde\gamma_f\Ct_{f} (s_\ast^1)=\widetilde\gamma_f(\widetilde
s_\ast^1)=\widetilde s_\ast^2,
$$
$$
\Ct_{f^\alpha}\widetilde\beta_X(s_\ast^1)=\Ct_{f^\alpha}(s_\ast^2)=
\widetilde s_\ast^2.
$$
This gives commutativity of the diagram $\ast$ of morphisms, i.e.,
$$
\widetilde\gamma_f\Ct_{f} =\Ct_{f^\alpha}\widetilde\beta_X.
$$
The same can be done for the functor
$\widetilde{\gamma_f^{-1}}=\widetilde\gamma_f^{-1}$ and the second
commutative diagram $\ast\ast$ that finishes the proof of the
theorem

\subhead{9.2 Main result. Additional remarks}
\endsubhead

7.1. Let us look at
 the definition of
equivalence from the general perspective of category theory. Given
two functors $\varphi_1 : C_1 \to C_1 ^0$ and $\varphi_2 : C_2 \to
C_2 ^0$, we say that $C_1$ and $C_2$ are equivalent in respect to
$\varphi_1$ and $\varphi_2$, if there is an isomorphism $\psi :
C_1 ^0 \to C_2 ^0$ and functors $\psi_1 : C_1 \to C_2$, $\psi_2 :
C_2 \to C_1$ with the commutative diagrams
$$
\CD
 C_1@>\psi_1>> C_2\\
@V\varphi_1 VV @VV\varphi_2 V\\
C_1^0@>\psi>> C_2^0\\
\endCD
$$

$$
\CD
 C_1@<\psi_2<< C_2\\
@V\varphi_1 VV @VV\varphi_2 V\\
C_1^0@<\psi^{-1}<< C_2^0\\
\endCD
$$

Usual equivalence of categories is equivalence in respect to the
transition to skeletons of categories. In our situation we may say
that equivalence of knowledge bases means that there exists
equivalence of categories of description of knowledge in respect
to transition to the categories of knowledge content.

7.2. Return to the definition of knowledge bases with multimodels
$(G_1, \Phi_1,$$ F_1)$ and $(G_2, \Phi_2, F_2)$, and let the
bijection $\alpha : F_1 \to F_2$ determine equivalence of the
corresponding $KB_1$ and $KB_2$. Assume that two instances $f_1$
and $f_2$ from $F_1$ are connected by a commutative diagram
$$
\CD
\Hal_{\Theta}(\Phi_1) @>\Val_{f_1}>> R_{f_1} \\
@. @/SE/ \Val_{f_2}// @VV\gamma V\\
@. R_{f_2}\\
\endCD
$$

\noindent
 where $\gamma$ is a homomorphism of algebras. We want to
estimate the relation between $f_1 ^\alpha$ and $f_2 ^\alpha$.

Proceed from the diagrams

$$
\CD
 \Hal_{\Phi_1\Theta}@>\beta_f>> \Hal_{\Phi_2\Theta}\\
@V\Val_f VV @VV\Val_{f^\alpha} V\\
R_f@>\gamma_f>> R_{f^\alpha}\\
\endCD
$$

$$
\CD
 \Hal_{\Phi_1\Theta}@<\beta' _f<< \Hal_{\Phi_2\Theta}\\
@V\Val_f VV @VV\Val_{f^\alpha} V\\
R_f@<\gamma_f^{-1}<< R_{f^\alpha}\\
\endCD
$$

$$
\CD
 R_{f_1}@>\gamma_{f_1}>> R_{f_1^\alpha}\\
@V\gamma VV @VV\gamma^\alpha V\\
R_{f_2}@>\gamma_{f_2}>> R_{f_2^\alpha}\\
\endCD
$$

Here,
$$
\gamma\Val_{f^1}=\Val_{f^2}, \qquad
\gamma^\alpha=\gamma_{f_2}\gamma\gamma_{f_1}^{-1}
$$
and
$$
\gamma^\alpha\Val_{f_1^\alpha}=\gamma^\alpha\gamma_{f_1}\Val_{f_1}\beta'_{f_1}=
\gamma_{f_2}\gamma\Val_{f_1}\beta'_{f_1}=\gamma_{f_2}\Val_{f_2}\beta'_{f_1}=
\Val_{f_2^\alpha}\beta_{f_2}\beta'_{f_1},
$$
Hence,
$\gamma^\alpha\Val_{{f^1}^\alpha}=\Val_{f_2^\alpha}\beta_{f_2}\beta'_{f_1},$
i.e., the connection is twisted by the product
$\beta_{f_2}\beta'_{f_1}.$

 At last, let us note that from the
diagrams above follow the natural identities:

1. $\Val_{f}(u)=\Val_{f}(\beta'_{f}\beta_{f}(u))$ for every
$u\in\Hal_\Theta(\Phi_1)$.

2. $\Val_{f^\alpha}(u)=\Val_{f^\alpha}(\beta_{f}\beta'_{f}(u)$ for
every $u\in\Hal_\Theta(\Phi_2)$.

7.3. Note that the equivalence condition of two knowledge bases in
the case of finite multimodels can be formulated in terms of these
multimodels (cf. \cite{PTP}).

\proclaim{Definition 2} Let the models $(G_1,\Phi_1,f_1)$ and
$(G_2,\Phi_2,f_2)$ be given. Let $\Aut(f_1)$ and  $\Aut(f_2)$ be
the corresponding groups of automorphisms. The models
$(G_1,\Phi_1,f_1)$  and $(G_2,\Phi_2,f_2)$ are called automorphic
equivalent if there exists an isomorphism of algebras $\delta:
G_1\to G_2$ such that
 $$
 \Aut(f_2)=\delta \Aut(f_1)\delta^{-1}.
 $$
\endproclaim
\proclaim{Definition 3} Let the multimodels $(G_1,\Phi_1,F_1)$ and
$(G_2,\Phi_2,F_2)$ be given. These multimodels are called
automorphic equivalent if there exists a bijection $\alpha :
F_1\to F_2$ such that for every $f\in F_1$ the models
$(G_1,\Phi_1,f)$  and $(G_2,\Phi_2,f^\alpha)$ are automorphic
equivalent.
\endproclaim

It is natural to define an isomorphism of multimodels with the
same set of relations $\Phi_1$ and $\Phi_2$. An isomorphism of
multimodels implies their automorphic equivalence. Evidently, the
inverse statement is not true.

Let the knowledge bases $KB_1=KB(G_1,\Phi_1,F_1)$ and
$KB_2=KB(G_2,\Phi_2,F_2)$ with the finite multimodels be given.

\proclaim{Theorem 19}The knowledge bases $KB_1=KB(G_1,\Phi_1,F_1)$
and $KB_2=KB(G_2,$ $\Phi_2,$ $ F_2)$ are informationally
equivalent if and only if the corresponding models are automorphic
equivalent.
\endproclaim

The proof of this theorem uses the Galois-Krasner theory and
follows from Theorems 14 and 18.


 Theorem 19 implies an algorithm
for the informational equivalence verification (see \cite{PK}).


 \bye


\frenchspacing \Refs\nofrills{\bf { Bibliography}}
\widestnumber\key{XXXXX}

\ref \key {\bf Ba} \by Bancillon F. \paper On the completeness of
query language for relational databases \jour Lecture Notes in CS
\yr 1978\vol 64 \pages 112--123
\endref


\ref \key {\bf Be} \by Beniaminov E.\paper Galois theory of
complete relational subalgebras of algebras of relations, logical
structures, symmetry  \jour NTI, ser. 2\yr 1980, \pages 237--261
\endref

\ref \key {\bf BJ} \by Bulatov A., Jeavons, P., \paper An
algebraic approach to multi-sorted constraints\jour Proceedings of
CP'03, to appear\yr 2004\pages 15pp
\endref




\ref \key {\bf G} \by  Ganter B., Mineau G. \book  Ontology,
metadata, and semiotics \publ  Lecture Notes in AI,
Springer-Verlag \vol 1867 \yr 2000 \pages 55--81
\endref

\ref \key {\bf H} \by Halmos P.R. \book Algebraic logic \yr 1969
\publ New York
\endref


\ref \key {\bf HMT} \by Henkin L., Monk J. D., Tarski A. \book
Cylindric Algebras \yr 1985 \publ North-Holland Publ. Co.
\endref

\ref \key {\bf JCP} \by Jeavons, P., Cohen, D.,; Pearson, J.
\paper Constraints and universal algebra \jour Ann. Math.
Artificial Intelligence \vol 24\year 1999\pages 51--67
\endref


\ref \key {\bf Kr} \by Krasner M. \paper  Generalisation et
anlogues de la theorie de Galois\jour Comptes Rendus de Congress
de la Victorie de l'Ass. Franc. pour l'Avancem. Sci.  \yr 1945,
\pages 54--58
\endref


\ref \key {\bf L} \by  Lenat D. \book Steps to Sharing Knowledge
\publ Toward Very Large Knowledge Bases, edited by N.J.I. Mars.
IOS Press \yr 1995
\endref

\ref \key {\bf ML} \by MacLane S. \book Categories for the working
mathematicians \yr 1971 \publ Springer
\endref


\ref \key {\bf NP} \by Nikolova D., Plotkin B. \paper  Some Notes
on Universal Algebraic Geometry \jour Proc. of Kuros Algebraic
Conference, Moscow 1998, Walter de Gryiter, \yr 1999, to appear
\endref


\ref \key {\bf Pl1} \by Plotkin B.I. \book Universal algebra,
algebraic logic and databases \yr 1993 \publ Kluwer
\endref






\ref \key {\bf Pl4} \by Plotkin B.I. \book  Algebra, categories
and databases \yr 1999 \publ Handbook of algebra, v.2, Elsevier,
Springer \pages 81-148
 \endref


\ref \key {\bf Pl5} \by Plotkin B.I.  Seven Lectures in Universal
Algebraic Geometry, Revised, Preprint, \book Arxiv: math,
GM/0204245, \year 2002, 87pp.
 \endref

\ref \key {\bf PT} \by Plotkin T. \book   Relational databases
equivalence problem \yr 1996 \publ Advances of databases and
information systems, Springer \pages 391-404
 \endref

\ref \key {\bf PTP} \by Plotkin B.I., Plotkin T. \paper
Geometrical aspect of databases and knowledge bases\jour Algebra
Universalis \vol 46 \yr 2001 \publ Birkhauser Verlag, Basel \pages
131-161
 \endref


\ref \key {\bf PK} \by Plotkin T., Knjazhansky M. \book
Informational Equivalence of Knowledge Bases Verification \yr 2004
\publ Scientific Israel, to appear
\endref


\ref \key {\bf S} \by   Sowa J. \book  Knowledge Representation:
Logical, Philosophical, and Computational Foundations \publ Brooks
Cole Publishing Co., Pacific Grove, CA \yr 2000
\endref
\bye